\documentclass[preprint, 3p]{elsarticle}

\pdfoutput = 1
\usepackage{lipsum}
\makeatletter
\def\ps@pprintTitle{%
	\let\@oddhead\@empty
	\let\@evenhead\@empty
	\def\@oddfoot{}%
	\let\@evenfoot\@oddfoot}
\makeatother

\usepackage{amssymb}
\usepackage{latexsym}

\usepackage{url}
\usepackage{xcolor}

\usepackage{hyperref}

\definecolor{newcolor}{rgb}{.8,.349,.1}
\journal{Medical Image Analysis}

\usepackage{graphicx}
\usepackage{tabularx}
\usepackage{longtable}
\usepackage{nicefrac}
\usepackage{placeins}
\usepackage{arydshln}

\usepackage{amsmath}
\usepackage{multirow}
\usepackage{tikz}
\usepackage{floatrow}
\usetikzlibrary{matrix,chains,positioning,decorations.pathreplacing,arrows}
\usepackage{amsfonts}
\usepackage{makeidx}  
\usepackage{amsmath, mathtools}
\usepackage{graphicx}
\usepackage{enumerate}
\usepackage{array}
\usepackage{caption}
\usepackage{subcaption}

\usepackage{stfloats} 
\usepackage{booktabs}
\usepackage{verbatim}
\usepackage{algorithm2e}
\graphicspath{{figures/}}
\usepackage[acronym,nomain]{glossaries}
\usepackage{float}
\usepackage{xcolor}
\usepackage{tabularx}
\usepackage{comment}

\usepackage{multicol}

\usepackage{tikz}
\usetikzlibrary{calc}
\usetikzlibrary{spy}
\usetikzlibrary{decorations.pathreplacing}
\usetikzlibrary{shapes}
\usetikzlibrary{through}
\usetikzlibrary{intersections}
\usetikzlibrary{decorations.markings}
\usetikzlibrary{3d}
\usetikzlibrary{positioning}

\usepackage{fancyhdr}
\usepackage{datetime}

\usepackage{quoting,xparse}
\usetikzlibrary{arrows, shapes, trees}
\usepackage{geometry}
\usepackage{hyperref}
\usepackage{ulem}

\usepackage{textcomp}
\usepackage{diagbox}

\usepackage{siunitx}
\usepackage[draft]{todonotes}   

\newcommand{\R}[0]{\ensuremath{\mathbb{R}}}

\DeclareCaptionLabelFormat{andtable}{#1~#2  \&  \tablename~\thetable}

\newcommand{\tabref}[1]{Tab.~\ref{#1}}

\newcommand\clearrow{\global\let\rowmac\relax}

\newcolumntype{?}{!{\vrule width 1pt}}

\usepackage[symbol]{footmisc}
\def\correspondingauthor{\footnote{Corresponding author: gernot.plank@medunigraz.at, Medical University of Graz, Neue
        Stiftingtalstraße 6(MC1.D.)/IV,
        8010, Graz, Austria}}

\NewDocumentCommand{\bywhom}{m}{
    {\nobreak\hfill\penalty50\hskip1em\null\nobreak
        \hfill\mbox{\normalfont(#1)}%
        \parfillskip=0pt \finalhyphendemerits=0 \par}%
}

\NewDocumentEnvironment{pquotation}{m}
{\begin{quoting}[
        indentfirst=true,
        leftmargin=\parindent,
        rightmargin=\parindent]
    \bywhom{#1}
\end{quoting}}

\clearrow

\SetCommentSty{mycommfont}

\makeglossaries

\definecolor{ColRev1}{RGB}{192,0,0}
\definecolor{ColRev2}{RGB}{0,96,160}
\definecolor{ColRev3}{RGB}{229,122,0}
\definecolor{ColRev4}{RGB}{230,0,153}
\definecolor{ColChange}{RGB}{0,160,0}
\definecolor{forestgreen}{RGB}{20,160,80}

\usepackage[acronym]{glossaries}

\newacronym{lv}{LV}{left ventricle}
\newacronym{rv}{RV}{right ventricle}
\newacronym{la}{LA}{left atrium}
\newacronym{ra}{RA}{right atrium}
\newacronym{cs}{CS}{coronary sinus}
\newacronym{lpv}{LPV}{left pulmonary veins}
\newacronym{lipv}{LIPV}{left inferior pulmonary vein}
\newacronym{lspv}{LSPV}{left superior pulmonary vein}
\newacronym{rpv}{RPV}{right pulmonary veins}
\newacronym{ripv}{RIPV}{right inferior pulmonary vein}
\newacronym{rspv}{RSPV}{right superior pulmonary vein}
\newacronym{laa}{LAA}{left atrial appendage}
\newacronym{raa}{RAA}{right atrial appendage}
\newacronym{svc}{SVC}{superior vena cava}
\newacronym{ivc}{IVC}{inferior vena cava}
\newacronym{tv}{TV}{tricuspid valve}
\newacronym{mv}{MV}{mitral valve}

\newacronym{bb}{BB}{Bachmann's bundle}
\newacronym{fo}{FO}{fossa ovalis}
\newacronym{pm}{PM}{pectinate muscle}
\newacronym{san}{SAN}{sino-atrial node}
\newacronym{hps}{HPS}{His-Purkinje System}
\newacronym{ic}{IC}{Inter-atrial Conduction}

\newacronym{af}{AF}{atrial fibrillation}

\newacronym{ep}{EP}{electrophysiology}
\newacronym{ecg}{ECG}{electro-cardiogram}
\newacronym{egm}{EGM}{electrogram}
\newacronym{eam}{EAM}{electro-anatomical mapping}
\newacronym{pwd}{PWD}{P-Wave Duration}

\newacronym{uac}{UAC}{universal atrial coordinates}
\newacronym{rd}{R-D}{Reaction-Diffusion}
\newacronym{re}{R-E}{Reaction-Eikonal}
\newacronym{bem}{BEM}{Boundary Element Method}
\newacronym{fem}{FEM}{Finite Element Method}
\newacronym{relf}{RELF}{Reaction-Eikonal Lead Field}
\newacronym{pca}{PCA}{Principal Component Analysis}
\newacronym{dbc}{DBC}{Dirichlet Boundary Conditions}

\newacronym{scn}{SCN}{SpatialConfiguration-Net}
\newacronym{ct}{CT}{computed tomography}
\newacronym{mri}{MRI}{Magnetic Resonance Imaging}

\newacronym{rmse}{RMSE}{relative mean square error}
\newacronym{ld}{LD}{Laplace-Dirichlet}

\begin{document}
	
	\begin{frontmatter}
		
		\title{An efficient end-to-end computational framework for the generation of ECG calibrated volumetric models 
			of human atrial electrophysiology}
		\author {Elena Zappon$^{a,b}$}
		\author {Luca Azzolin$^{c}$}
		\author {Matthias A.F. Gsell$^{a}$}
		\author {Franz Thaler$^{d,a,e}$}
		\author {Anton J. Prassl$^{a}$}
		\author {Robert Arnold$^{a}$}
		\author {Karli Gillette$^{f,g,a,b}$}
		\author {Mohammadreza Kariman$^{a}$}
		\author {Martin Manninger-W\"{u}nscher$^{h}$}
		\author {Daniel Scherr$^{h}$}
		\author {Aurel Neic$^{c}$}
		\author {Martin Urschler $^{e,b}$}
		\author {Christoph M. Augustin$^{a,b}$}
		\author {Edward J. Vigmond$^{i}$}
		\author {Gernot Plank$^{a,b}$\correspondingauthor{} \\ \ \\
			$^a$ Division of Biophysics, Medical University of Graz, Graz, Austria \\
			$^b$ BioTechMed-Graz, Graz, Austria \\
			$^c$ NumeriCor Gmbh, Graz, Austria\\
			$^d$ Institute of Visual Computing, Graz University of Technology, Graz, Austria \\
			$^e$ Institute for Medical Informatics, Statistics and Documentation, Medical University of Graz, Graz, Austria \\
			$^f$ Scientific Computing and Imaging Institute, University of Utah, USA\\
			$^g$ Department of Biomedical Engineering, University of Utah, USA\\
			$^h$ Clinical Department of Cardiology, Medical University of Graz, Graz, Austria\\
			$^i$ University of Bordeaux, CNRS, Bordeaux, France\\}
		
		\begin{abstract}
			Computational models of atrial \gls{ep} are increasingly utilized for applications such as the development of advanced mapping systems, personalized clinical therapy planning, and the generation of virtual cohorts and digital twins. These models have the potential to establish robust causal links between simulated \textit{in silico} behaviors and observed human atrial EP, enabling safer, cost-effective, and comprehensive exploration of atrial dynamics. However, current state-of-the-art approaches lack the fidelity and scalability required for regulatory-grade applications, particularly in creating high-quality virtual cohorts or patient-specific digital twins. Challenges include anatomically accurate model generation, calibration to sparse and uncertain clinical data, and computational efficiency within a streamlined workflow.
			
			This study addresses these limitations by introducing novel methodologies integrated into an automated end-to-end workflow for generating high-fidelity digital twin snapshots and virtual cohorts of atrial EP. These innovations include: (i) automated multi-scale generation of volumetric biatrial models with detailed anatomical structures and fiber architecture; (ii) a robust method for defining space-varying atrial parameter fields; (iii) a parametric approach for modeling inter-atrial conduction pathways; and (iv) an efficient forward EP model for high-fidelity \gls{ecg} computation.
			
			We evaluated this workflow on a cohort of 50 atrial fibrillation patients, producing high-quality meshes suitable for reaction-eikonal and reaction-diffusion models and demonstrating the ability to simulate atrial ECGs under parametrically controlled conditions. These advancements represent a critical step toward scalable, precise, and clinically applicable digital twin models and virtual cohorts, enabling enhanced patient-specific predictions and therapeutic planning.
		\end{abstract}
		
		
	\end{frontmatter}
	
	

	\section{Introduction}
	\label{sec:intro}
	Computational models of atrial \gls{ep} are increasingly being considered in a variety of applications, ranging from the industrial development of devices such as \gls{eam} systems \citep{reddy2023electrographic,szili2023electrographic}
	to the stratification and planning of clinical therapy \citep{boyle2019computationally}.
	These applications are built upon the mechanistic nature of biophysical models of atrial \gls{ep}, 
	and are based on the tacit assumption 
	that simulated \textit{in silico} behaviors closely correspond to the real human atrial \gls{ep} observed in patients.
	If sufficient regulatory-strength evidence of such a close causal relation is provided,
	\emph{in silico} models can be used to explore human atrial \gls{ep} more comprehensively, 
	in a safer and more cost-effective manner compared to the current paradigm based on preclinical animal testing and clinical trials \citep{dossel2012computational,schotten2011pathophysiological}. 
	However, current state-of-the-art models supporting evidence of such a tight causal and quantitatively accurate relation are limited. 
	The vast majority of computational studies use simplified anatomical models \citep{azzolin2023:_augmenta,labarthe2014:_bilayer,roney2023:_bia_vol} with uncalibrated default parameters \citep{boyle2019computationally, sakata2024:_assessing},
	and any comparison to directly observable quantities, such as \gls{ecg}, is usually limited \citep{ferrer2015detailed,loewe2016influence} or even not accounted for \citep{roney2023:_bia_vol}.
	This can be largely attributed to the limited capabilities of current modeling technologies in terms of anatomical model generation as well as calibration and simulation technologies.
	A direct comparison to observations would reveal major discrepancies between physical and virtual spaces,
	thus undermining the credibility of the model.
	
	The degree of fidelity needed, as well as the metrics used for measuring it, is application-dependent \citep{bhagirath2024:_bits}.
	Beyond the minimum requirement of a mechanistic relation -- that is, the models of atrial \gls{ep} are able to qualitatively recapitulate all the mechanisms of atrial \gls{ep} at play --
	in silico models can be calibrated to be representative of a group of patients of interest, 
	also covering anatomical and functional variability \citep{niederer2020creation}.
	Such representation of a patient cohort,
	rather than an individual, may facilitate generic interpretation of observations in real physical atria, and prediction of \gls{ep} responses to therapeutic interventions.
	Such sets of functionally-similar models referred to as virtual cohorts, have recently started to be considered for in silico trials and for safety and efficacy testing of new devices or therapies \citep{viceconti2020:_in_silico}.
	Nevertheless, creating virtual cohorts of in silico models will necessitate a shift in cardiac modeling, moving from a limited set of custom models to efficient and scalable workflows capable of generating large volumes of models quickly \citep{niederer2020creation}.
	
	Most demanding are clinical applications geared toward precision cardiology, that is, to tailor therapies to individual patients. 
	There, models are sought that replicate cardiac anatomy and structures and quantitatively calibrate to match functional observations from an individual patient's heart in a one-to-one manner.
	Such functionally equivalent models, where a particular stimulus or perturbation leads to the same emergent response in virtual and real space at a single time point, are referred to as digital twin snapshots. When the models are continuously or periodically updated with measurements, they become true digital twins \citep{bhagirath2024:_bits,corral2020digital,hopman2023right}.
	
	However, 
	the ability of current atrial \gls{ep} modeling pipelines 
	for creating high-fidelity digital twin snapshots or virtual cohorts at scale 
	in a sufficiently efficient and robust manner 
	is severely limited. 
	The key challenge is to describe the electrical sources in the atria accurately enough to predict the electrical potential field in its surroundings where all observable measurements are recorded. 
	In humans \emph{in vivo}, electrical measurements are in the form of \glspl{egm} by devices or \gls{eam} systems, or as \gls{ecg} at the body surface.
	While conceptually simple, the implementation of such a pipeline is vastly demanding, posing a long list of formidable challenges many of which remain unaddressed.
	In general, modeling pipelines for creating cardiac digital twin models of atrial \gls{ep} 
	are separated into two distinct stages, an anatomical and a functional twinning stage \citep{azzolin2023:_augmenta,labarthe2014:_bilayer,roney2023:_bia_vol}.
	
	At the anatomical stage,
	multi-label segmentation of tomographic images is performed to identify all relevant domains \citep{payer2017multi}
	which are turned into multi-label computational meshes \citep{crozier2016image,prassl09:_tarantula}
	to accurately represent biatrial anatomy.
	Involved procedures are notoriously laborious, 
	requiring numerous manual operator interventions by trained experts and significant computational resources \citep{azzolin2023:_augmenta,roney2023:_bia_vol}, 
	to obtain anatomically accurate representations of sufficient mesh quality for a given type of cardiac EP simulation.
	Specifically for atrial anatomies, these are often simplified \citep{harrild2000:_atria} 
	and represented as manifolds only \citep{azzolin2023:_augmenta,labarthe2014:_bilayer,roney2023:_bia_vol},
	thus limiting the achievable quantitative accuracy in representing electrical sources and associated potential fields.  
	
	At the functional twinning stage, 
	the fundamental core challenge of calibrating an atrial \gls{ep} model to clinical data
	is to infer high-dimensional space-varying parameter fields governing the \gls{ep} behavior 
	from limited sparse clinical recordings
	that are afflicted with substantial observational uncertainties \citep{coveney2022calibrating,whittaker2020calibration}.
	This poses a number of technological problems.
	Model functionalization refers to conceiving a framework for comprehensively describing 
	a sufficiently high-dimensional parameter space 
	defined over geometrically complex objects, such as the heart \citep{bayer2018universal,roney2019uac}, 
	that encapsulates all relevant factors governing atrial \gls{ep}, and the genesis of the associated extracellular potential field.
	These fields must be exposed to unattended algorithmic manipulation to facilitate parameter sweeps in order to minimize the mismatch between simulated and observed data within an optimization procedure \citep{grandits2021:_geasi,grandits2023:_geodesic_bp}.
	A computationally efficient yet accurate forward \gls{ep} model for generating the observed electrical recordings, 
	i.e. \glspl{ecg}, \glspl{egm} or \glspl{eam},  is required
	to cope with the computational burden of a large number of model evaluations incurring during optimization \citep{gillette2021:_framework,neic17:_reaction_eikonal,pezzuto2017:_ecg}.
	Computational workflows meeting these criteria and reported so far in the literature \citep{azzolin2023:_augmenta,boyle2019computationally,roney2023:_bia_vol}, are primarily limited in terms of anatomical and structural representation, functional calibration capabilities, as well as automation, efficiency, robustness, reproducibility, and numerical accuracy.
	
	In this study, we address these limitations by developing novel methodologies that are essential for the scalable generation of atrial digital twin snapshots and virtual cohorts, calibrated by \glspl{ecg} and, potentially, \glspl{egm} or \glspl{eam}. 
	These comprise: 
	\begin{enumerate}[i)]
		\item An automated flexible multi-scale approach for generating anatomically accurate volumetric biatrial models with comprehensive parametric incorporation of atrial structures and fiber architecture, and
		of sufficiently high mesh quality to be suitable for most widely used forward \gls{ep} models;
		\item A robust method for generating a volumetric atrial reference frame 
		for defining space-varying atrial parameter fields, and their unattended manipulation in optimization sweeps;
		\item A novel method for the flexible parametric representation of \gls{ic} pathways;
		\item  An efficient and clinically-compatible \gls{ecg} forward model 
		that computes high fidelity gold standard \glspl{ecg} and \glspl{egm} with close to real-time performance. 
	\end{enumerate}
	
	We report on the integration of these methodological developments into a highly automated and efficient end-to-end workflow 
	for building high-fidelity digital twin snapshots of atrial \gls{ep} from clinical data. 
	We elucidate in detail the methodological underpinnings of all processing stages
	-- from image analysis over model generation to \gls{ecg} and \gls{egm} prediction --
	and evaluate the performance of the end-to-end workflow on a clinical cohort of patients
	treated for \gls{af}, by processing 50 patient datasets randomly selected from the Graz \gls{af} registry.
	For each patient, both coarse and fine-resolution meshes were generated 
	at $\approx$ \SI{0.90}{\milli\meter} and \SI{0.25}{\milli  \meter} for reaction-eikonal \citep{neic17:_reaction_eikonal}
	and reaction-diffusion mono- or bidomain models, respectively, together with a set of automatically generated volumetric \gls{uac}.
	We demonstrate the ability to compute atrial \glspl{ecg} by performing sweeps on parametrically-controlled and easy-to-generate inter-atrial pathways, and we show the possibility of automatically varying the position and shape of the \gls{san}, allowing fast calibration of the model at hands.
	
	The structure of the paper is as follows: Section \ref{Sec:methods} presents a comprehensive overview of the proposed end-to-end workflow, encompassing the anatomical twinning stage (Section \ref{Subsec:anatomical_twinning}) and the functional modeling stage (Section \ref{Subsection:functional_modeling_stage}). Specifically, our approach for generating anatomically accurate volumetric biatrial models, including anatomical structures and fibers, is detailed in Sections \ref{Subsubsection:Imaging}-\ref{Subsubsec:anatomical_structure}. Sections \ref{Subsubsec:UAC} and \ref{Subsubsec:Torso} describe the computation of the \gls{uac} and torso model surrounding the biatrial mesh. Our novel method for computing interatrial connections is outlined in Section \ref{Subsubsection:cables}, while the parameter sweeping process for model calibration is detailed in Section \ref{Subsubsec:baseline_sim}. The numerical results are provided in Section \ref{Sec:results}, with the discussion, limitations, and conclusions presented in Sections \ref{Sec:discussion}, \ref{Sec:limitation}, and \ref{Sec:conclusion}, respectively.

	
	\section{Methods} 
	\label{Sec:methods}
\begin{figure}[!t]
		\centering
		\includegraphics[width=0.9\textwidth]{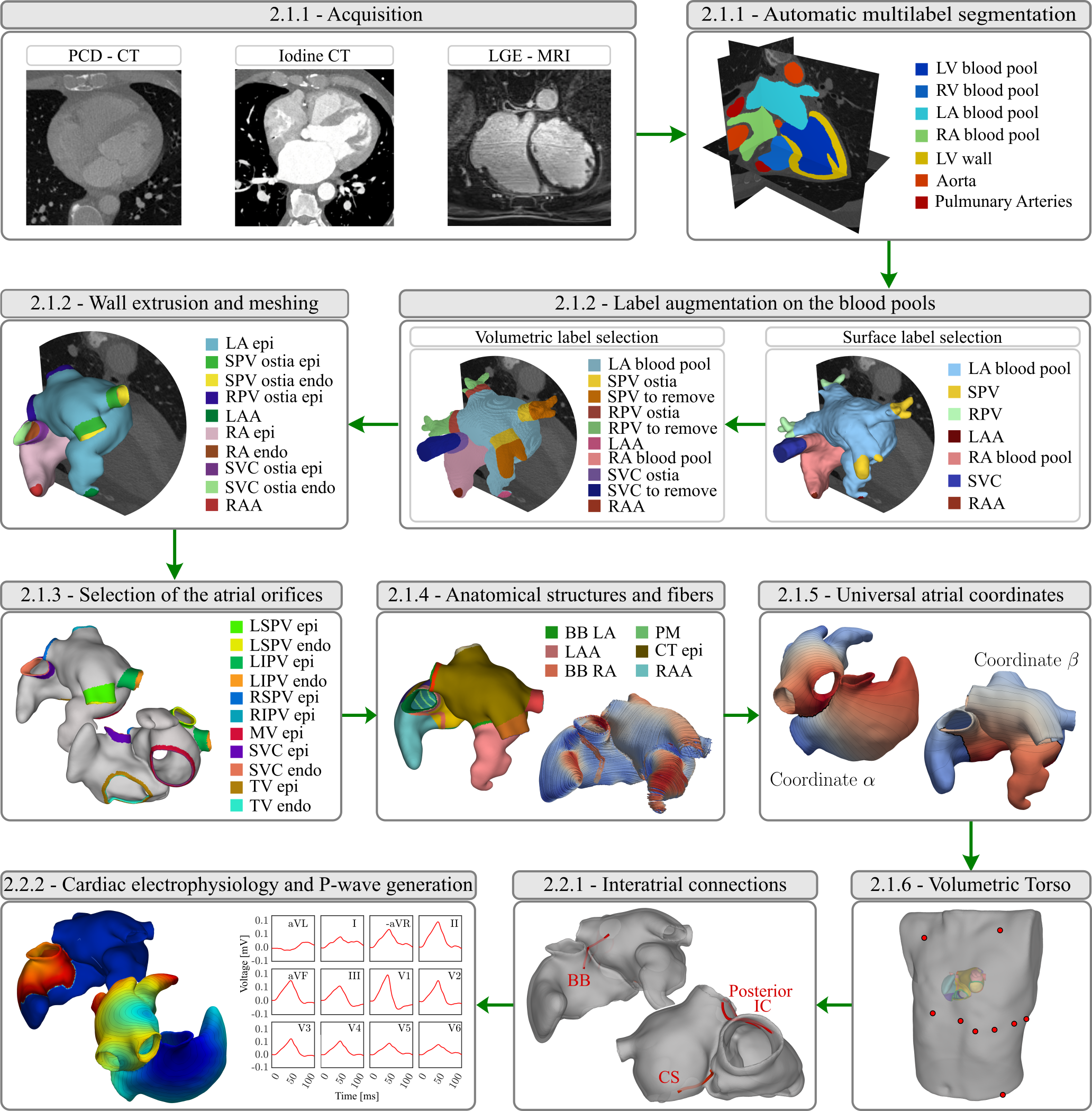}
		\caption{Schematic outline of the end-to-end framework for the generation of \gls{ecg}-calibrated volumetric models of patient-specific human atria. After image acquisition, the workflow comprises nine steps: automatic multilabel segmentation, automatic label augmentation on the blood pools, extrusion of the volumetric bilayer walls, automatic selection of atrial orifices, automatic annotation of anatomical structures and fiber generation, generation of \gls{uac}, registration and generation of a torso volume conductor, generation of \gls{ic}s, and cardiac electrophysiological simulation and P-wave generation. We moreover provide the paper section index where each step is detailed.}
		\label{Fig:scheme}
	\end{figure}

	We provide a schematic overview of an end-to-end workflow for the generation of anatomically accurate volumetric biatrial-torso model
	in Figure \ref{Fig:scheme}. 
	
	The workflow comprises two major processing stages, (i) an anatomical twinning stage comprising the generation of the atrial anatomical model with structural labels and fiber architecture, the generation of an atrial reference frame for spatial parameter swiping, and a torso anatomical model,
	and (ii) a functional modeling stage, for 
	the definition of inter-atrial pathways, and
	the setting up of a forward \gls{ep} model for representing and calibrating atrial action potential propagation, 
	and for generating the associated clinically observable signals, such as the \glspl{ecg} and the \glspl{egm}.
	Both stages are closely linked to the choice of forward \gls{ep} mathematical model, 
	which influences, for instance, the targeted spatial mesh resolution of the biatrial anatomical model.
	
	The anatomical twinning stage is organized following an almost automatic and sequential step approach:
	\begin{enumerate}[i)]
		\vspace{-0.2cm}\item image segmentation: by exploiting a convolutional neural network, tomographic scans, either \gls{ct} or \gls{mri}, are automatically segmented;
		\vspace{-0.2cm}\item label augmentation on the atrial blood pool: 
		veins and appendages on the atrial blood pools are identified and labeled, 
		for defining the border between atrial myocardium and attached vasculature;
		\vspace{-0.2cm} \item extrusion of the volumetric atrial walls: 
		by accounting for the identified labels on the blood pools, a sequence of prescribed operations is employed on the segmentation to extrude the endo- and epicardial walls endowed with a selected thickness, automatically opening veins and valves, and generating a corresponding smooth volumetric mesh;
		\vspace{-0.2cm}\item selection of atrial orifices: a fully automated process is carried out to identify the complete set of vein ostia and valves on the atrial walls, separated in endocardial and epicardial tissue;    
		\vspace{-0.2cm}\item automatic labeling of anatomical structures and fiber generation:
		the remaining anatomical regions are automatically annotated, and the atrial fiber architecture is computed;
		\vspace{-0.2cm}\item \glspl{uac} generation: a system of normalized volumetric universal coordinates is computed;
		\vspace{-0.2cm}\item integration of a torso volumetric conductor: a reference torso anatomical model is automatically integrated with the biatrial anatomical model 
		and lead fields are computed for all considered electrode positions
		for computing all clinical \gls{ecg} signals of interest.
	\end{enumerate}
	
	The functional modeling stage comprises:
	\begin{enumerate}[i)]
		\setcounter{enumi}{7}
		\vspace{-0.2cm}\item generating \glspl{ic}: a flexible set of conduction pathways that cannot be delineated from images are computed to facilitate inter-atrial impulse conduction 
		\vspace{-0.2cm}\item generate \gls{ep} clinical outputs: employing either \gls{rd} monodomain formulation, or a \gls{re} model to represent the transmembrane voltage, $V_m$, in the atrial model, a parameter spaces of \gls{ep} parameter is selected to be used for the model calibration. The P-wave of the 12-lead \glspl{ecg} on the torso are computed using the respective lead fields  \citep{gillette2022personalized,potse2018scalable}. The same approach can be potentially employed to compute \glspl{egm} in the blood pools.
	\end{enumerate}
	
	
	
	The entire workflow is implemented in a single user-friendly software 
	building on the package \texttt{meshtool} \citep{neic2020automating} and on \texttt{python} code. 
	For meshing operations the software TetGen (Weierstrass Institute, Berlin, Germany) 
	and NetGen (CerbSim GmbH, Vienna, Austria) have been licensed 
	and integrated with \texttt{meshtool}. 
	An interactive mode is moreover available, allowing for the verification of the labeling process, the visual control of all processing stages, and the correction of any potential software errors, through the use of the software
	\href{https://numericor.at/rlb/wordpress/products/}{CARPentry Studio} (Numericor GmbH, Graz, Austria). 
	All required interactive steps are readily supported by the freely available starter version of \href{https://numericor.at/rlb/wordpress/resources/}{CARPentry Studio for academia}.
	Each step of our workflow will be extensively described in the following sections.
	
	\subsection{Anatomical twinning stage}
	\label{Subsec:anatomical_twinning}
	\subsubsection{Imaging data acquisition and segmentation}
	\label{Subsubsection:Imaging}
	For developing and testing the workflow, 50 patients 
	diagnosed with \gls{af} and scheduled for \gls{af} ablation therapy were selected.
	Iodanized contrast \gls{ct} scans were acquired at an isotropic resolution of \SI{0.4}{\milli \meter} 
	as part of the routine standard of care at the Medical University of Graz hospital 
	for patients included in the local \gls{af} ablation registry. 
	This registry was approved by the ethics committee of the Medical University of Graz (reference number 26-217 ex 13/14) 
	and all patients gave written informed consent. 
	
	Each \gls{ct} dataset was segmented using the \gls{scn} \citep{thaler2021efficient}, 
	an automated multi-label segmentation method based on convolutional neural networks. 
	Prior to input to the \gls{scn}, 3D image stacks were cropped around the center of the heart, 
	re-sampled to an isotropic resolution of \SI{1.2}{\milli \meter}, 
	and Gaussian filtered with $\sigma=$\SI{1}{\milli \meter}.
	The \gls{scn} was trained to recognize and label seven cardiac domains, 
	comprising all four major cardiac \gls{lv}, \gls{rv}, \gls{la}, \gls{ra} blood pools, 
	the \gls{lv} myocardium, and the vascular blood pools in the aorta and the pulmonary artery.
	Segmentation quality was assessed interactively, to ascertain topological soundness. 
	Topological errors, if present, such as connections between the blood pools of the \gls{la} appendage and the \gls{lspv}, 
	or between the \gls{cs} and the \gls{la} were identified, 
	and manually corrected using the \texttt{ITK-SNAP} software \citep{py06nimg}. 
	
	
	\subsubsection{Automated label augmentation on the blood pools}
	\label{Subsubsection:LabelAugmentation}
	Label augmentation is performed to identify the atrial appendages and veins, 
	which will later be used in the wall extrusion process. 
	Among all the volumetric labels defined through segmentation, 
	only the \gls{la} and \gls{ra} blood pools are targeted at this stage, while all other anatomical labels serve as auxiliary references.
	
	The blood pool label augmentation is carried out in two stages: 
	the first stage automatically defines key landmarks on the surface of the blood pools; 
	the second stage effectively labels the entire blood pool volumes. 
	
\begin{figure}[!t]
		\centering
		\includegraphics[width=0.9\textwidth]{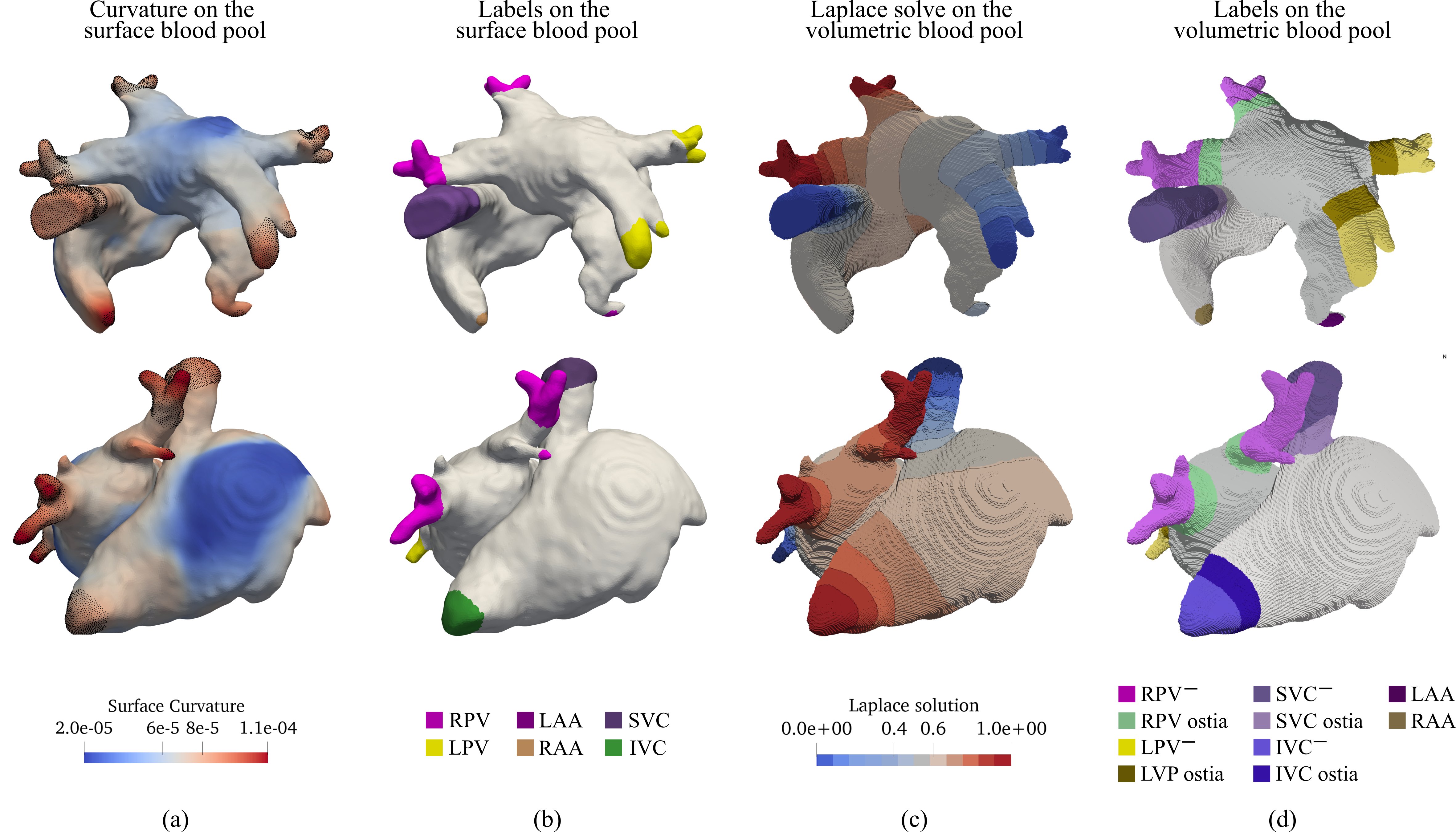}
		\caption{(a) Point-wise curvature on the blood pool surface mesh. The figure also highlights the correlation between the surface region of high curvature and anatomical structures to mark. 
			(b) Minimal set of labels selected based on surface curvature, including \gls{lpv}, \gls{rpv} and \gls{laa} on the \gls{la}, 
			and \gls{ivc}, \gls{svc} and \gls{raa} on the \gls{ra}. 
			(c) Instance of \gls{ld} solution and isosurfaces on the blood pool volumetric mesh. 
			(d) Vein ostiae and the discardable tissue labeled on the blood pool volumetric mesh.}
		\label{Fig:image_process}
	\end{figure}

	In the first stage, beginning with the segmentation, triangular surface meshes 
	of the \gls{la} and \gls{ra} blood pools are generated. 
	Key landmarks to identify include the \gls{lpv}, \gls{rpv}, and the \gls{laa} in the \gls{la}, 
	and the \gls{svc}, \gls{ivc}, \gls{cs}, and \gls{raa} in the \gls{ra}. 
	The automatic landmarking process is based on the computation of the curvature of the surface mesh, 
	as presented in \cite{azzolin2023:_augmenta} (see Figure \ref{Fig:image_process}(a)). 
	The curvature is defined point-wise, allowing for clustering of the mesh nodes based on curvature thresholding. 
	Clusters of nodes with the highest curvature -- five on the \gls{la} and four in the \gls{ra} -- are considered for anatomical structure selection and marked using the following algorithm (see Figure \ref{Fig:image_process}(b)):
	\begin{itemize}
		\item \textbf{\gls{la}:} The first landmark selected is the vertex of the \gls{laa}, identified as the central point of the \gls{la} cluster of node closest to the ascending aorta, or, if the aorta is not labeled, 
		as the cluster with the highest curvature. 
		A point cloud algorithm based on \gls{pca} is then employed to divide the remaining clusters into two groups: 
		the two clusters nearest to the \gls{laa} are designated as \gls{lpv}, 
		while the others are identified as \gls{rpv}.
		\item \textbf{\gls{ra}:} Similarly to the \gls{la}, the vertex of the \gls{raa} is first selected as the central point of the \gls{ra} cluster closest to the ascending aorta or the one with the highest curvature. 
		\gls{pca} is then applied to the two largest clusters in terms of nodes to identify the vena cavae. 
		The \gls{svc} is designated as the cluster nearest to the \gls{raa}, while the other is marked as \gls{ivc}. 
		Among the remaining clusters with the highest curvature, the one nearest to the \gls{ivc} is labeled as \gls{cs}.
	\end{itemize}
	
	The second stage entails 
	the generation of a volumetric hexahedral mesh of \gls{la} and \gls{ra} blood pools 
	from the initial volumetric segmentation of the \gls{scn}.
	For each atrium, a series of \gls{ld} problems is solved 
	to define a set of measures of distance between the previously defined veins and appendages.
	Dirichlet boundary conditions are generated by projecting the landmarks defined on the blood pool surface mesh 
	onto the volumetric mesh (see Figure \ref{Fig:image_process}(c)). 
	Isosurfaces of the \gls{ld} solves are then employed to label the veins ostiae 
	according to prescribed thresholds. 
	The remaining portions of the veins are considered discardable tissue, as they will not be part of the final anatomical model. Therefore, they are marked for removal. 
	Hereon, we will define this landmark by employing the name of the vessel followed by a negative sign, 
	e.g.\ the portion of discardable tissue of the \gls{svc} will be denoted as \gls{svc}$^-$  (see Figure \ref{Fig:image_process}(d)). 
	The same procedure is also applied to the \gls{cs}. 
	The labels defined on the volumetric hexahedral mesh are finally projected back onto an image stack
	that will then contain all augmented fine-grained labeling information.
	
	The high variability in the atrial anatomy occasionally led the automatic algorithm 
	to mislabel some structures. 
	Specifically, errors in the positioning of the \gls{ivc} and \gls{cs} were witnessed, 
	or in the identification of a reasonable portion of discardable \gls{ivc} blood pool volume. 
	In such cases, the workflow visualization modality allows for a swift manual correction.
	
	As the atrial walls are thin, their accurate segmentation from current routine clinical imaging
	with \gls{ct} and \gls{mri} is, in general, not feasible \citep{bishop2016three,dewland2013use}. 
	In our workflow, we therefore use \emph{a priori} knowledge \citep{azzolin2020influence,beinart2011left,varela2017novel,whitaker2016role} 
	on the \gls{ra} and \gls{la} endo- and epicardial width to construct the atrial walls.
	A two steps process is applied. First, atrial walls are grown by image operations 
	applied to the image stack, followed by a mesh generation step 
	to create conformal unstructured meshes with smooth boundaries.
	
	Using the augmented volumetric segmentation, atrial walls are grown by imposing a series of rule-based extrusion operations at the interface of two different landmarks, or at the interface of marked tissue with the image background. 
	Extrusion operations may comprise both growth in an inward direction, \emph{i.e.} erosion, 
	or outward, corresponding to negative and positive growth along the surface normal of the labeled domain.
	The spatial granularity of extrusion is prescribed in physical length units, i.e. mm, 
	which is converted to discrete steps corresponding to the voxel size of \SI{400}{\micro \meter} of the segmented image. 
	
	
	Specifically, volumetric representations of the atrial walls are then generated by the following rules:
	i) an inclusion set of labels where interfaces between differently labeled domains are allowed to grow into each other, comprising the blood pool labels of \gls{ra}, and \gls{la},
	the ostiae of both caval veins, \gls{svc} and \gls{ivc}, and the pulmonary veins, \gls{rpv} and \gls{lpv}; 
	ii) an exclusion set comprising domains labeled as not pertaining to the atrial blood pool, 
	to impede extrusion along these interfaces. This set includes the \gls{lv} and \gls{rv} blood pools, 
	the \gls{cs}, the \gls{svc}$^-$, \gls{ivc}$^-$, as well as right and left pulmonary veins, \gls{rpv}$^-$, \gls{lpv}$^-$ (see Figure \ref{Fig:image_process}). 
	The extrusion operations are performed then on the inclusion set, with the exclusion set serving as boundaries
	that impede any further growth.
	First, an inward extrusion is performed to create a minimal endocardial surface layer, 
	followed by an outward dilation to thicken the endocardial layer to the prescribed width,
	and a further outward extrusion to create an epicardial layer. 
	Distinct values for the thickness of the of \gls{la} and \gls{ra} endocardial and epicardial layers,
	as well as the wall width of each vein can be specified independently.
	As the exclusion set impedes any inward or outward extrusion,
	growth is then blocked at the interfaces of \gls{ra} and \gls{rv}, and \gls{la} and \gls{lv} blood pools, respectively, 
	to create the orifices of \gls{tv} and \gls{mv}. Similarly, atrial walls are not grown at the interface of the \gls{ra} and \gls{la} blood pool with the tissue of the veins and \gls{cs} to be dismissed, thus effectively creating the openings of all in- and outflow anatomical structures. 
	Additional rules may be applied to mitigate effects due to limited segmentation accuracy.
	For instance, the growth of the epicardial wall may be constrained around the ostiae of the \gls{lspv} 
	in cases where \gls{laa} and the \gls{lspv} are in close proximity, to avoid merging. 
	
	After the extrusion phase, the walls of both atria are volumetrically defined on the augmented image stack. 
	The generation of a volumetric anatomical mesh starts with extracting a surface mesh
	that encloses all atrial wall labels.
	This initial surface mesh, conforming to a jagged voxel representation,
	is subjected to several remeshing and smoothing steps to obtain a smooth representation of all atrial walls, with smooth transitions between wall segments of different widths. Topological corrections, mesh quality checks, and improvements are moreover implemented to ascertain topological correctness and sufficiently high mesh quality. Finally, the surface mesh is resampled to match a prescribed target resolution. 
	Typically, for \gls{re} models a resolution of \SI{0.9}{\milli\meter}, sufficient to resolve endo- and epicardial atrial wall layers, was selected.
	For \gls{rd} simulation, a higher resolution of \SI{0.25}{\milli \meter} was used 
	to resolve slowly propagating wavefronts without producing numerical artifacts \citep{bishop24:_devarp_trap}.
	As a final step, volumetric meshes of the atrial walls are created by the inward meshing of the labeled surfaces, 
	again followed by a topological correction, reindexing, refinement, 
	and mesh quality-enhancing procedures.
	Mesh quality of the final volumetric mesh is measured as in \cite{karabelas2018:_towards} 
	based on \cite{batdorf1997computational,kanchi20073d,knupp2022worst} (referred to Figure \ref{Fig:dbcs}-left).
	
	Thus, at the end of this stage, the biatrial anatomy is represented by an unstructured volumetric tetrahedral mesh endowed with basic anatomical labels, including endo- and epicardial domains of \gls{ra}, \gls{la}, \gls{rpv}, \gls{lpv}, \gls{ivc} and \gls{svc}. 
	
	
	\subsubsection{Automated selection of the atrial orifices}
	\label{Subsubsection:orifices}
	The atrial orifices comprising \gls{lpv}, \gls{rpv}, \gls{mv}, \gls{tv}, 
	\gls{svc}, \gls{ivc}, and \gls{cs} are individually labeled as separated anatomical entities.
	Specifically, following \cite{azzolin2023:_augmenta},    
	orifices are detected as mesh nodes of the biatrial surface 
	that are shared between adjacent elements with disjoint endo- or epicardial labels of the same anatomical structure. 
	Each orifice is therefore represented as a cluster $\Gamma$ of mesh nodes 
	arranged along a ring, and
	is readily identified based on anatomical labels, 
	except for the differentiation into superior and inferior of the \gls{rpv} and \gls{lpv}. 
	There, each entity, \gls{rpv} and \gls{lpv}, consists of two clusters which are discriminated based on the Euclidian distances between their centroid and those of the \gls{ivc}, \gls{svc}, \gls{mv}, \gls{mv} as well as the tips of \gls{raa} and \gls{laa}. 
	For instance, the centroid of the \gls{lspv} is closer to the centroid of the \gls{svc} and to the \gls{raa} tip. 
	Given the identified orifices, the anatomical regions of the \gls{la} labeled as \gls{rpv} and \gls{lpv} are updated 
	to separate the inferior and superior pulmonary veins.
	
	At this stage, the sets of nodes spanning the orifices $\Gamma_{ivc}$, $\Gamma_{svc}$, $\Gamma_{tv}$, $\Gamma_{cs}$ on the \gls{ra}, and
	$\Gamma_{rpv} = \Gamma_{ripv} \cup \Gamma_{rspv}$, $\Gamma_{lpv} = \Gamma_{lpv} \cup \Gamma_{lpv}$, and $\Gamma_{mv}$, are known 
	and can be used as boundary conditions 
	(see Figure \ref{Fig:dbcs} for a graphical representation 
	of the identified clusters), along with the corresponding anatomical regions, as represented in Figure \ref{Fig:wall_to_fibers}(a). 
	

	
	
	
	
	\subsubsection{Automated labeling of anatomical structures and fiber generation}
	\label{Subsubsec:anatomical_structure}
\begin{figure}[!t]
		\centering
		\includegraphics[width=0.9\textwidth]{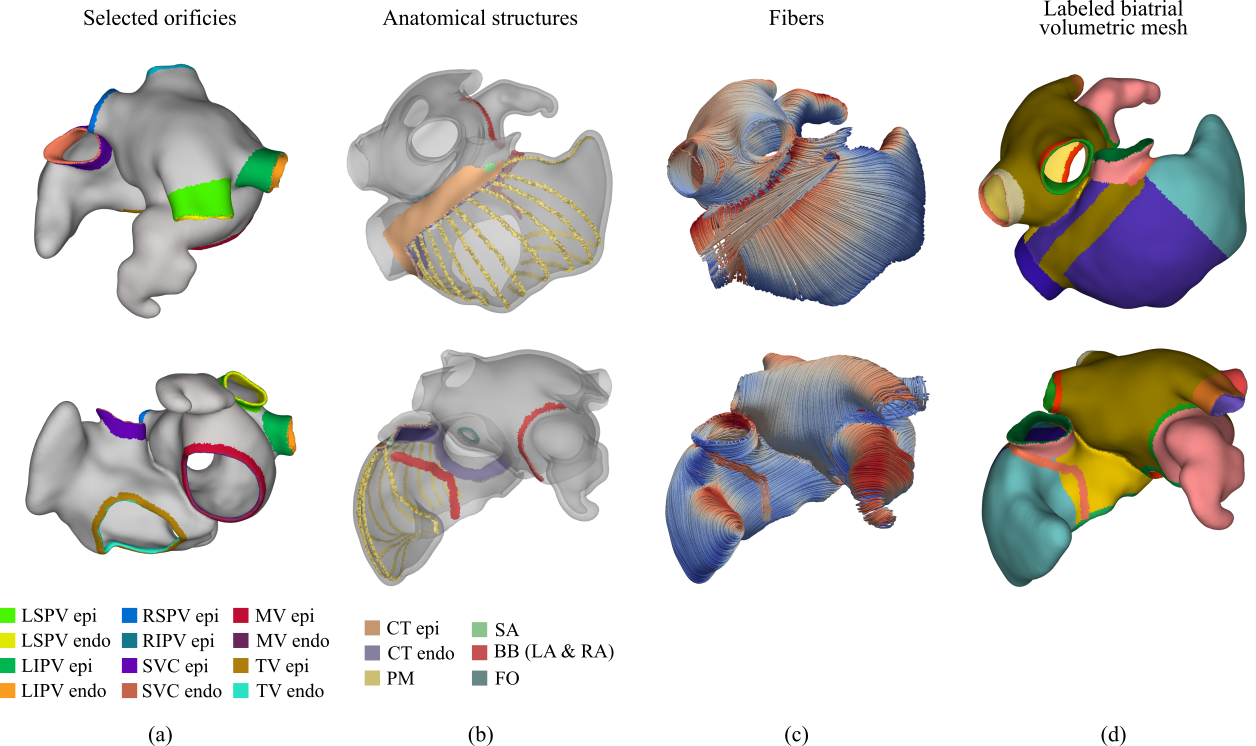}
		\caption{(a) Identification and labeling of the atrial orifices. 
			Different labels are used for endocardial and epicardial tissue layers. 
			(b) Anatomical structures identified on the biatrial anatomy, 
			including \gls{ct}, \gls{san}, \gls{pm}s, \gls{fo}, and part of the \gls{bb}. 
			(c) Generated rule-based atrial fiber architecture. 
			(d) Final biatrial volumetric mesh, augmented with labeled anatomical structures.}
		\label{Fig:wall_to_fibers}
	\end{figure}

	Following \cite{azzolin2023:_augmenta,ZhengAzzolinSanchezDosselLoewe2021},
	a further refined sub-classification of anatomical structures is performed 
	on the \gls{ra} and \gls{la} walls to define anatomical structures 
	and tissue regions of known differences in \gls{ep} behaviors. 
	Such structures -- including the \gls{san}, \gls{ct} and \glspl{pm}, the \gls{bb}, 
	and the muscular rim of the \gls{fo}, as well as \gls{raa} and \gls{laa} -- 
	are automatically annotated on a per-rule basis defining 
	location and width of the sought-after anatomical structures (refer to Figure \ref{Fig:wall_to_fibers}b).
	
\begin{figure}[!t]
		\centering
		\includegraphics[width=0.9\textwidth]{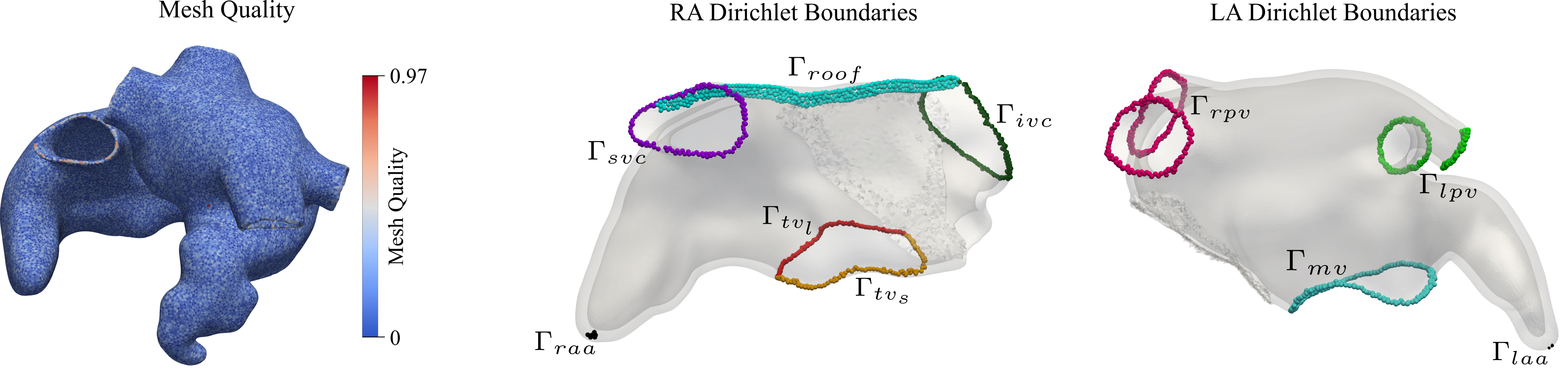}
		\caption{Left: Mesh quality computed on the biatrial model. 
			Right: Boundaries on the \gls{ra} and \gls{la} 15 \gls{ld} problems employed to determine atrial anatomical structures and tissues known to have different \gls{ep} properties.}
		\label{Fig:dbcs}
	\end{figure}

	The implemented rules are based on solutions of multiple instances of \gls{ld} problems 
	on the domain $\Omega_{myo}$, i.e. the right $\Omega_{\mathrm{RA}}$ or left $\Omega_{\mathrm{LA}}$ atrium, with varying boundary conditions. 
	Specifically, given a generic unknown $\varphi$, this involves solving
	\[ \begin{cases}
	-\Delta \varphi = 0 &\text{in } \Omega_{myo},\\
	\varphi = \varphi_i &\text{on } \Gamma_i,~i=1,\dots,m\\
	\nabla \varphi \cdot \mathbf{n} = 0 &\text{on } \Gamma_n,
	\end{cases}
	\]
	for $m = 2,3,4$ suitable \gls{dbc} $\varphi_i \in \R$, 
	set on generic partitions of the atria boundaries $\Gamma_i$, and Neumann boundary conditions $\Gamma_n$, 
	such that $\cup_{i = 1}^{m}\Gamma_i \cup \Gamma_n = \partial \Omega_{myo}$. The \gls{ld} solutions are used as measures of distance between the orifices and/or, the appendages. 
	Moreover, a plane passing through the center of the \gls{ivc} and \gls{svc}, and the center of mass of the \gls{ra} is computed to separate the \gls{ra} into lateral and septal regions. Hereon, we will call the intersection of such plane with the \gls{ra} wall as \gls{ra} roof.
	Similarly, the \gls{tv} annulus is partitioned into 
	a lateral and a septal tricuspid valve annulus, 
	$\Gamma_{tv_{l}}$ and $\Gamma_{tv_{s}}$, respectively. 
	A representation of the defined set of boundaries is shown in Figure \ref{Fig:dbcs}. 
	
	Compared to previous work \citep{azzolin2023:_augmenta,ZhengAzzolinSanchezDosselLoewe2021}, 
	three additional \gls{ld} solutions were computed on the \gls{ra}: 
	i) $\varphi_{v3}$ with \glspl{dbc} $\varphi_1 = 0$ and $\varphi_2 = 1$ imposed at the rings of \gls{svc} and \gls{ivc}, respectively; 
	ii) $\varphi_{r2}$ with $\varphi_1 = 0$ at the rings of the \gls{svc} and \gls{ivc}, and on the \gls{ra} roof, 
	and $\varphi_2 = 1$ at the ring of the \gls{mv}; and 
	iii) $\varphi_{w2}$ with \glspl{dbc}s $\varphi_1 = 0$ on the roof, $\varphi_2 = 1$ and $\varphi_3 = -1$ on the septal and lateral parts, respectively, at the \gls{mv} annulus. 
	The complete set of \glspl{dbc} is summarized in Table \ref{Table:dbcs}.
	These additional solutions are designed to improve the definition of the \gls{ct} and \glspl{pm}.
	
	\begin{table}[!t]
		\centering
		\begin{tabular}{lccccccccc}
			\toprule
			Atrium &$\varphi$ &$\varphi_1$ &$\Gamma_1$ &$\varphi_2$ &$\Gamma_2$ &$\varphi_3$ &$\Gamma_3$ &$\varphi_4$ &$\Gamma_4$ \\
			\hline \\[-2ex]
			\multirow{9}{*}{\gls{ra}} &$\varphi_{trans}$ &0 &$\Gamma_{endo}$ &1 &$\Gamma_{epi}$ &&&&\\
			&$\varphi_{ab}$ &-1 &$\Gamma_{raa}$ &0 &$\Gamma_{svc}$ &1 &$\Gamma_{tv_{l}}\cup \Gamma_{tv_{s}}$ &2 &$\Gamma_{ivc}$\\
			&$\varphi_{v}$ &0 &$\Gamma_{svc}\cup\Gamma_{raa}$ &1 &$\Gamma_{ivc}$ &&&&\\
			&$\varphi_{v2}$ &0 &$\Gamma_{ivc}$ &1 &$\Gamma_{raa}$ &&&&\\
			&$\varphi_{v3}$ &0 &$\Gamma_{svc}$ &1 &$\Gamma_{ivc}$ &&&&\\
			&$\varphi_{r}$ &0 &$\Gamma_{roof}$ &1 &$\Gamma_{tv_{l}}\cup \Gamma_{tv_{s}}$ &&&&\\
			&$\varphi_{r2}$ &0 &$\Gamma_{svc}\cup\Gamma_{roof}\cup\Gamma_{ivc}$ &1 &$\Gamma_{tv}$ &&&&\\
			&$\Gamma_{w}$ &-1 &$\Gamma_{tv_{l}}$ &1 &$\Gamma_{tv_{s}}$ &&&&\\
			&$\Gamma_{w2}$ &-1 &$\Gamma_{tv_{l}}$ &0 &$\Gamma_{roof}$ &1 &$\Gamma_{tv_{s}}$ &&\\
			\hline \\[-2ex]
			\multirow{6}{*}{\gls{la}} &$\varphi_{trans}$ &0 &$\Gamma_{endo}$ &1 &$\Gamma_{epi}$ &&&&\\
			&$\varphi_{ab}$ &-1 &$\Gamma_{laa}$ &0 &$\Gamma_{lpv}$ &1 &$\Gamma_{mv}$ &2 &$\Gamma_{rpv}$\\
			&$\varphi_{ab2}$ &0 &$\Gamma_{rpv}$ &1 &$\Gamma_{laa}$ &&&&\\
			&$\varphi_{v}$ &0 &$\Gamma_{lpv}$ &1 &$\Gamma_{rpv}$ &&&&\\
			&$\varphi_{r}$ &0 &$\Gamma_{rpv}\cup\Gamma_{lpv}\cup\Gamma_{mv}$ &1 &$\Gamma_{mv}$ &&&&\\
			&$\varphi_{r2}$ &0 &$\Gamma_{rpv}\cup\Gamma_{lpv}$ &1 &$\Gamma_{mv}$ &&&&\\
			\bottomrule
		\end{tabular}
		\caption{\glspl{dbc} for the 15 \gls{ld} problems employed to determine atrial anatomical structures and tissues known to have different \gls{ep} properties.}
		\label{Table:dbcs}
	\end{table}
	
	The labeling of the \gls{bb} follows \cite{ZhengAzzolinSanchezDosselLoewe2021}. 
	Using the additional solutions, we marked the area of the exit site of the \gls{san}, and, optionally, the muscular rim of the \gls{fo}, 
	and the anterior-central band of the \gls{bb}. 
	The \gls{san} is defined as a sphere with a radius of \SI{2.5}{\milli\meter} and centered on the node closest to the \gls{svc} with minimum $\varphi_{v2}$. 
	The \gls{fo} is defined as an annulus with thickness of \SI{2}{\milli\meter}, 
	located in the center of the atrial septum.
	The start and end points of the anterior-central band of the \gls{bb} 
	are selected as two points in the middle of the \gls{bb} lateral bands 
	of both \gls{ra} and \gls{la}.
	Additionally, an intermediate point is chosen as the epicardial node of the \gls{la} nearest to the midpoint of the \gls{ra} lateral \gls{bb}. 
	A geodesic path is, moreover, created using Dijkstra's algorithm between the starting and intermediate points, and between the intermediate and end points. 
	Around this path, a region of \SI{2}{\milli\meter} width is marked as \gls{bb}. 
	Unlike in \cite{azzolin2023:_augmenta,ZhengAzzolinSanchezDosselLoewe2021}, 
	no additional volumetric mesh is added to define the central bundle of the \gls{bb}. 
	While our workflow allows for marking the complete \gls{bb}, 
	in this work, we replace the anterior-central band of the \gls{bb}, 
	as well as all other inter-atrial conduction pathways,
	by conducting cables (see Section \ref{Subsubsection:cables}),
	increasing the flexibility of the modeling approach.
	
	Finally, the fiber architecture is computed using the rule-based method 
	presented in \cite{Piersanti2021}, 
	including the enhancements previously reported in \cite{azzolin2023:_augmenta,ZhengAzzolinSanchezDosselLoewe2021}. 
	A graphical overview of all anatomical structures, tissues, and fibers, 
	is given in Figure \ref{Fig:wall_to_fibers}.
	
	\subsubsection{Universal atrial coordinates generation}
	\label{Subsubsec:UAC}
\begin{figure}[!t]
		\centering
		\includegraphics[width=0.9\textwidth]{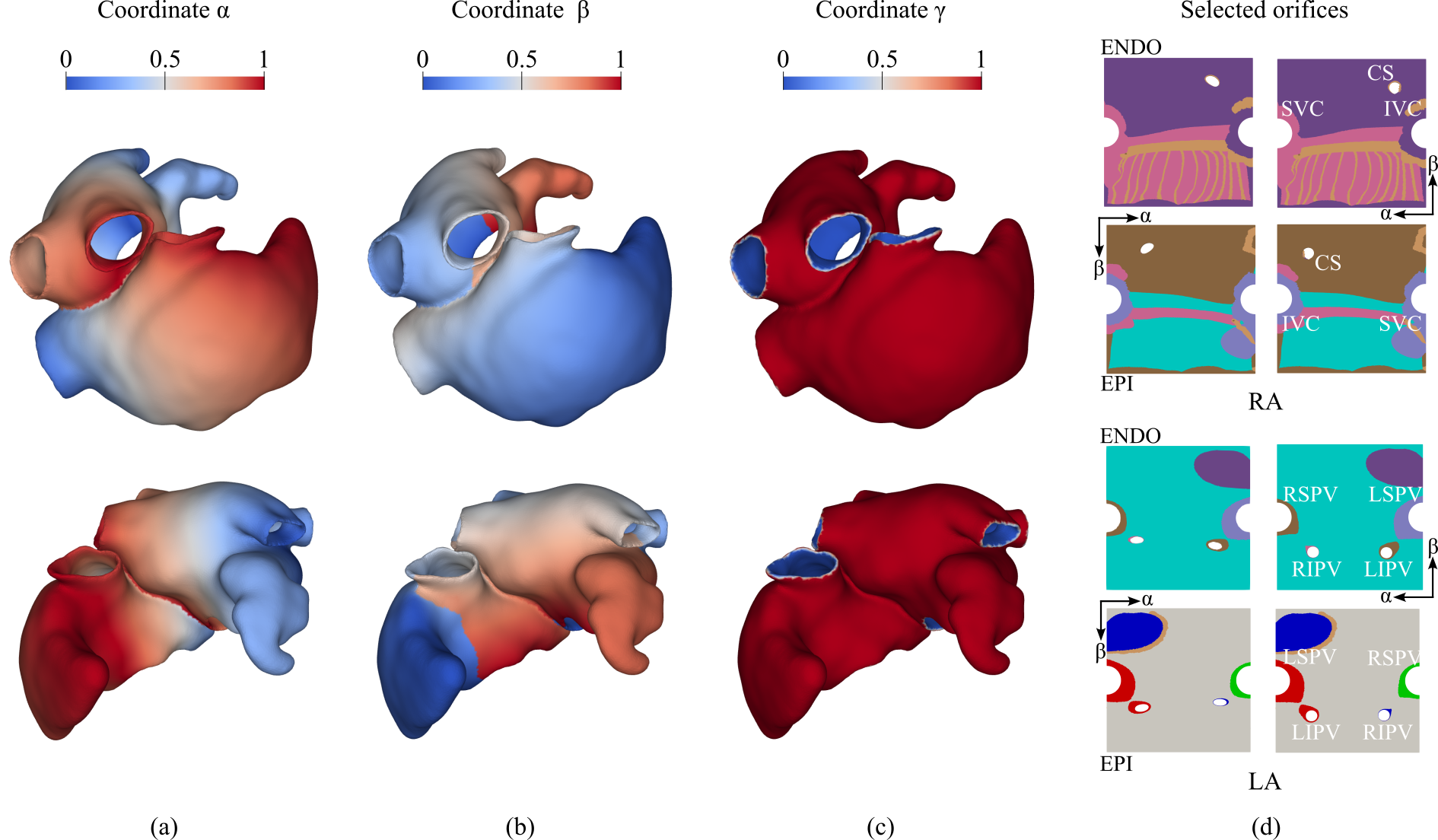}
		\caption{(a) Distribution of the 
			coordinate $\alpha$, representing the $\gls{svc}$-to-\gls{ivc} coordinate for the $\gls{ra}$, and the lateral-to-septal 
			coordinate for the \gls{la}. (b) Distribution of the coordinate $\beta$, representing the lateral-to-septal 
			coordinate for the \gls{ra}, and the posterior-to-anterior coordinate for the \gls{la}. (c) Distribution 
			of the coordinate $\gamma$, representing the endocardial-to-epicardial coordinate for both atria. (d) Projection of the labels on the space generated by the \gls{uac} before (left) and after (right) the solution of the linear elasticity problem. }
		\label{Fig:uac}
	\end{figure}

	To compute the \glspl{uac} directly on the volumetric models 
	we expand upon the method outlined in \cite{roney2023:_bia_vol,roney2021constructing}. For each atrium, we define four coordinates:
	i) $\alpha_{\mathrm{RA}}$, the \gls{ivc}-to-\gls{svc} coordinate on the \gls{ra},
	and $\alpha_{\mathrm{LA}}$, the lateral-to-septal coordinate on the \gls{la} (see illustration in Figure \ref{Fig:uac}(a));
	ii) $\beta_{\mathrm{RA}}$, the lateral-to-septal coordinate on the \gls{ra} 
	that starts at the lateral \gls{tv} annulus, runs through the roof, and ends at the septal \gls{tv} annulus,
	and $\beta_{\mathrm{LA}}$, the posterior-to-anterior coordinate on the \gls{la}
	that starts at the posterior \gls{mv} annulus, runs through the roof, and ends at the anterior \gls{mv} annulus  
	(see Figure \ref{Fig:uac}(b));
	iii) $\gamma_{\mathrm{RA}}$ and $\gamma_{\mathrm{LA}}$, the relative transmural distance 
	between endocardium and epicardium in \gls{ra} and \gls{la}, respectively (see Figure \ref{Fig:uac}(c));
	iv) a binary value encoding to which atrium the other coordinates belong, always assigning shared points at the \gls{la}--\gls{ra} interface to the \gls{la}.
	To compute the four coordinates, for each atrium, we first solve three \gls{ld} problems,
	imposing \gls{dbc} at pre-defined boundary surfaces. Secondly, to ascertain equal locations of the \gls{cs}, \gls{lipv}, and \gls{ripv} in the \gls{uac} space across different atrial models, a linear elasticity problem is solved.
	
\begin{figure}[!t]
		\centering
		\includegraphics[width=0.7\textwidth]{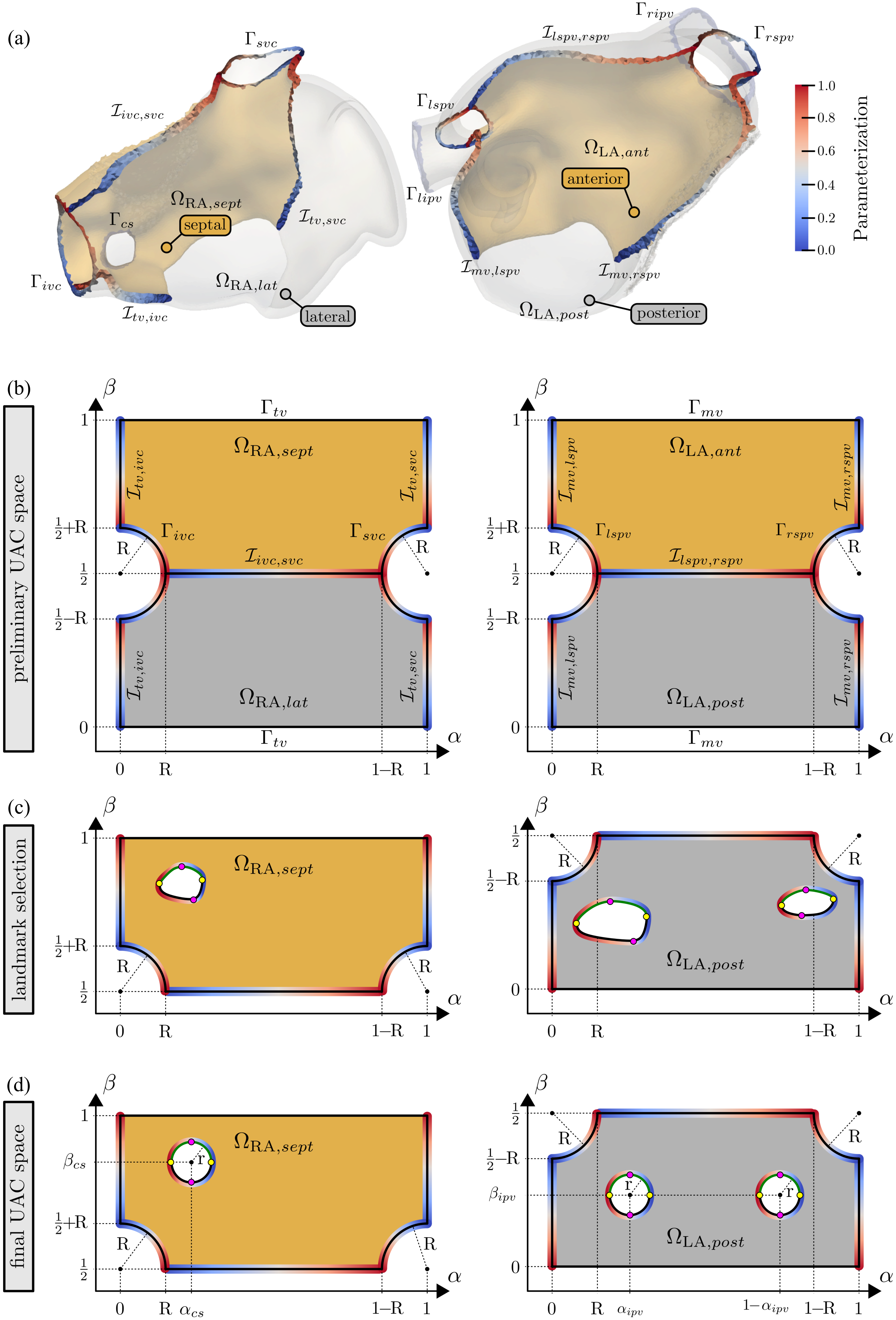}
		\caption{(a) Interfaces and boundary surfaces in the \gls{ra} and \gls{la} and their parametrization.
			(b) Dirichlet values at the boundaries and interfaces used to compute the preliminary 
			$\alpha$ and $\beta$ components for the \gls{ra} and \gls{la} are shown. (c) Selected landmarks and the parametrization of the \gls{cs}, the \gls{lipv}, and the \gls{ripv}. (d) Final $\alpha$ and $\beta$ components after moving the orifices. The color are coded from blue $\!=\!0$ to red $\!=\!1$. }
		\label{Fig:UCAparametrization}
	\end{figure}

	The selection of the boundary surfaces is fully automated and based on the anatomical labels.
	Specifically, the marked atrial orifices are utilized to define a set of reproducible points on both the vessel openings and the atrial walls. 
	Such points are selected according to the general shape of the atria, to follow the roof of the \gls{ra}, between the \gls{svc} and the \gls{ivc}, 
	and of the \gls{la}, between the \gls{lspv} and \gls{rspv}, and to set lateral and septal boundaries in both atria. 
	The points are then connected with quadrilateral surfaces. 
	The intersections of such surfaces with the atrial volumetric walls define six auxiliary surface meshes, three for each atrium, 
	effectively representing the interfaces separating the
	\gls{ra} into a septal 
	$\Omega_{\mathrm{RA},sept}$ and a lateral  $\Omega_{\mathrm{RA},lat}$ part, 
	and the \gls{la} in an anterior $\Omega_{\mathrm{LA},ant}$ and a posterior $\Omega_{\mathrm{LA},post}$ part. 
	Moreover, to avoid splitting of the appendages, the \gls{raa} is assigned to the lateral side of the \gls{ra}, 
	and the \gls{laa} is assigned to the anterior 
	side of the \gls{la}. 
	Therefore, for the \gls{ra}  we obtain two interfaces, $\mathcal{I}_{tv,ivc}$ and $\mathcal{I}_{tv,svc}$, 
	from the \gls{tv} to the \gls{ivc} and \gls{svc}, respectively, 
	and one interface $\mathcal{I}_{\mathrm{ivc,svc}}$ along the roof from \gls{ivc} to \gls{svc}. 
	Similarly, the two interfaces $\mathcal{I}_{mv,lspv}$ and $\mathcal{I}_{mv,rspv}$,  
	from the \gls{mv} to the \gls{lspv} and \gls{rspv}, respectively, are computed for the \gls{la}, 
	together with the interface $\mathcal{I}_{lspv,rspv}$ along the roof from \gls{lspv} to \gls{rspv}. 
	The defined interfaces are furthermore used to partition the \gls{tv}, \gls{ivc}, and \gls{svc} into septal and lateral parts, 
	and the \gls{mv}, \gls{lspv}, and \gls{rspv} into anterior and posterior (see Figure \ref{Fig:UCAparametrization} for a graphical representation). 
	Such partition will be used in the following when solving the linear elasticity problem.
	
	Both the interfaces and the atrial openings are then parametrized with a value $s \in [0, 1]$ with a Bi-Eikonal approach \citep{schuler2021cobiveco}.
	The $\alpha$ and $\beta$ coordinates are computed by solving \gls{ld} problems 
	on the partitioned $\Omega_{\mathrm{RA}}$ and $\Omega_{\mathrm{LA}}$, 
	with \gls{dbc} assigned based on $s$, as defined in Table \ref{Tab:UACdirichlet}.
	For the $\gamma$ component, only two Dirichlet boundary conditions are imposed, 
	$0$ at the endocardium and $1$ at the epicardium.
	
	\begin{table}
		\centering
		\begin{tabularx}{\textwidth}{X|cc|cc}
			\toprule
			& \multicolumn{2}{c|}{$\Omega_{\mathrm{RA},sept}\;\;$ and $\;\;\Omega_{\mathrm{LA},ant}$} & 
			\multicolumn{2}{c}{$\Omega_{\mathrm{RA},lat}\;\;$ and $\;\;\Omega_{\mathrm{LA},post}$} \\
			& $\alpha$ & $\beta$ & $\alpha$ & $\beta$ \\
			\midrule
			$\Gamma_{tv}$, $\Gamma_{mv}$
			& -- & 1 
			& -- & 0 \\[0.15cm]
			$\Gamma_{ivc}$, $\Gamma_{lspv}$
			& $\mathrm{R}\,\cos\Big((1\!-\!s)\,\frac{\pi}{2}\Big)$ 
			& $\frac{1}{2}\!+\!\mathrm{R}\,\sin\Big((1\!-\!s)\,\frac{\pi}{2}\Big)$
			& $\mathrm{R}\,\cos\Big((1\!-\!s)\,\frac{\pi}{2}\Big)$ 
			& $\frac{1}{2}\!+\!\mathrm{R}\,\sin\Big((3\!+\!s)\,\frac{\pi}{2}\Big)$ \\[0.15cm]
			$\Gamma_{svc}$, $\Gamma_{rspv}$
			& $1\!+\!\mathrm{R}\,\cos\Big((1\!+\!s)\,\frac{\pi}{2}\Big)$ 
			& $\frac{1}{2}\!+\!\mathrm{R}\,\sin\Big((1\!-\!s)\,\frac{\pi}{2}\Big)$
			& $1\!+\!\mathrm{R}\,\cos\Big((1\!+\!s)\,\frac{\pi}{2}\Big)$ 
			& $\frac{1}{2}\!+\!\mathrm{R}\,\sin\Big((3\!+\!s)\,\frac{\pi}{2}\Big)$ \\[0.15cm]
			$\mathcal{I}_{tv,ivc}$, $\mathcal{I}_{mv,lspv}$
			& $0$ & $\Big(\frac{1}{2}\!+\!\mathrm{R}\Big)+(1\!-\!s)\,\Big(\frac{1}{2}\!-\!\mathrm{R}\Big)$ 
			& $0$ & $s\,\Big(\frac{1}{2}\!-\!\mathrm{R}\Big)$ \\[0.15cm]
			$\mathcal{I}_{tv,svc}$, $\mathcal{I}_{mv,rspv}$
			& $1$ & $\Big(\frac{1}{2}\!+\!\mathrm{R}\Big)+(1\!-\!s)\,\Big(\frac{1}{2}\!-\!\mathrm{R}\Big)$ 
			& $1$ & $s\,\Big(\frac{1}{2}\!-\!\mathrm{R}\Big)$ \\[0.15cm]
			$\mathcal{I}_{ivc,svc}$, $\mathcal{I}_{lspv,rspv}$
			& $\mathrm{R}\!+\!(1\!-\!2\,\mathrm{R})\,s$ & $\frac{1}{2}$ 
			& $\mathrm{R}\!+\!(1\!-\!2\,\mathrm{R})\,s$ & $\frac{1}{2}$ \\[0.15cm]
			\bottomrule
		\end{tabularx}
		\caption{Dirichlet boundary conditions for the Laplace problems that are solved to obtain 
			the preliminary \gls{uac} space in the septal and lateral \gls{ra} and the posterior and 
			anterior \gls{la} where $s \in [0, 1]$ is a parametrization value of the corresponding
			boundary surface and interface, respectively. In our framework,  $\mathrm{R}=0.1$.}
		\label{Tab:UACdirichlet}
	\end{table}
	\begin{table}
		\begin{tabularx}{0.6\textwidth}{l|cc}
			\toprule
			& $\alpha$ & $\beta$ \\
			\midrule
			$\Gamma_{cs,roof}$
			& $\alpha_{cs}\!+\!\mathrm{r}\,\cos\Big((2\!-\!s)\,\pi\Big)$
			& $\beta_{cs}\!+\!\mathrm{r}\,\sin\Big((2\!-\!s)\,\pi\Big)$ \\[0.15cm] 
			$\Gamma_{cs,tv}$
			& $\alpha_{cs}\!+\!\mathrm{r}\,\cos\Big(s\,\pi\Big)$ 
			& $\beta_{cs}\!+\!\mathrm{r}\,\sin\Big(s\,\pi\Big)$ \\[0.15cm]
			$\Gamma_{lipv,roof}$
			& $\alpha_{ipv}\!+\!\mathrm{r}\,\cos\Big(s\,\pi\Big)$ 
			& $\beta_{ipv}\!+\!\mathrm{r}\,\sin\Big(s\,\pi\Big)$ \\[0.15cm]
			$\Gamma_{lipv,mv}$
			& $\alpha_{ipv}\!+\!\mathrm{r}\,\cos\Big((2\!-\!s)\,\pi\Big)$
			& $\beta_{ipv}\!+\!\mathrm{r}\,\sin\Big((2\!-\!s)\,\pi\Big)$ \\[0.15cm] 
			$\Gamma_{ripv,roof}$
			& $(1\!-\!\alpha_{ipv})\!+\!\mathrm{r}\,\cos\Big(s\,\pi\Big)$
			& $\beta_{ipv}\!+\!\mathrm{r}\,\sin\Big(s\,\pi\Big)$ \\[0.15cm]
			$\Gamma_{ripv,mv}$
			& $(1\!-\!\alpha_{ipv})\!+\!\mathrm{r}\,\cos\Big((2\!-\!s)\,\pi\Big)$
			& $\beta_{ipv}\!+\!\mathrm{r}\,\sin\Big((2\!-\!s)\,\pi\Big)$ \\[0.15cm]
			\bottomrule
		\end{tabularx}
		
		\caption{Final position of the \gls{cs}, \gls{lipv}, and \gls{ripv} orifices with parameterization value $s \in [0, 1]$.}
		\label{Tab:UACdirichletLinElas}
	\end{table}
	
	The linear elasticity problem is solved in the \gls{uac} space, by projecting the first three coordinates onto two thin cuboids, one for each atrium, 
	where $\alpha$ and $\beta$ are described along the edges of the square face, and $\gamma$ represents the thickness of the cuboid. 
	The \gls{ivc}, \gls{svc}, \gls{rspv}, and \gls{lspv} are moreover recast in circumferential form, on the edges of the square face representing the $\beta$ coordinate.    
	The square faces of each cuboid are then representing the endocardial and epicardial surfaces of the atrial walls.
	An updated Lagrangian formulation is then used to solve a linear elasticity problem, 
	prescribing rule-based Dirichlet boundary conditions for the \gls{cs}, \gls{ripv}, and \gls{lipv}.    
	To achieve this, the parametrization of the \gls{cs}, \gls{ripv}, and \gls{lipv} is carried out as follows. 
	Beginning with the mesh nodes of each orifice, i.e., $\Gamma_{ripv}$, $\Gamma_{lipv}$, and $\Gamma_{cs}$ as defined in Section \ref{Subsubsection:orifices}, 
	a node selection is grown across the biatrial surface in all directions. 
	The outermost selected nodes form two rings, one on the endocardium and one on the epicardium. 
	For each ring, four points are identified, corresponding to the minimum and maximum values of the $\alpha$ and $\beta$ components. 
	and for each of the four endocardial-epicardial point pairs, the shortest path is extracted on the selected surface.
	The \gls{cs} node selection is then divided into a roof and a \gls{tv} section 
	and parameterized using a \gls{ld} solve, with Dirichlet boundary conditions applied to the eight previously extracted paths. 
	Similarly, the \gls{lipv} and \gls{ripv} node selections are divided into a roof and a \gls{mv} section and parameterized in the same manner. 
	The final position of the \gls{cs}, \gls{ripv}, and \gls{lipv} after convergence of the linear elasticity problem is defined in Table \ref{Tab:UACdirichletLinElas} 
	with $\mathrm{R}=0.1$, $\mathrm{r}=0.04$, $\alpha_{cs}=0.2$, $\beta_{cs}=0.8$, $\alpha_{ipv}=0.25$, and 
	$\beta_{ipv}=0.25$) and in Figure \ref{Fig:uac}(d), where the anatomical labels are also reported.

	\subsubsection{Torso volume conductor model}
	\label{Subsubsec:Torso}
\begin{figure}[!t]
		\centering
		\includegraphics[width=0.9\textwidth]{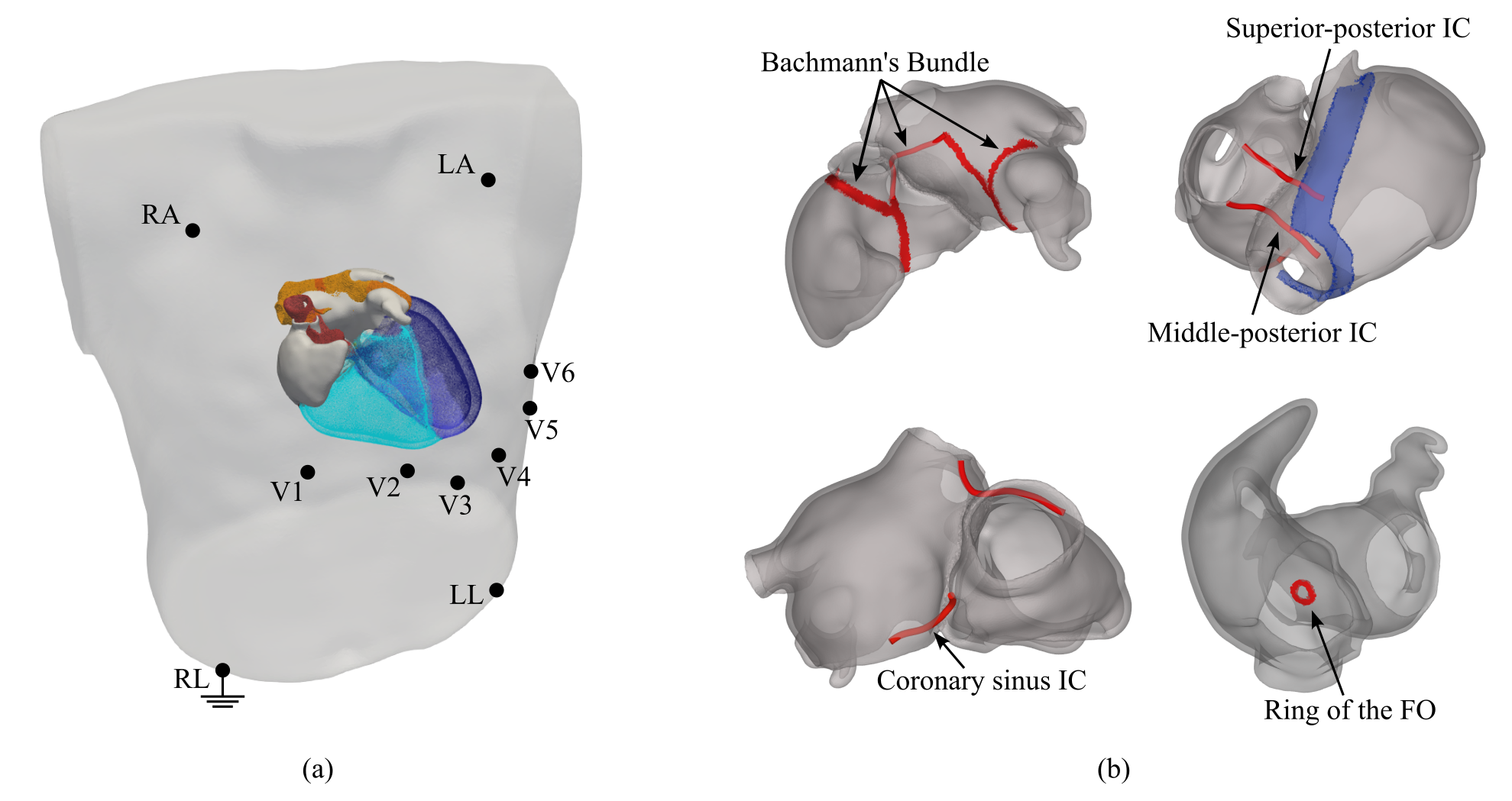}
		\caption{(a) Combined conformal atria-torso mesh generated by registering a patient atrial model, 
			generated with our workflow, to a selected template torso model with positioned \gls{ecg} electrodes. The ventricles for the template heart are also represented. 
			(b) Representation of the electrical connections between the \gls{ra} and \gls{la}, including four \gls{ic} and the \gls{fo}.}
		\label{Fig:torso_and_ica}
	\end{figure}

	Since most clinical observations of atrial \gls{ep} are obtained by sampling the potential fields generated within the volume conductor, such as \glspl{ecg} or \glspl{egm}, volumetric models must be embedded in a torso domain to enable the prediction of these observations.
	However, the field of view of the \gls{ct} scans of the patients within the Graz \gls{af} registry was restricted to the heart to minimize the patient's exposure to x-ray radiation.
	Thus a full torso view including all \gls{ecg} recording sites was not available.
	Thus, to predict atrial P-waves in the \gls{ecg}, an automated procedure was implemented for registering a template heart-torso model with known electrode locations around the generated atria. The template was selected from a previously generated virtual cohort \citep{gillette2021automated, gillette2022personalized}.
	More precisely, the atrial surfaces of the template were registered to the generated patient models by employing a global optimal 3D iterative closest point (GO-ICP) algorithm \citep{yang2015go,yang2013go}.
	The transformation was then stored and applied to the surface of the template torso
	along with the known positions of the electrodes \citep{gillette2021automated, gillette2022personalized} (see Figure \ref{Fig:torso_and_ica}(a)).
	
	A volumetric mesh was subsequently generated between the surfaces of the torso and the atria, then merged with the volumetric mesh of the atria.
	For simplicity, other organs represented in the template model, such as the lungs, were excluded to avoid potential meshing issues caused by intersections with the atrial geometry.
	
	The obtained atria-torso model is suitable for the accurate simulation of \glspl{ecg} and \glspl{egm},
	either based on a bidomain model for simulating the entire extracellular potential field, $\phi_{\rm e}(\mathbf{x},t)$  \citep{bishop2011bidomain,neic17:_reaction_eikonal,ZAPPON2024112815}, 
	or on a lead field approach \citep{geselowitz1989theory,multerer2021uncertainty,potse2018scalable},
	yielding time traces of the extracellular potential sampled at discrete points $\phi_{\rm e,\mathbf{x}}(t)$. 
	
	
	\subsection{Functional modeling stage}
	\label{Subsection:functional_modeling_stage}
	\subsubsection{Parametric modeling of inter-atrial connections}
	\label{Subsubsection:cables}
	Action potential propagation between \gls{ra} and \gls{la} is limited to 
	a few electrically excitable muscular strands of tissue called \glspl{ic}, 
	that serve as conduction pathways between the otherwise electrically insulated atria.
	Known pathways include \gls{bb}, 
	the sheath of the \gls{cs}, the muscular rim of the \gls{fo} 
	and, posteriorly, a superior-posterior and a middle-posterior bridge.
	Histologically, \glspl{ic} arise sub-epicardially and branch outside of the atrial walls \citep{Margo2013,platonov2002morphology}. 
	While the \gls{bb} \citep{Bachmann1916,Margo2013,lemery2003anatomic} and the \gls{cs} musculature \citep{chauvin2000anatomic,antz1998electrical} have been shown 
	to be the most relevant electrical \gls{ic} \citep{sakamoto2005interatrial,knol2019bachmann}, 
	the specific number and location of the complete set of \gls{ic} remains unclear \citep{platonov2002morphology,platonov2007interatrial,sakamoto2005interatrial}.    
	
	In this study, we considered all five aforementioned \glspl{ic} as the only electrical connections between the \gls{ra} and \gls{la}.
	Except for the rim of the \gls{fo}, which forms an electrically conductive bridge within the volumetric biatrial mesh,
	\gls{ra} and \gls{la} are electrically fully decoupled using a nodal splitting approach \citep{costa2014:_fibrotic}. 
	All other \glspl{ic} traversing the space outside the atrial walls 
	are modeled as discrete electrically conducting strands, 
	that are anchored within the \gls{ra} and the \gls{la} at given sites defined through \glspl{uac} 
	(see Figure \ref{Fig:torso_and_ica}(b)), and constrained to the sub-epicardium.
	These strands are defined as auto-generated cables as previously utilized in \cite{gillette2022personalized}, and their electrical behavior is modeled with the approach used to define the \gls{hps} in the ventricles in \cite{neic17:_reaction_eikonal,vigmond2006:_sawtooth}. The cables used to define the \glspl{ic} were generated along the shortest path 
	connecting the anchoring location in \gls{ra} and \gls{la}, respectively.
	
	
	Of the five represented \glspl{ic}, the \gls{bb} is composed of both an intra- and an inter-atrial component, 
	where the intra-atrial component is modeled as tissue embedded within the atrial epicardial wall,
	as described in Section \ref{Subsubsec:anatomical_structure}, 
	and only the inter-atrial central component is modeled as a cable. 
	The \gls{bb} emerges at the junction of the \gls{ra} body and the \gls{svc}, 
	in close proximity to the sino-atrial node (\gls{san}), 
	and runs in an anterior direction around the \gls{raa}, where it branches into two bands, 
	one towards the \gls{ra} vestibule, 
	and the inter-atrial central \gls{bb} running toward the anterior \gls{la} wall, anchoring in the \gls{la} epicardium.
	From there another intra-atrial band runs epicardially towards the \gls{laa}, encircles the \gls{laa}, 
	and terminates superiorly between the \gls{laa} and the \gls{lspv}, 
	and inferiorly, at the \gls{mv} annulus \citep{Margo2013}. 
	
	The two posterior \gls{ic}s were defined as in \cite{loewe2016influence,wachter2015mesh}. 
	The superior-posterior bridge emerges in proximity of the \gls{ivc} at the 90\% of the \gls{ct}, 
	and connects the \gls{ra} with the \gls{la} wall near the \gls{ripv}, on the septal side. 
	The middle-posterior \gls{ic} starts at the junction of the \gls{ivc} with the atrial body, 
	on the \gls{ct}, and ends at the junction of the \gls{ripv}, on the inferior side of the \gls{ripv}.
	
	The musculature of the \gls{cs} was represented by a single cable providing a pathway 
	from the \gls{ra} wall, near the \gls{cs} ostium, to the posterior \gls{la} at the annulus of the \gls{mv}. 
	The origin of the cable in the \gls{ra} was chosen superior to the \gls{cs}, 
	at a distance between 3 and \SI{8}{\milli\meter} \citep{chauvin2000anatomic} from the \gls{cs} ostium, 
	while the \gls{la} anchoring site was selected according to \cite{wachter2015mesh}.
	
	The muscular rim of the \gls{fo} was modeled as tissue using previously defined labels, 
	endowed with slow conduction properties \citep{kharbanda2019current,harrild2000computer,padala2021anatomy}. 
	
	
	\subsubsection{Baseline calibration of the atrial activation sequence and P-wave simulation}
	\label{Subsubsec:baseline_sim}
	The genesis of the atrial P-wave is governed by the spatio-temporal evolution 
	of electrical depolarization wavefronts traversing the atria.
	At each point in the atria, the propagation of the wavefronts 
	is governed by (i) a conduction velocity tensor 
	that is determined by intra- and extracellular conductivities along the eigenaxes, (ii)
	the tissue composition condensed into a scalar bidomain surface-to-volume ratio, 
	and (iii) the tissue's \gls{ep} properties related to cellular dynamics.
	Beyond these intrinsic conduction properties, the site of initiation of a depolarization wavefront, 
	typically at an exit site of the \gls{san} or at a pacing site,
	and the location of anchoring sites of inter-atrial conduction pathways 
	and their associated conduction velocities influence the distribution of electrical sources during depolarization
	and, thus, the morphology of the P-wave.
	As such, the parameter vector controlling atrial activation is high dimensional.
	While all these parameters are exposed for unattended manipulation in parameter sweeps, their calibration as a whole is computationally not feasible with current technology,
	and, further, observable data are very sparse and insufficient to constrain the calibration procedure.
	As such, we only demonstrate the ability of our fully parametric 3D biatrial model 
	to sweep over selected key parameters to generate physiologically meaningful sets of atrial P-waves.
	
	Our workflow universally supports all common \gls{ep} modeling approaches 
	that build on the volumetric representation of atrial electrical sources, 
	including all \gls{rd} mono- and bidomain models \citep{bishop2011bidomain,potse2006comparison}, 
	as well as variants of \gls{re} models \cite{neic17:_reaction_eikonal}, 
	including plain Eikonal models with $V_{\rm m}$-recovery \citep{pezzuto2017:_ecg}.
	To demonstrate the feasibility of a sampling-based model calibration using the P-wave,
	we chose two popular approaches, 
	the computationally expensive \gls{rd} monodomain model, 
	requiring a high average spatial resolution of $\Delta x\approx$ \SI{0.25}{\milli\meter},  
	and a lightweight \gls{re} model at $\Delta x\approx$ \SI{0.90}{\milli \meter}. Both models are combined with the lead field approach for the P-wave computation. 
	We refrain here from providing a mathematical model formulation and refer to previous studies where all methodological underpinnings of the forward \gls{ep} modeling have been described in detail \citep{bishop2011bidomain,gillette2021:_framework,neic17:_reaction_eikonal}. 
	All simulations are executed using CARPentry \citep{neic17:_reaction_eikonal} for \gls{rd} and \gls{re} simulation runs, respectively, as well as openCARP \citep{plank2021opencarp} only suitable for \gls{rd} simulations.
	
	
	\paragraph{Baseline parameter setting for electrophysiological simulations}
	The dimensionality of the parameter space governing the spatial variation in conduction velocity 
	is reduced as follows. 
	First, the overall fiber arrangement is kept constant,
	and target conduction velocities are prescribed in the \gls{ct}, \gls{pm}s, \gls{ic}, \gls{fo}, 
	and the wider \gls{la} and \gls{ra} walls,  
	on a per-region basis 
	using reported ranges \citep{deng2012simulation}, eventually scaled to ascertain that the total atrial activation times computed with the \gls{re} model
	fall within the physiological ranges \citep{lemery2007normal}.
	
	Anisotropy in conduction velocity between longitudinal and transverse to the fiber directions was set to a value of 1.3 within the reported range of 1.0 to 1.6 \citep{gray1996incomplete,hansson1998right}.
	Atrial cellular dynamics were represented by the Courtemanche model \citep{courtemanche1998ionic},
	parameterized for various \gls{ep} regions following \cite{azzolin2023:_augmenta}
	with increased sodium channel conductivity to obtain the prescribed conduction velocities.
	To identify the electrical conductivities for the \gls{rd} model at the targeted mesh resolution that match the \gls{re} conduction velocities, the ForCEPSS framework \citep{gsell2024:_forcepss} was employed, using reported values \citep{roberts1979influence} for initialization. 
	The chosen parameter settings are summarized in Table \ref{tab:calibration}.
	Conduction velocities in the \gls{ic} cables were considered parameters, 
	tuned to obtain physiological inter-atrial activation delays and, thus,
	the separation of the contributions of \gls{ra} and \gls{la} to the P-wave.
	
	All simulations were conducted for a duration of \SI[parse-numbers = false]{150}{\milli\second} 
	to cover the total atrial activation time and, thus, the entire duration of the P-wave. 
	All \glspl{ecg} were filtered with a \SI{150}{Hz} low-pass and a \SI{0.5}{Hz} high-pass, 
	and scaled by a factor of 0.2 to obtain \gls{ecg} signals with the observed magnitudes \citep{gillette2022personalized}.
	
	
	\begin{table}[!t]
		\centering
		\begin{tabular}{lllllll}
			\toprule
			& $v_{\rm l}$ & $v_{\rm t}$ &$g_{\rm{i,l}}$ & $g_{\rm {e,l}}$ &$g_{\rm{i,t}}$ & $g_{\rm {e,t}}$\\
			&m/s &m/s &S/m &S/m &S/m &S/m \\ 
			\hline \\[-2ex]
			\gls{ra}           & 0.97 & 0.74 & 0.583  & 0.742  & 0.232  & 1.162 \\
			\gls{la}           & 0.98 & 0.76 & 0.570  & 0.726  & 0.220  & 1.102 \\
			\gls{ct}           & 1.21 & 0.92 & 0.904  & 1.112  & 0.416  & 0.899 \\
			\glspl{pm}         & 1.30 & 0.99 & 1.053  & 1.339  & 0.408  & 2.042 \\
			\gls{bb} (atrial)  & 1.40 & 1.08 & 1.229  & 1.564  & 0.479  & 2.399 \\
			Rim \gls{fo}       & 0.33 & 0.24 & 0.080  & 0.102  & 0.034  & 0.168 \\
			\bottomrule
		\end{tabular}
		\caption{Conduction velocities in longitudinal and transverse directions, $v_l$ and $v_t$, 
			prescribed in the \gls{re} model along with calibrated conductivities 
			yielding matching velocities in the \gls{rd} models at the target average mesh resolution of 
			$\Delta x\approx$ \SI{0.25}{\milli \meter}. 
			The bidomain surface-to-volume ration was chosen as $\beta=$\SI{1400}{\centi \meter}$^{-1}$. }
		\label{tab:calibration}
	\end{table}
	
	
	We assessed the accuracy of the \gls{re} model relative to a gold standard \gls{rd} monodomain model for simulating atrial activation and associated P-waves.
	Using the baseline setup (see Table \ref{tab:calibration}), 
	a normal sinus activation with associated P-wave was simulated on two models of the same biatrial geometry
	discretized at different spatial resolutions of 
	\SI{0.9}{\milli\meter} and \SI{0.25}{\milli\meter} for \gls{re} and \gls{rd}, respectively,
	according to their specific numerical requirements \citep{neic2017eikonal}.
	
	
	\paragraph{Effect of the \gls{bb} insertion in the \gls{la} on the P-wave}
	
	Activation of the \gls{la} is mediated through discrete \glspl{ic}
	with \gls{bb} being a primary pathway \citep{knol2019bachmann,sakamoto2005interatrial}.
	A site of insertion, $\mathbf{x}$, of an \gls{ic} defines an earliest activation site, $\varepsilon(\mathbf{x},t)$, of the \gls{la}
	and, thus, constitutes a major factor in shaping the \gls{la} contribution to the P-wave morphology.
	With traditional modeling approaches, relying on discrete meshing of the \glspl{ic}, elucidating such effects is challenging, 
	as any change in the insertion site requires major remeshing of the model.
	In this work, we aim at investigate the impact of the location of the central \gls{bb} insertion on the \gls{la} activation and P-wave, by sweeping over the spatial parameter $\mathbf{x}_{BB_{LA}}$ representing the \gls{bb} insertion point on the \gls{la} wall, that is $\varepsilon(\mathbf{x}_{BB_{LA}},t_{la})$. The temporal parameter $t_{la}=l_{bb}/v_{bb}$ represents the activation time of $\varepsilon$, and depends on the length of the cable $l_{bb}$, and the designated conduction velocity $v_{bb}$. Here, we chose to prescribed $v_{bb}$, while varying $l_{bb}$ dependently on the insertion point $\mathbf{x}_{BB_{LA}}$. Therefore, changes in $\mathbf{x}_{BB_{LA}}$ will impact both the \gls{la} activation pattern, that is, the P-wave shape, and the initial activation time of the \gls{la}, corresponding to a shift in the portion of the P-wave related to the \gls{la} activation.
	The importance of these parameters is investigated concerning the sensitivity of the P-wave shape, and their ability to yield a large enough physiological envelope of P-waves to cover the true measured P-wave.
	
	However, using a single \gls{ic} cable to mediate \gls{bb} conduction limits the model ability to produce single point-like \gls{la} activation sites, and potentially the \gls{ecg} shape generation.    
	As the number of inter-atrial cables is not limited in our modeling approach, we investigate a more realistic volumetric-like approximation of the \gls{bb} band,
	by using multiple cables.
	Specifically, we select three points, $\mathbf{x}{BB_{LA}}^i$ ($i=1,2,3$), on the anterior wall of the \gls{la} as locations for anchoring three reference cables representing strands of the \gls{bb}, all of which originate from the midpoint of the \gls{ra} lateral band of the \gls{bb}. 
	As per the single-cable case, each location $\mathbf{x}{BB_{LA}}^i$ was then varied within a circular area of $\approx$ \SI{5}{\milli\meter}, defining a set of alternative entry points ${ \mathbf{x}{BB_{LA_{j}}}^i}{j=1}^{N{i}}$, and corresponding automatically generated cables.
	By running a \gls{re} simulation for each tuple of cables $(\mathbf{x}{BB_{LA_{j}}}^1,\mathbf{x}{BB_{LA_{k}}}^2,\mathbf{x}{BB_{LA_{\ell}}}^3)_{j=1,\dots,N_1,k=1,\dots,N_2,\ell=1,\dots,N_3}$, 
	the set of P-waves was finally computed.
	
	\paragraph{Effect of the discrete RA endocardium on the P-wave}
	The endocardium of the \gls{ra} consists of discrete structures 
	in the form of muscular strands comprising \gls{ct}, \gls{bb} and \glspl{pm},
	that are partially attached to a thin and smooth epicardial layer \citep{corradi2011atria,ho2009importance}. 
	
	The effect of the added endocardial tissue in atrial \gls{ep} simulation, 
	especially the inter-pectinate tissue, on the P-wave has never been study before.  
	To investigate the role of the \gls{ra} endocardium on the P-wave, 
	\gls{re} model simulations were run using two different parameter settings for each biatrial model. 
	In one setting, the \gls{ra} endocardium was considered cardiac conductive tissue, 
	while in the other, the tissue labeled as \gls{ra} endocardium was treated as part of the bath surrounding the heart. In the former case, conduction velocity tensors, as presented in Table \ref{Tab.RE_vs_RD}, were assigned to the \gls{ra} endocardium tissue. In the latter case, this tissue was marked as torso tissue.
	
	
	\paragraph{Calibration of the atrial activation sequence based on the P-wave of the 12 lead \gls{ecg}}
	We investigate the feasibility of calibrating the atrial activation sequence to the corresponding patient P-wave by sampling key parameters to generate a P-wave envelope covering the observed P-wave data.
	Due to the high dimensionality of the parameter space, we restrict sampling to three key parameters that are expected to be the most influential in altering the different phases of the P-wave morphology \citep{loewe2015left}. These are:
	\begin{itemize}
		\item the location of the leading pacemaker site within the \gls{san}, responsible for the initial \gls{ra} activation. 
		This site is known to be variable, allowing modulation of the onset phase of the P-wave.
		The \gls{san} is an elongated and highly heterogeneous \gls{ra} region
		with action potential properties, cell size, and capacity as well as conductance 
		changing from the center to periphery \citep{honjo_correlation_1996}, 
		and along inferior-superior and caudal-cranial gradients \citep{munoz_onset_2011}. 
		For calibration purposes, we avoided modeling the physiological function of the \gls{san} with biophysical detail and instead used only the spatial location of the exit site $\mathbf{x}_{\text{SAN}}$ as a parameter.
		Specifically, the \gls{san} was modeled as a focal activation site, with a short, cigarette-shaped band of \SI{2.2}{\centi\meter} in length \citep{boyett2000sinoatrial,monfredi2010anatomy}, located in the inferior-central part of the \gls{ct}. 
		Moreover, the \gls{san} was electrically isolated by a non-conductive region along the boundary of the \gls{ct}, replicating the effective block observed in physiological studies \citep{dobrzynski2005computer}. 
		
		\item Calibration to the overall duration of earlier phase of the P-wave 
		was based on varying conduction velocities globally within the \gls{ra}.
		Velocity in the longitudinal direction $v_l$ was uniformly sampled in steps of \SI{0.5}{\meter/\second} 
		between \SI{0.6}{} and \SI[parse-numbers = false]{1.1}{\meter/\second}.
		Assuming rotational isotropy, velocity in the transverse direction $v_t$
		was sampled in steps of \SI{0.0375}{\meter/\second} 
		between \SI[parse-numbers = false]{0.45}{} and \SI[parse-numbers = false]{0.7875}{\meter/\second}, but kept smaller to the longitudinal velocity $v_l$, in accordance with \cite{deng2012simulation}.
		
		\item Calibration of the shape and duration of the late phase of the P-wave, was based on manual variation of the anchoring sites of the \glspl{ic} and of the global conduction velocities in the \gls{la}. Specifically, the \gls{bb} was modeled as a three-cable band initiating at the same origin in the \gls{ra}, and anchoring in the \gls{la} at three independent insertion sites. The other three \glspl{ic} were also varied in the anchoring position on the \gls{la}, and conduction velocities in the cables were manually calibrated to control the timing of the initial \gls{la} activation. 
	\end{itemize}

	\paragraph{Data analysis}
	We performed both qualitative and quantitative analyses of \gls{ep} and P-wave variations resulting from the previously described modeling assumptions. The overall tissue activation was qualitatively analyzed by recording the total activation time of both atria. Changes in P-waves were quantitatively investigated using three different metrics. 
	
	In the presence of a reference signal $\phi_r^j$, either simulated or given as clinical \gls{ecg} data, we analyzed amplitude variation with the \gls{rmse} in percentage, as described in \cite{keller2010ranking,ZAPPON2024112815}. For each lead $j$, the \gls{rmse} is expressed as:
	$$ \text{\gls{rmse}}(\varphi_v^j)[\%] = \sqrt{\frac{\sum_{n=1}^N(\phi_v^j(n)-\phi_r^j(n))^2}{\sum_{n=1}^N (\phi_r^j(n))^2}} \cdot 100.$$
	We computed the \gls{rmse} when comparing the \gls{re} and \gls{rd} models when exploring the effect of the variation of the \gls{san} and \gls{ra} conduction velocity on the P-wave.
	
	When analyzing the variations of the \gls{bb}, whether using one or three cables, none of the computed \gls{ecg} signals matched the real \gls{ecg} data and could not be used as a reference to calculate the \gls{rmse}. For these test cases, we first defined an average P-wave, $\bar{\varphi}^j$, for each lead $j$. We then evaluated the overall morphological and amplitude P-wave variation by computing the mean absolute distance (MAD) of each resulting P-wave $\varphi_v^j$ from the average one as:
	\begin{equation}
	\label{eq:MAD}
	\text{MAD}(\varphi_v^j) = \frac{\sum_{t=1}^N |\varphi_v^j(t) - \bar{\varphi}^j(t)|}{N}, 
	\end{equation}
	where $\varphi_v^j(t)$ and $\bar{\varphi}^j(t)$ are the varying computed signal and the average signal at time $t$,  and $N$ is the total number of \gls{ecg} time samples. The MAD is then averaged over all leads.
	
	Additionally, in all cases, we accounted for the variation of \gls{pwd}, obtained using a sloping approach \citep{tan2000detection}. The \gls{pwd} was measured for each lead. Moreover, the \gls{pwd} averaged over all leads was given.
	
	\section{Results}
	\label{Sec:results}
	\subsection{Model generation workflow performance evaluation}
	We evaluated our atrial model generation workflow with respect to the required processing time 
	and the achieved degree of automation.
	50 contrast \gls{ct} datasets of \gls{af} patients were processed to generate biatrial models of target resolution of $\approx$ \SI{0.25}{\milli \meter} and $\approx$ \SI{0.9}{\milli \meter}, and the execution times of individual stages of the automated workflow were measured.
	After each stage, visual checks were performed to detect processing errors. 
	Errors that required manual correction were recorded for each stage and at both target resolutions, and are summarized in Table \tabref{Tab:timings}.
	Benchmark results are reported for the execution of the workflow on a compute workstation equipped with an AMD Ryzen Threadripper PRO5965wx processor,
	using 16 CPU cores. 
	
\begin{figure}[!t]
		\centering        \includegraphics[width=0.9\textwidth]{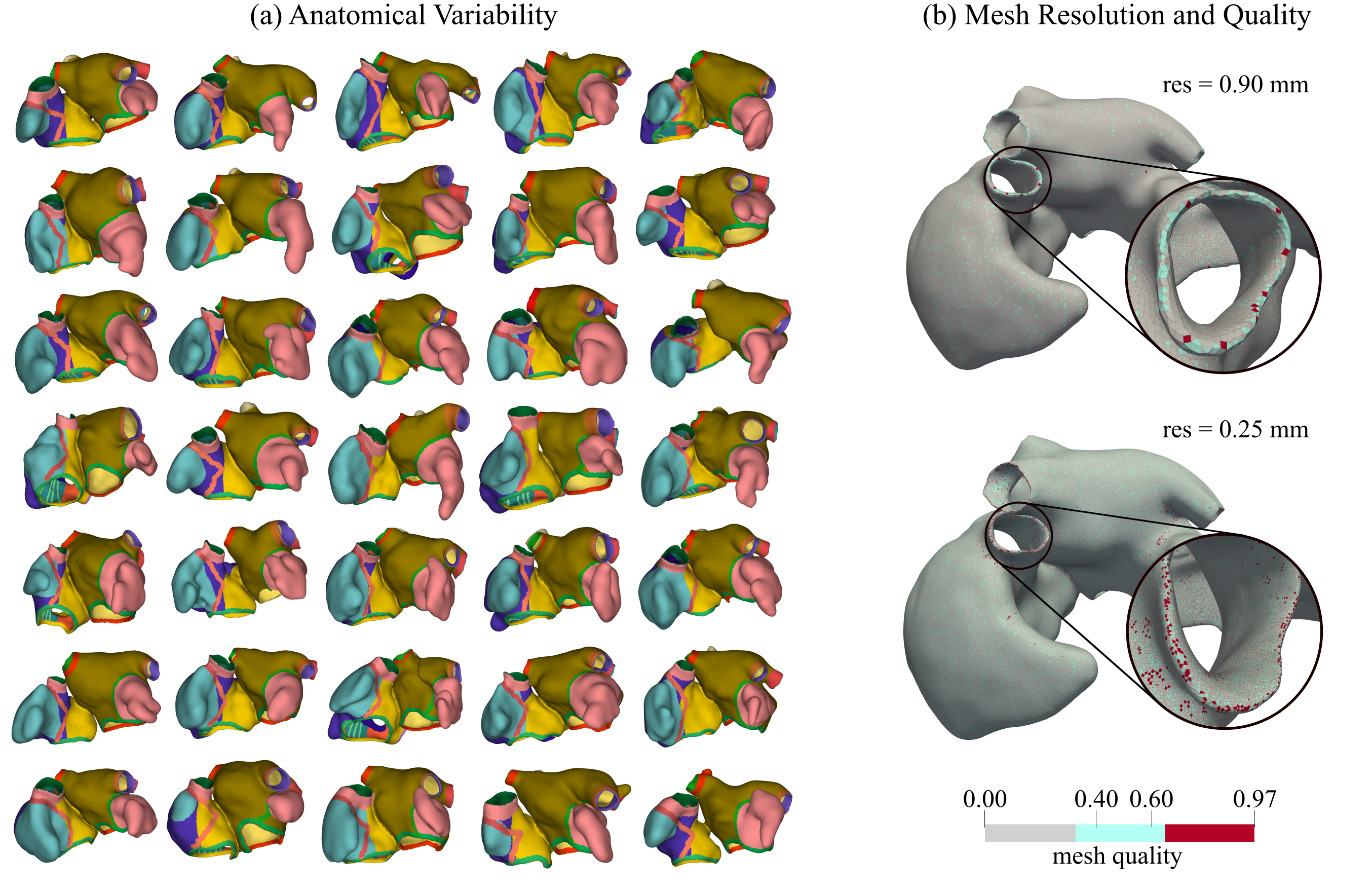}
		\caption{(a) Representation of 35 out of 50 generated biatrial models. Only a subset of the generated geometries is represented to allow for a better visualization. (b) Examples of computed meshes of resolution of $\approx$ \SI{0.90}{\milli\meter} and $\approx$ \SI{0.25}{\milli\meter} and corresponding obtained mesh quality.}
		\label{fig:_processing_results}
	\end{figure}

	Overall, all 50 cases in both resolutions were processed automatically in the majority of cases (38), 
	with only minimal user intervention required in 22 cases.
	The workflow produced anatomically highly detailed computational meshes, 
	with fine-grained domain annotation, fiber arrangement, and an anatomical reference frame (referred to Figure \ref{fig:_processing_results}).
	For the lower resolution \gls{re} model, the overall processing per model lasted, on average, less than ten minutes. 
	Higher resolution \gls{rd} compatible meshes were more costly to generate, specifically, the \gls{uac} generation stage where more than $\approx$84\% of the costs incurred.
	Importantly, a fully automatic meshing processes was achieved in all cases, yielding meshes free of topological errors and of, overall, excellent mesh quality. Worst element quality was always below 0.99, according to the quality metric \citep{karabelas2018:_towards},  
	which is considered a critical threshold in simulations using \gls{rd} solvers such as openCARP. Lower quality elements close to $0.99$ clustered mostly in high curvature regions, around the orifices of the thin-walled atria (referred to Figure \ref{fig:_processing_results}). 
	
	\subsubsection{Automation failure}  
	A manual intervention was needed to separate the \gls{lspv} from the \gls{laa} in 19 cases, and the \gls{cs} from the \gls{la} in 2 cases. The incorrect segmentation was however not due to the \gls{scn}, but to the low contrast in the acquired images. 
	At the second stage, during the labeling of the blood pools, manual intervention was required to switch the \gls{cs} with the \gls{ivc} landmarks in 21 cases.
	The final marking of the veins and the corresponding tissue to discard from the blood pool, including the \gls{cs}, \gls{laa}, and \gls{raa}, was automatically done by the workflow in 28 cases. For the remaining 22 cases, manual intervention was only needed to define or adjust the tissue to discard from the \gls{ivc} to ensure a better opening.

	\begin{table}[!t]
		\centering
		\begin{tabular}{lccc}
			\toprule
			& \multicolumn{2}{c}{\textbf{Mesh resolution/Model}} & \textbf{Manual Correction} \\
			&\textbf{\SI{0.9}{\milli\meter}}/\gls{re} & \textbf{\SI{0.25}{\milli\meter}}/\gls{rd} & \# Cases \\ 
			\hline \\[-2ex]
			Automatic multilabel segmentation &$\approx$ \SI[parse-numbers = false]{15}{\second} & $\approx$ \SI[parse-numbers = false]{15}{\second} & 19 \\
			Label augmentation on the blood pools  &\SI[parse-numbers = false]{297}{} $\pm$ \SI[parse-numbers = false]{180}{\second} &\SI[parse-numbers = false]{297}{} $\pm$ \SI[parse-numbers = false]{180}{\second} & 22 \\
			Walls extrusion and meshing  &\SI[parse-numbers = false]{187}{} $\pm$ \SI[parse-numbers = false]{172}{\second} &\SI[parse-numbers = false]{628}{} $\pm$ \SI[parse-numbers = false]{557}{\second} & 0\\
			Selection of the atrial orifices  &\SI[parse-numbers = false]{4}{} $\pm$ \SI[parse-numbers = false]{6}{\second} &\SI[parse-numbers = false]{39}{} $\pm$ \SI[parse-numbers = false]{36}{\second} & 0\\
			Anatomical structures and fibers  &\SI[parse-numbers = false]{30}{} $\pm$ \SI[parse-numbers = false]{35}{\second} &\SI[parse-numbers = false]{1084}{} $\pm$ \SI[parse-numbers = false]{628}{\second} & 0\\
			Universal atrial coordinates  &\SI[parse-numbers = false]{8}{} $\pm$ \SI[parse-numbers = false]{14}{\second} &\SI[parse-numbers = false]{175}{} $\pm$ \SI[parse-numbers = false]{37}{\minute} & 0 \\
			\bottomrule
			\\[-2ex]
			Total processing time & \SI[parse-numbers = false]{9}{} $\pm$ \SI[parse-numbers = false]{7}{\minute} & \SI[parse-numbers = false]{209}{} $\pm$ \SI[parse-numbers = false]{51}{\minute} & --  
		\end{tabular}
		\caption{Average timings for the generation of all the 50 biatrial geometries 
			at target resolutions of $\approx$ \SI{0.90}{\milli\meter} and $\approx$ \SI{0.25}{\milli\meter}
			and the number of models requiring manual corrections are given for each processing stage.}
		\label{Tab:timings}
	\end{table}

	\subsection{EP simulation workflow performance evaluation}
	The efficiency of our workflow in the function twinning stage, that is of simulating atrial \gls{ep} for high fidelity \gls{ecg} generation, was tested for both \gls{re} and \gls{rd} models for all 50 patients.
	Execution times of individual processing stages -- comprising torso generation, setting up of \glspl{ic}, computation of the \gls{ecg} lead fields, and the simulation of an entire atrial activation sequence initiated at the \gls{san} -- were measured (refer to Table \tabref{Tab.Times_re_rd} and Figure \ref{Fig.RD_vs_RE}).
	As shown previously in detail \citep{gillette2021:_framework}, owing to its relaxed mesh resolution dependency, 
	the \gls{re} model is significantly more lightweight,
	facilitating the setup of an \gls{ep} model in $\approx$ \SI{3}{\minute} and the computation of a full biatrial activation sequence with a high fidelity \gls{ecg} in $\approx$ \SI{1}{\second}.
	All steps of the workflow are computationally more costly due to the stricter mesh resolution requirements, requiring $\approx$ \SI{32}{\minute} for setup, and $\approx$ \SI{19}{\minute} 
	for computing the activation sequence and the \glspl{ecg}.
	While formally the same torso surface was employed, extra costs were incurred due to the required mesh conformity at the atrial surfaces,
	leading to a markedly higher number of elements in the \gls{rd} case. 
	
	
	
	\begin{table}[!t]
		\centering
		\begin{tabular}{lll}
			\toprule
			& \multicolumn{2}{c}{\textbf{Mesh resolution / Model}} \\
			&\SI{0.9}{\milli\meter} / \gls{re} &\SI{0.25}{\milli\meter} / \gls{rd}\\ 
			\hline \\[-2ex]
			Volumetric torso &$\approx$ \SI[parse-numbers = false]{2}{\minute} \SI[parse-numbers = false]{25}{\second} & $\approx$ \SI[parse-numbers = false]{22}{\minute} \SI[parse-numbers = false]{30}{\second}\\
			Interatrial connections  &$\approx$  \SI[parse-numbers = false]{19}{\second} & $\approx$ \SI[parse-numbers = false]{39}{\second}\\
			Lead field &$\approx$ \SI[parse-numbers = false]{35}{\second} &$\approx$ \SI[parse-numbers = false]{9}{\minute} \\
			\gls{ep} simulation  &$\approx$  \SI[parse-numbers = false]{27}{\second} &$\approx$ \SI[parse-numbers = false]{19}{\minute}\\
			\bottomrule
		\end{tabular}
		\caption{Execution times of setup and simulation using \gls{re} or \gls{rd} models.}
		\label{Tab.Times_re_rd}
	\end{table}
	
	The computational costs of creating the \gls{ic} cables also showed slight differences between the coarse \SI[parse-numbers = false]{0.9}{\milli\meter} mesh and the fine \SI[parse-numbers = false]{0.25}{\milli\meter} mesh, averaging \SI[parse-numbers = false]{19}{\second} and \SI[parse-numbers = false]{39}{\second}, respectively. However, this difference is negligible compared to the total costs of generating the biatrial and torso geometry, as well as the CPU time required for the computation of the lead field solution, which is of \SI[parse-numbers = false]{35}{\second} for the \gls{re} mesh, and \SI[parse-numbers = false]{9}{\minute} for the \gls{rd} grid, and the  \gls{ep} simulation, corresponding to \SI[parse-numbers = false]{27}{\second} of computation for the \gls{re} model, and to \SI[parse-numbers = false]{19}{\minute} for the \gls{rd} gold-standard monodomain model.
	
	\subsection{Assessing discrepancy between \gls{re} and \gls{rd} model}
\begin{figure}[!t]
		\centering \includegraphics[width=0.95\textwidth]{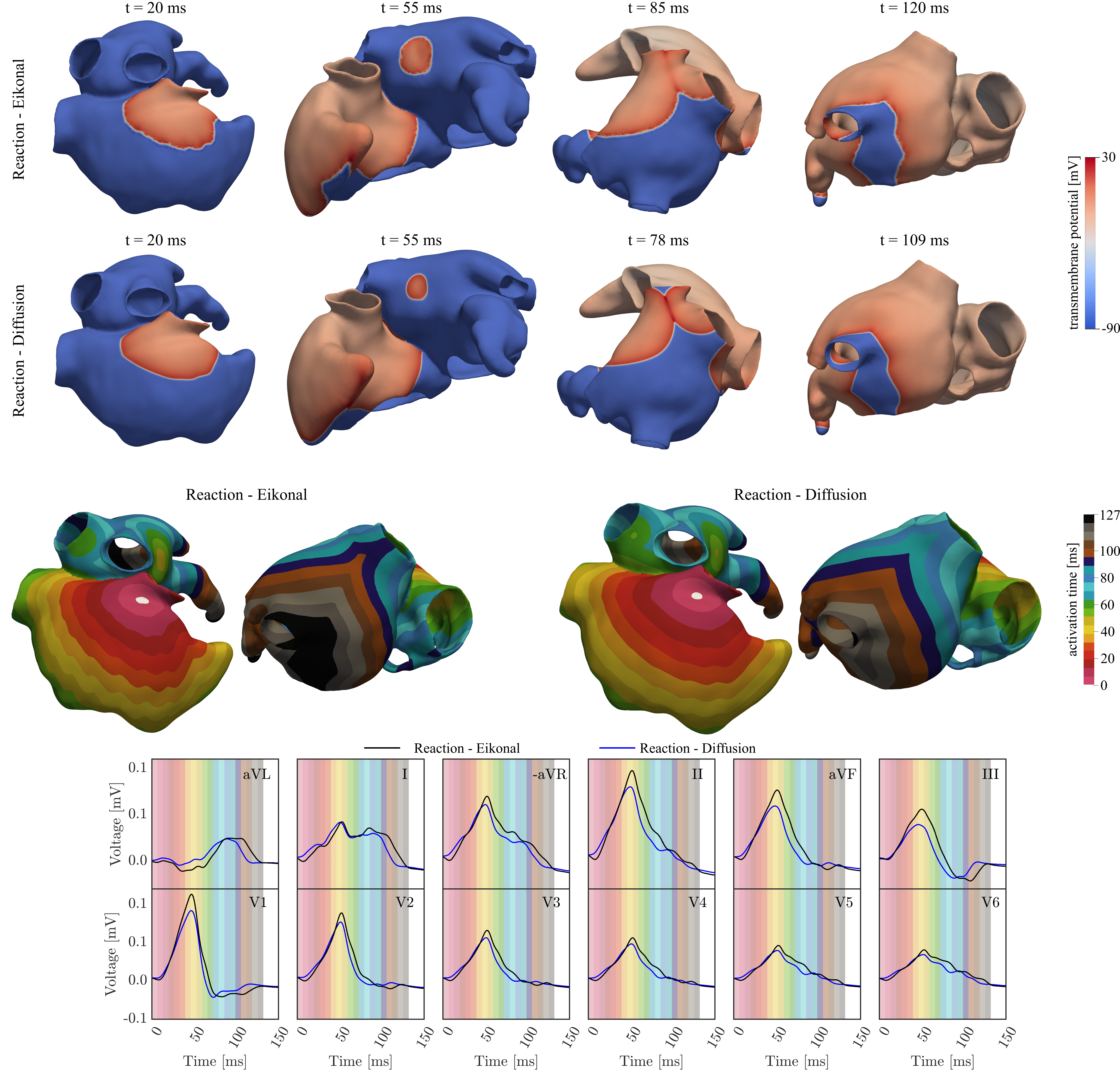}
		\caption{First/Second row: Transmembrane potential $V_m$ propagation computed with the \gls{re} and \gls{rd} models, at different time instants. Third row: Activation map obtained by solving the \gls{re} and the \gls{rd} model. Fourth row: \glspl{ecg} obtained by solving the \gls{re} (black) and \gls{rd} (blue) models coupled with the lead field. The background colors refer to the activation map.  By coloring the background with time-bands of colors corresponding to the activation map, the activation map is compared with the \gls{ecg}, highlighting the relation between regional atrial activation and P-wave.}
		\label{Fig.RD_vs_RE}
	\end{figure}

	To evaluate the computational fidelity of the \gls{re} model compared to a gold standard \gls{rd} monodomain model \citep{nagel2022comparison}, 
	we compare the spatiotemporal distribution of electrical sources, $V_{\rm m}(\boldsymbol{x},t)$, 
	the activation maps, $\tau(\boldsymbol{x})$, and the P-wave generated by both models,
	employing the baseline parameter settings in Section \ref{Subsubsec:baseline_sim}. 
	Results are illustrated for a representative test case in Figure \ref{Fig.RD_vs_RE}.

	Overall, the model parameters calibrated to match \gls{rd} and \gls{re} conduction velocities led to nearly identical activation patterns. 
	Minor differences emerged due to imperfect matching of conduction velocities in the \gls{ra} and \gls{la} tissue,
	as well as in the \gls{ic} cables,
	and differences in anterograde and retrograde activation by the \glspl{ic} cables mediated by an electrotonic source-sink mismatch at the interface between cable and tissue (refer to Figure \ref{Fig.RD_vs_RE}, top panels). 
	These combined effects led to a slightly longer total activation time of the atria,
	which manifested in minor differences in the \gls{pwd}, 
	with an average variation of \SI[parse-numbers = false]{2}{\milli\second}.
	
	Differences in P-wave morphology and magnitude were negligible,
	with a maximum amplitude difference in lead II of 5\%, and an average \gls{rmse} across all leads of 3.71\% (refer to Table \ref{Tab.RE_vs_RD}).
	Magnitude differences could be attributed to the difference in spatial resolution,
	as these disappeared when running the \gls{re} simulation on the higher resolution mesh (not shown).
	
	\begin{table}[!t]
		\centering
		\begin{tabular}{lccccccccccccc}
			\toprule
			Lead &aVL &I &-aVR &II &aVF &III &V1 &V2 &V3 &V4 &V5 &V6 &Average\\ 
			\hline \\[-2ex]
			\gls{rmse} [\%] &3.56 &3.44 &3.54 &5.16 &4.87 &5.17 &4.64 &3.66 &3.09 &1.69 &2.39 &2.34 &3.71\\
			\gls{pwd} \gls{re} [ms] &128 &131 &131 &125 &129 &89 &126 &131 &123 &126 &131 &131 &125.1\\
			\gls{pwd} \gls{rd} [ms] &121 &124 &121 &124 &129 &90 &120 &131 &123 &126 &131 &131 &122.6\\
			\bottomrule
		\end{tabular}
		\caption{\gls{rmse} in percentage, between \gls{ecg} signals computed with the \gls{re} and \gls{rd} models coupled with the lead field. The \gls{pwd} of the \gls{ecg} traces obtained with the \gls{re} and the \gls{rd} models coupled with the lead field are also reported. }
		\label{Tab.RE_vs_RD}
	\end{table}

	\subsection{Impact of Bachmann's Bundle insertion on the P-wave}
\begin{figure}[!t]
		\centering \includegraphics[width=0.95\textwidth]{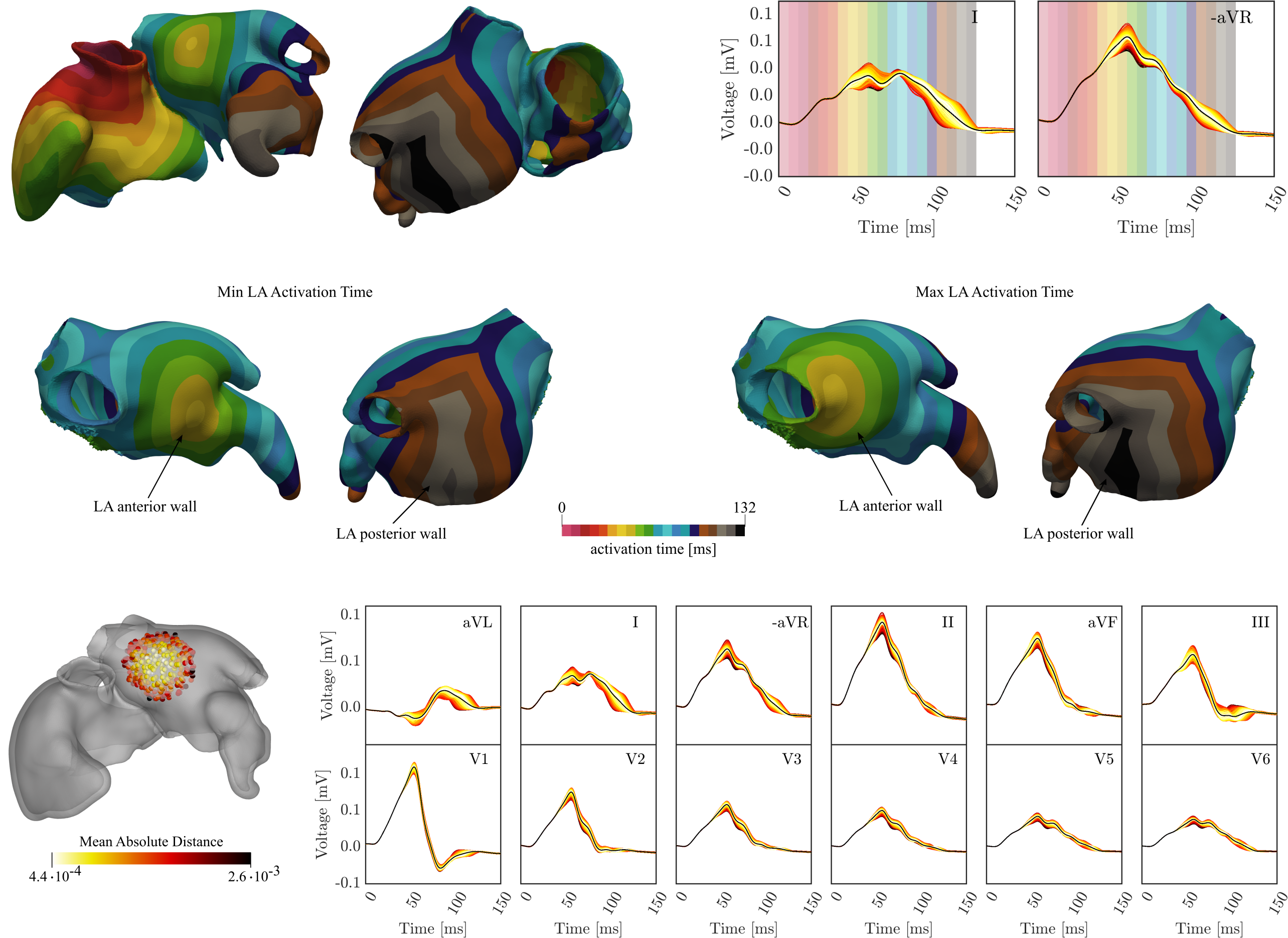}
		\caption{Influence of \gls{bb} \gls{la} insertion site.  First row: 
			Biatrial activation map obtained for \gls{bb} reference insertion site along with P-waves in leads I and -aVR. The background of P-wave traces is color coded according to the activation map 
			to highlight the relation between 
			the dipole layer of depolarization wave fronts 
			and the associated amplitude in the P-wave. 
			Second row: \gls{la} activation maps corresponding to the minimum and maximum total activation time. 
			Third row: P-wave variation due to varying the \gls{bb} insertion site in the \gls{la}. 
			Insertion site an \gls{ecg} are color coded corresponding to the absolute distance from the mean \gls{ecg} (black), 
			averaged over time and the 12 leads.}
		\label{Fig.Tests_1BB}
	\end{figure}

	The role of \gls{bb} insertion site $\mathbf{x}_{BB_{LA}}$ on the \gls{la} activation and P-wave morphology was investigated.
	The use of a cable-based \gls{ic} formulation readily facilitated a parametric sweeping of $\mathbf{x}_{BB_{LA}}$,
	where $\mathbf{x}_{BB_{LA}}$ was varied within a radius of $\approx$ \SI[parse-numbers = false]{5}{\milli \meter} 
	around the reference location, $\hat{\mathbf{x}}_{BB_{LA}}$, as used above.
	Shifts in $\mathbf{x}_{BB_{LA}}$ altered the location and timing of wavefronts collisions with waves initiated through the other \glspl{ic}.
	These were predominantly the rim of the fossa ovalis, activating the \gls{la} almost synchronously with the \gls{bb},
	and the posterior-superior \gls{ic}, activating the carina of the \glspl{lpv} with a delay of \SI[parse-numbers = false]{10}{\milli \second}. 
	This entailed a change in the total activation time of the \gls{la}, and induced a spread in \gls{pwd} and magnitude of the P-wave, with very limited effect on the corresponding morphology (refer to Figure \ref{Fig.Tests_1BB}). 
	
	Variation of $\mathbf{x}_{BB_{LA}}$ alone led to differences in total activation time of \SI[parse-numbers = false]{23}{\milli \second}, 
	in \gls{pwd} between \SI{12}{\milli \second} to \SI{24}{\milli \second} in individual leads, with a minimum and maximum average \gls{pwd} of 104 and \SI[parse-numbers = false]{120}{\milli \second} across all leads, respectively,    
	and in magnitude to a mean absolute distance \eqref{eq:MAD} between 4.4$\times$10$^{-4}$ and  2.6$\times$10$^{-3}$ mV.
	A quantitative summary on \gls{pwd} variability is given in Table \ref{Tab.PWD}.

	\begin{table}[!t]
		\centering
		\begin{tabular}{lccccccccccccc}
			\toprule
			Lead &aVL &I &-aVR &II &aVF &III &V1 &V2 &V3 &V4 &V5 &V6 &Average\\ 
			\hline \\[-2ex]
			Min \gls{pwd} [ms] &93 &93 &108 &110 &104 &107 &99 &106 &106 &106 &106 &106 &104,2\\
			Max \gls{pwd} [ms] &117 &115 &124 &126 &119 &120 &118 &118 &118 &124 &124 &123 &120,5\\
			(Max - Min) \gls{pwd} [ms] &22 &16 &13 &16 &24 &15 &19 &12 &12 &18 &15 &14 &16,3\\
			\bottomrule
		\end{tabular}
		\caption{Maximum and minimum \gls{pwd}, and their differences, for each lead and averaged over all leads.}
		\label{Tab.PWD}
	\end{table}
	
	\subsection{Effect of a fan-like \gls{la} insertion Bachmann's Bundle on the P-wave}
\begin{figure}[!t]
		\centering \includegraphics[width=0.9\textwidth]{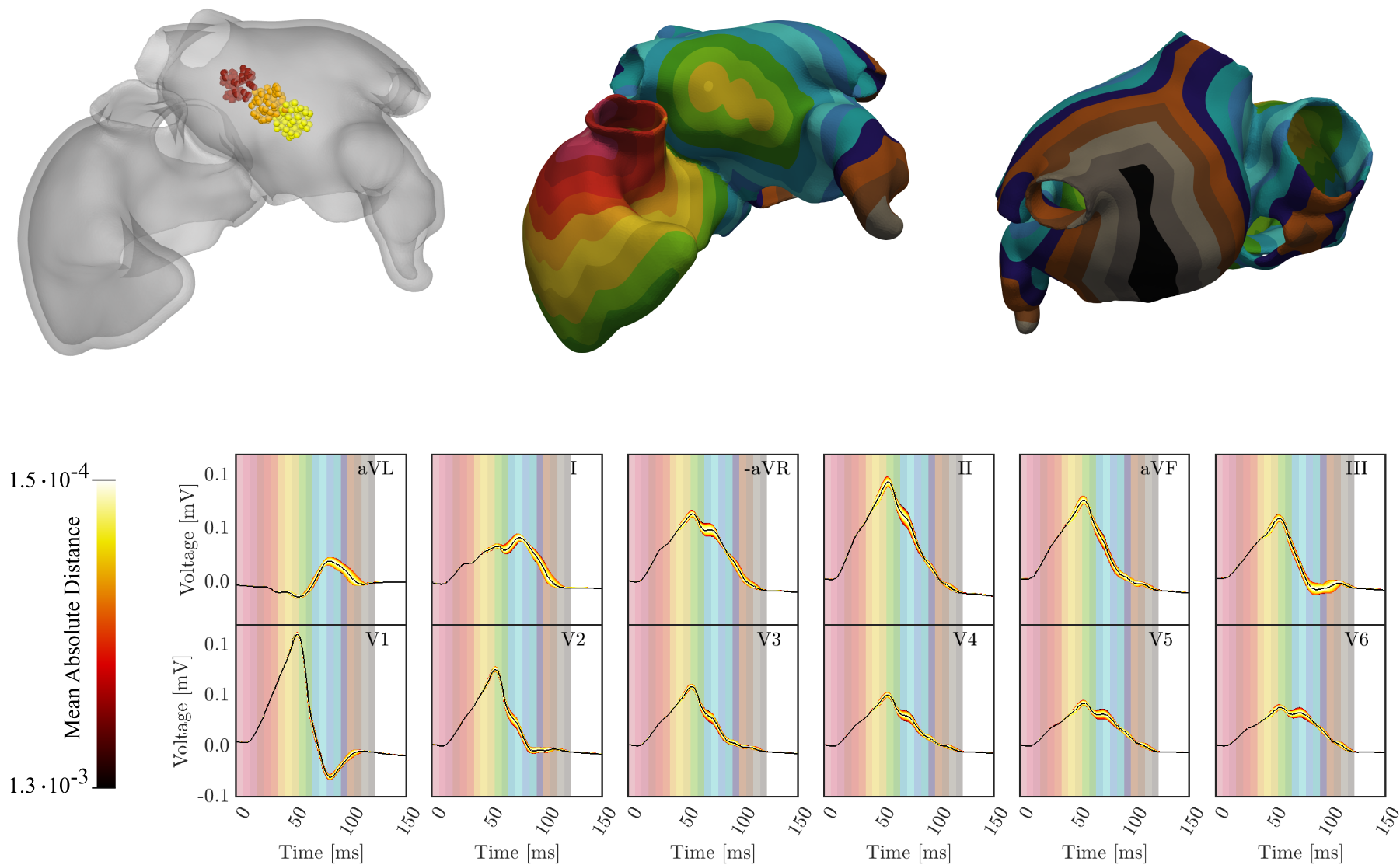}
		\caption{First row-left: Investigated \gls{la} entry site areas of the three cables representing the \gls{bb}. The different colors identify the three different regions of search. First row-center and right: An example of activation map obtained with the three-cable \gls{bb}. Second row: P-wave variation obtained by moving the \gls{la} entry sites of the three cables representing the \gls{bb}. The signals coloration correspond to the absolute distance from the mean \gls{ecg}, averaged over time and the leads. The mean \gls{ecg} is shown in black. By coloring the background of the lead with time-bands of colors corresponding to the activation map, the activation map is compared with the \gls{ecg}, highlighting the relation between regional atrial activation and P-wave.}
		\label{Fig.Three_bb}
	\end{figure}

	Our modeling approach flexibly supports an arbitrary number of \glspl{ic} cables which can be bundled to increase source strength for generating \glspl{ecg} and \glspl{egm},
	or to model inter-atrial coupling with higher anatomical complexity.
	This ability is showcased for investigating the influence of a distributed fan-like insertion of \gls{bb} into the \gls{la} with three independent coupling locations.
	An exemplary activation map along with the resulting distribution of P-waves is shown in Figure \ref{Fig.Three_bb}.
	As the same velocity was assumed in all cables the initial activation at the three \gls{la} insertion sites 
	was not synchronous but led to an elongated area of early activation along the location of the wide \gls{bb}.
	
	The total activation time was observed between 118 and \SI[parse-numbers = false]{125}{\milli\second}, comparable smaller than the \gls{bb} single-cable case. 
	
	The highest variation in the P-wave was observed between 80 and \SI[parse-numbers = false]{132}{\milli \second}, corresponding to the activation of the posterior wall of the \gls{la}.
	The average range of variation of the \gls{pwd} was registered from 102 to \SI[parse-numbers = false]{110}{\milli \second}, 
	thus smaller than the single-cable case. 
	The maximum and minimum \gls{pwd} over the leads was also reduced compared to the single-cable case, with maximum variation in lead aVL and II.    
	
	\begin{table}[!t]
		\centering
		\begin{tabular}{lccccccccccccc}
			\toprule
			Lead &aVL &I &-aVR &II &aVF &III &V1 &V2 &V3 &V4 &V5 &V6 &Average\\ 
			\hline \\[-2ex]
			Min \gls{pwd} [ms] &80 &100 &107 &110 &95 &107 &94 &98 &110 &110 &109 &109 &102,4\\
			Max \gls{pwd} [ms] &89 &114 &115 &118 &99 &115 &104 &108 &113 &113 &115 &116 &109,9\\
			(Max - Min) \gls{pwd} [ms] &9 &14 &8 &8 &4 &8 &10 &10 &3 &3 &6 &7 &7,5\\
			\bottomrule
		\end{tabular}
		\caption{Maximum and minimum \gls{pwd}, and their differences, for each lead and averaged over all leads.}
		\label{Tab.PWD_threeBB}
	\end{table}
	
	
	\subsection{Effect of the RA endocardium on the P-wave}
	\begin{figure}[!t]
		\centering \includegraphics[width=0.95\textwidth]{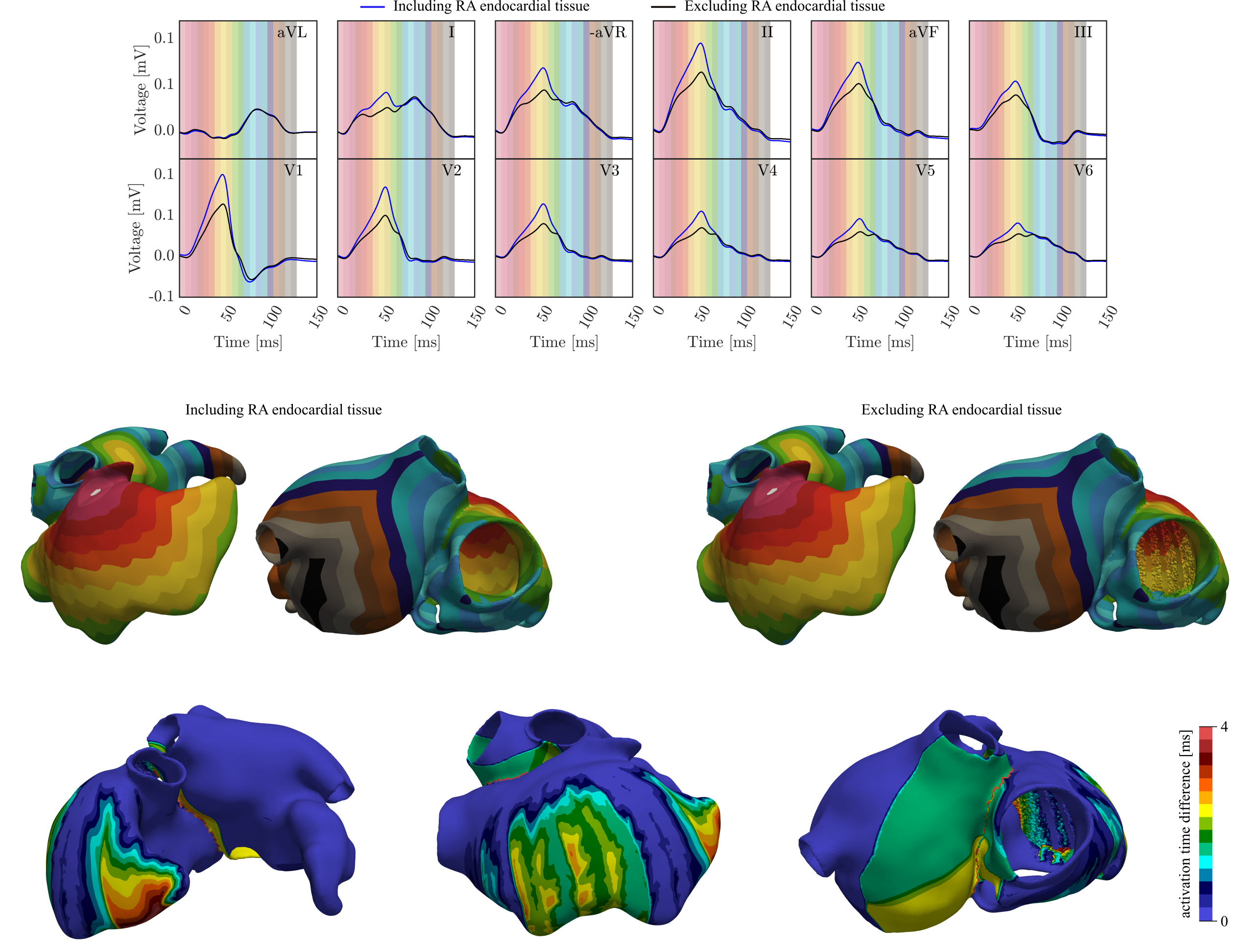}
		\caption{First row: P-wave obtained when including (blue) and excluding (black) the \gls{ra} endocardial tissue outside of the \gls{pm}s and \gls{ct}. The P-wave of all leads is compared to the activation map by coloring the background. Second row: Activation maps obtained when including (left) and excluding (right) the \gls{ra} endocardial tissue. Third row: Absolute difference between the activation maps obtained when including and excluding the \gls{ra} endocardial tissue.}
		\label{Fig.Tests_PM}
	\end{figure}
	
	The \glspl{ecg} obtained by simulating the atria, both including and excluding the \glspl{ra} endocardial layer outside of the \gls{pm} and \gls{ct}, are shown in Figure \ref{Fig.Tests_PM}. A decrease in wave amplitude between 14 and \SI[parse-numbers = false]{71}{\milli \second}, corresponding to the activation of the \gls{ra}, was observed in all leads. The \gls{rmse} in percentage is reported in Table \ref{Tab.PWD_PM}. The major amplitude variation, exceeding 1\%, was observed in leads II, V1, and V2, which capture the signal propagation towards the \gls{raa} and the \gls{pm}. Variations between 0.49\% and 0.91\% were observed in all other leads, except for aVL, where the variation was almost null. The \gls{pwd}, as reported in Table \ref{Tab.PWD_PM}, remained unchanged in all leads.
	
	\begin{table}[!t]
		\centering
		\begin{tabular}{lccccccccccccc}
			\toprule
			Lead &aVL &I &-aVR &II &aVF &III &V1 &V2 &V3 &V4 &V5 &V6 &Average\\ 
			\hline \\[-2ex]
			\gls{rmse} [\%] &0,08 &0,63 &0,91 &1,19 &0,88 &0,58 &1,31 &1,18 &0,85 &0,72 &0,58 &0,49 &0,78\\
			PWD [ms] &114 &119 &120 &121 &97 &89 &111 &86 &88 &99 &101 &121 &105,5\\
			\bottomrule
		\end{tabular}
		\caption{\gls{rmse} between the \gls{ecg}s obtained when excluding the \gls{ra} endocardial tissue outside of the \gls{pm}s and \gls{ct}, and the \gls{ecg}s signal obtained including the \gls{ra} endocardial tissue, and \gls{pwd}, for each lead and averaged over all leads.}
		\label{Tab.PWD_PM}
	\end{table}
	
	The \gls{ecg} variations were related to the computed activation maps and their absolute difference, depicted in Figure \ref{Fig.Tests_PM}. When the \gls{ra} endocardial wall was included, the signal propagated faster on the posterior part of the \gls{raa}, on the posterior-lateral wall of the \gls{ra} near the \glspl{pm}, and on the septum. A total activation delay of about \SI[parse-numbers = false]{2.1}{\milli \second} was observed at the tip of the \gls{raa}, and \SI[parse-numbers = false]{2.7}{\milli \second} on the epicardium between the \glspl{pm}. The \gls{bb}, \gls{ra} anterior wall, and the \gls{ct} were activated at the same speed. However, the removal of the \gls{ra} endocardial layer between the \gls{ct} and the septum caused a delay in the activation of the septum and the \gls{ra} wall near the \gls{cs}. Although this delay did not affect the overall \gls{ra} activation, it slightly impacted the activation of the posterior wall of the \gls{la}, causing a total delay of up to \SI[parse-numbers = false]{2}{\milli \second} in the posterior-inferior region. Nonetheless, the total activation map of the \gls{la} remained almost unchanged, consistent with the computed \gls{ecg}s.
	
	\subsection{Calibration of the atrial activation sequence based on the P-wave of the 12 lead \gls{ecg}}
	
\begin{figure}[!t]
		\centering \includegraphics[width=0.95\textwidth]{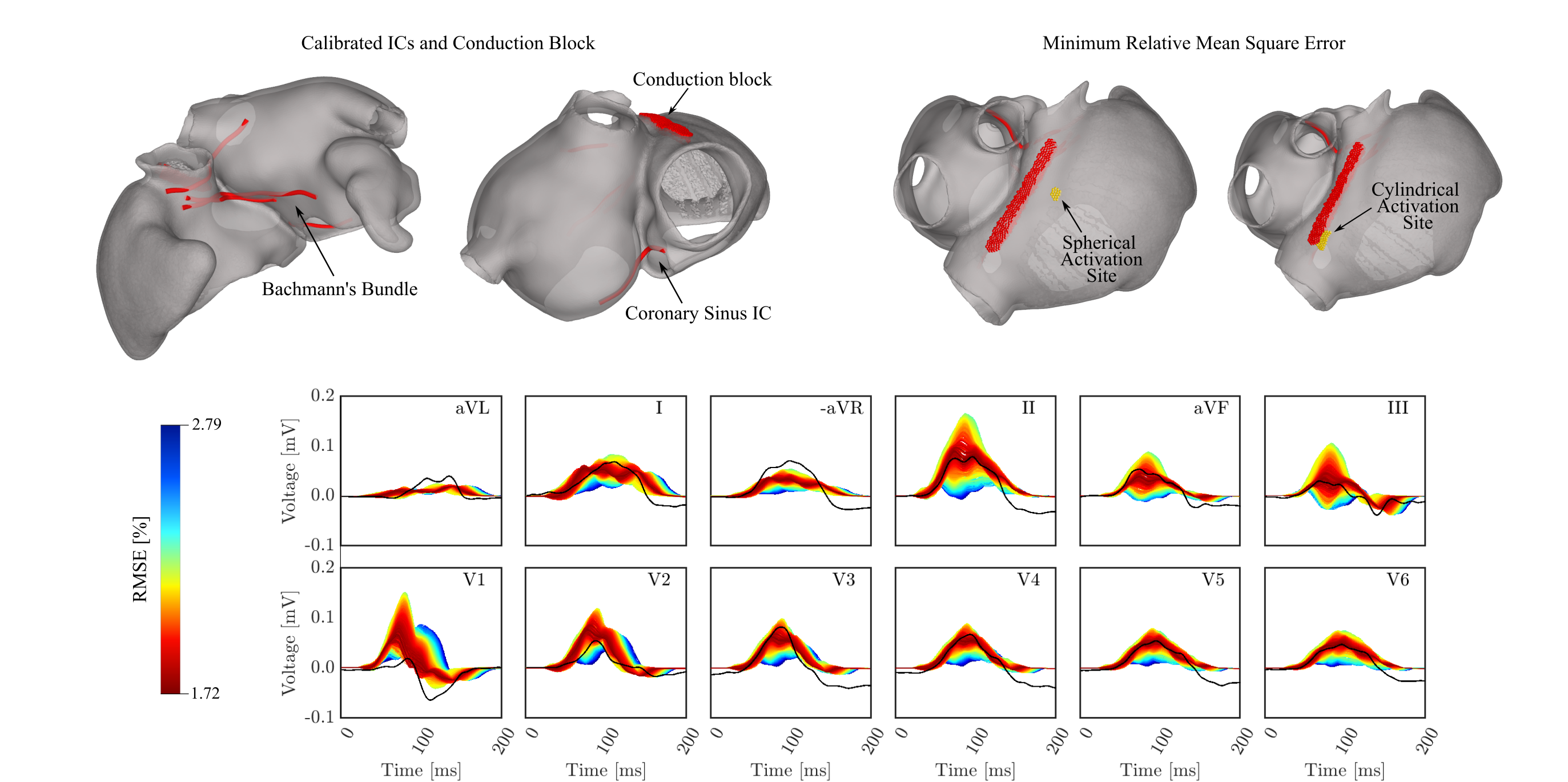}
		\caption{P-wave based calibration of the atrial activation sequence. 
			Top row, left: Constant \glspl{ic} configuration 
			obtained by interactive pseudo-calibration. 
			Top row - right: Spatial sampling ranges of \gls{san} exit sites 
			with optimal location that minimized the mismatched as measured by the \gls{rmse}. 
			Bottom row: Physiological envelope of P-waves (colored traces) obtained by varying 
			\gls{san} exit site, and the velocities $v_f$ and $v_s$ in the \gls{ra} 
			together the clinical P-wave (black).
			Traces are color-coded according to their distance to the clinical P-wave
			measured as \gls{rmse}.
		} \label{Fig.Fig_SAN_sampling}
	\end{figure}

	One anatomical model of a patient was selected along with the recorded P-wave under normal sinus rhythm. 
	To keep the calibration procedure tractable, the parameter space to be explored was restricted to the initial exit site at the \gls{san}, the conduction velocity in the \gls{ra},
	and the location and timing of three insertion sites of \gls{bb} in the \gls{la}. 
	In a first pseudo-calibration step for a fixed \gls{san} exit site and baseline velocities in the \gls{ra}, 
	the location of \gls{bb} was varied in interactive simulation runs to obtain a close approximation of the terminal half of the P-wave. 
	One configuration yielding a good morphological fit in all leads, assessed only by visual inspection, was selected (see Figure \ref{Fig.Fig_SAN_sampling}, top row).

	
	Keeping the \gls{bb} fixed, the physiological envelope of the P-wave signals was computed by sampling over $\mathbf{x}_{SAN}$, $v_f$, and $v_s$ (see Figure \ref{Fig.Fig_SAN_sampling}, mid row). 
	The P-wave envelope under this sampling covered the observed clinical signal in most leads quite well, 
	with the exceptions of leads aVL, -aVR and V1. 
	While the discrepancy in the low amplitude P-waves in leads aVL and -aVR were rather small,
	this was not the case in lead V1 where the biphasic P-wave could not be approximated.
	The terminal negative deflection of the measured P-wave was likely mediated by fibrotic tissue in this AF patients with impaired electrophysiologicalgls{ep} excitability and repolarization properties.
	
	Quantitatively, the P-waves were well approximated, 
	with a \gls{rmse} in the range of $[1.72\%,2.79\%]$.
	The minimum \gls{rmse} was achieved 
	when the \gls{san} was positioned at 62\% along the \gls{ivc}-\gls{svc} axis and 35\% along the lateral-to-septal axis (see Figure \ref{Fig.Fig_SAN_sampling}, top row), 
	with conduction velocities of $v_f = \SI[parse-numbers = false]{1.1}{\meter/\second}$ and $v_s = \SI[parse-numbers = false]{0.45}{\meter/\second}$. 
	The \gls{rmse} values, as well as the \gls{pwd} for each lead, 
	obtained with the optimal parameter set that minimized \gls{rmse} are presented in Table \ref{Tab.PWD_SAN}.
	On average, simulated \gls{pwd} was longer than clinical \gls{pwd}, with \SI[parse-numbers = false]{139}{\milli\second} versus \SI[parse-numbers = false]{130}{\milli\second}, suggesting that calibration could be further improved 
	by including \gls{ic} and \gls{la} velocities in the sampled parameter space.

	\begin{table}[!t]
		\centering
		\begin{tabular}{lccccccccccccc}
			\toprule
			Lead &aVL &I &-aVR &II &aVF &III &V1 &V2 &V3 &V4 &V5 &V6 &Average\\ 
			\hline \\[-2ex]
			min \gls{rmse} [\%] &1.16 &1.68 &2.27 &2.06 &1.33 &1.02 &2.60 &1.51 &2.22 &1.87 &1.58 &1.37 &1.72 \\[1ex]
			\gls{pwd} [ms]  &&&&&&&&&&&&& \\
			Clinical data &139 & 147 & 123& 124 &127 &146 &134 & 125 &118 &132 &120 &129 &130.33\\
			Simulated \gls{ecg} &136 &137 &135 &150 &133 &151 &145 &142 &131 &129 &131 &151 &139.25 \\
			Diff with data &3 & 10 &12 &26 &6 &5 &11 &17 &13 &3 &11 &22 &8.92\\[1ex]
			\bottomrule
		\end{tabular}
		\caption{Minimum \gls{rmse}, in percentage, \gls{pwd} of the \gls{ecg} data, of the computed \gls{ecg} with minimum \gls{rmse}, and their difference, obtained by varying the position of the spherical or cylindrical activation impulse representing the \gls{san}, and the conduction velocity $v_f$ and $v_s$ of the \gls{ra}.}
		\label{Tab.PWD_SAN}
	\end{table}

	\section{Discussion}
	\label{Sec:discussion}
	In this study, we present a highly automated scalable end-to-end workflow for generating patient-specific anatomical digital twins of human atria for their efficient calibration based on clinical data, e.g. the \glspl{ecg}.
	Our novel workflow is comprehensive, 
	including the generation of an anatomical reference frame, that facilitates parametric encoding of all space-varying \gls{ep} model properties,
	as required for model calibration, 
	a detailed annotation of all relevant anatomical landmarks and structures,
	including a fiber architecture based on a detailed anatomy-informed set of rules, and a new method for lightweight and flexible \glspl{ic}.
	
	Our workflow, implemented in a single software building on meshtool \citep{neic2020automating},
	generates smooth biatrial anatomical representations of sufficient mesh quality 
	for a range of resolutions suitable for the entire spectrum of atrial \gls{ep} simulations.
	Both coarser meshing at $\approx$ \SI{0.9}{\milli \meter} for real-time calibration using \gls{re} type models
	\citep{neic17:_reaction_eikonal,pezzuto2017:_ecg},
	as well as finer meshing with sufficient resolution for fully mechanistic \gls{rd} modeling studies, is supported. 
	The workflow efficiently generates simulation-ready models within less than 10 minutes
	at a reference resolution of approximately \SI[parse-numbers = false]{0.90}{\milli \meter}.
	In contrast to previous approaches for building volumetric biatrial models,
	that used statistical shape models to generate virtual anatomies \citep{nagel2021:_biatrial_ssm,roney2023:_bia_vol},
	our approach generates accurate representations of an individual patient's atrial anatomy, 
	with fidelity being limited only by image quality and segmentation accuracy.
	A key advantage of our approach is its robustness, 
	as the description of anatomy is built on a volumetric image stack.
	This avoids the more common manifold extrusion procedure prone to topological mesh errors.
	
	Building on this anatomical model generation pipeline, the workflow is moreover extended for \gls{ecg} based model calibration, employing either  
	\gls{re}-based real-time \gls{ep} modeling, suitable for \gls{ep} calibration,
	or a full fidelity computationally more demanding \gls{rd} model, suitable for predictive simulations.
	Inter-atrial conduction is modeled by initially imposing electrical insulation between the \gls{ra} and \gls{la} \citep{costa2014:_fibrotic}, and then discretely reconnecting them using a novel, highly flexible, and physiologically constrained method. This approach employs auto-generated cables to establish inter-atrial conduction pathways, accurately representing all anatomically significant \glspl{ic} between the \gls{ra} and \gls{la}. Our approach facilitate the generation of pathways with prescribed        
	conduction velocities to connect arbitrary parametrically steerable locations within \gls{ra} and \gls{la} on the fly, 
	and, thus, avoids rigid and error-prone explicit meshing of inter-atrial bundles \citep{nagel2021:_biatrial_ssm}.
	
	Finally, we employ a forward \gls{ecg} generation framework based on the \gls{relf} model 
	that yields full fidelity \glspl{ecg} with real-time performance, 
	to support a clinically compatible model calibration procedure.
	To this end, volumetric atrial models are swiftly registered with reference torso models
	with corresponding electrodes, and integrated into volumetric conformal atria-torso meshes. 
	All combinations of cardiac \gls{ep} source and field models are supported.
	As mesh conformity is preserved, high fidelity \glspl{ecg} and \glspl{egm} can be computed
	over the entire solution domain based on either bidomain or lead field formulation, 
	without being restricted to lower fidelity methods such as potential recovery \citep{bishop2011bidomain},
	or boundary element approaches, that require smooth low pass pass-filtered coarse representations 
	of the atrial sources \citep{schuler2019:_spatial}. 
	
	The ability to compute full fidelity \glspl{ecg} with close to real-time performance is demonstrated by exploring important parameters governing the genesis of the P-wave, which were impossible or notoriously difficult to explore with previous approaches, e.g., the anatomy of \glspl{ic}.
	We show the equivalence in activation patterns and \glspl{ecg} between the lower resolution \gls{re} model, and the high resolution \gls{rd} model, with discrepancies well below the overall model uncertainty.
	These combined features make our end-to-end workflow suitable for large-scale atrial \gls{ep} modeling studies.
	The achieved model efficiency supports a fast exploration of the high-dimensional parameter space 
	spanned by atrial anatomy, structure and \gls{ep}, 
	by facilitating unattended sweeps over important space-varying parameters
	which is needed for the automation of optimization loops 
	using \gls{ecg} or \gls{egm} observations as target for calibration.

	\subsection{Scalable generation of volumetric biatrial anatomy models}
	Robust computational workflows capable of generating comprehensive volumetric biatrial anatomical models, with fiber architecture and anatomical reference frames at scale remain challenging.  
	Only a few approaches have been reported to date, exhibiting significant variation in the degree of automation, anatomical and structural fidelity, flexibility in prescribing fiber architecture, and their ability to generate clinically observable signals such as \glspl{egm} and \gls{ecg} with sufficient accuracy.
	
	The majority of atrial modeling studies using larger cohorts relied upon a bilayer formulation \citep{labarthe2014:_bilayer}
	which greatly simplifies the meshing procedure and reduces computational expenses
	at the cost of decreased biophysical fidelity. 
	For instance, it has been demonstrated that the behavior of bilayer models 
	may differ strikingly from full 3D models
	with respect to propagation patterns and arrhythmia dynamics \citep{roney2021constructing}.
	Discrepancies may stem from various factors such as e.g.\ altered source-sink relations,      
	the impact of mediating transmural conduction through resistive coupling between the endocardial and epicardial layers using finite elements of incompatible dimensionality,
	or differences in inter-atrial conduction pathways.
	Another important limitation of bilayer models is their reduced accuracy in representing the atrial sources, and a lack of suitable approaches for computing higher-fidelity extracellular potential fields \citep{bishop2011bidomain}, due to incompatibility with standard \gls{fem}.
	As such, modeling clinical data, such as the \glspl{egm} or \gls{ecg}, with high fidelity is challenging. Moreover, modeling extracellular potential field-related physiological events, such as bath loading effects, cannot be represented at all \citep{bishop2011bathloading}.
	These factors combined limit their trustworthiness and application scope.
	
	Volumetric biatrial meshes generated in early pioneering studies used artisanal handcrafted non-scalable methods that were tractable only for constructing small-size cohorts, 
	often consisting of a single model only \citep{harrild2000:_atria,seemann2006:_atria}.  
	Later, model-generation pipelines were proposed, gradually refined and increasingly automated
	towards a more industrialized process for producing models \emph{en masse}.
	A first mostly complete workflow was reported in a study by \cite{krueger2012personalization} 
	which was applied to generate a cohort of nine biatrial models.
	However, as numerous manual interventions were required in segmentation and landmark selection, for fiber mapping, and anatomical region definition, scalability remained limited.
	
	A further automated approach using a segmentation-derived endocardial surface mesh as input
	was reported and tested in 29 patients in \cite{azzolin2023:_augmenta}.
	The approach bears similarities to ours in that 
	e.g.\ surface curvature measures are used to identify structures such as the pulmonary veins, 
	but also differs markedly in other processing steps 
	which are more challenging to perform when using manifold meshes only. 
	For instance, the identification of the mitral valve orifice,
	which is trivially identified as a label interface in our volumetric approach,
	required a rigid registration with a mean statistical shape model.
	Moreover, important anatomical landmarks such as the \gls{fo} and its rim 
	or the \gls{cs} were not identified,
	and additional manual procedures were required to identify anatomical regions such as the \gls{laa}.
	
	The largest atrial modeling cohort study comprising 1000 biatrial models 
	has been reported in \cite{roney2023:_bia_vol}. 
	While this is indicative of robustness and scalability, 
	the numerous interactive manual interventions required suggest a labor intense procedure.
	Moreover, the workflow was primarily tailored for generating patient-specific bilayer models,
	as only these were derived from clinical imaging data.
	Support for generating volumetric models was limited to volumetric remeshing of pre-existing manifold meshes built from a statistical shape model \citep{rodero2021linking}, 
	and auxiliary data such as fiber architecture and \glspl{uac} were transferred over from an atlas.
	A direct generation of anatomically accurate patient-specific 3D models from images,
	as per our approach, is currently not supported.
	Similarly, in the study by \cite{nagel2021:_biatrial_ssm}, 
	volumetric biatrial models were also generated using statistical shape model manifolds as input
	to a heterogeneous workflow that integrated five different mesh generation and manipulation tools.

	\subsection{Performance considerations}
	The scalability of model generation methods depends on two major factors,
	the degree of automation to limit time-consuming interactive processing,
	and the robustness and speed of the processing steps, from medical image to final mesh.
	Our workflow achieves the highest degree of automation reported so far, 
	with minimal to no interactive processing at all. 
	Using a \gls{scn} \citep{thaler2021efficient}
	all \gls{ct} data sets were automatically segmented within $<$\SI{10}{\second}.
	As all required labels were detected robustly and accurately from all contrast CT datasets used in this study, 
	with no need for manual label correction.
	The additional landmarking of veins and appendages can be fully automated, requiring minimal user input to verify or, if necessary, correct the auto-generated labels. This process includes an interactive step supported by a custom tool designed for quick label correction, reducing the overall operation time to less than \SI{1}{\minute}.
	For all subsequent steps, including the generation of the volumetric atrial walls,
	the assignment of a rule-based fiber architecture and the computation of \glspl{uac},
	full automation has been achieved. 
	Performance measured over 50 models yielded execution times of only $<$\SI{10}{\minute} and $<$\SI{30}{\minute}
	for generating an atrial model at mesh resolutions suitable for \gls{re} and \gls{rd} simulations, respectively.
	

	\subsection{Mesh quality}
	A common limitation of most previous approaches is the mesh quality of the generated models. To evaluate this, we randomly sampled publicly available datasets of biatrial anatomies and found degenerated elements in all the sampled meshes \citep{nagel2021:_biatrial_ssm,roney2023:_bia_vol}.
	These deficiencies went unnoticed as all studies employed \gls{rd} monodomain models
	which cope well with poor mesh quality.
	However, employing these meshes for higher fidelity bidomain-based \gls{egm} generation is often challenging, 
	since it requires the solution of sensitive elliptic problems, whose solution may be associated with slow convergence or even solver divergence in the presence of degenerated elements when employing classical solvers \citep{conley2016overcoming,schneider2018decoupling}. 
	Mesh quality issues and topological errors may arise from volumetric mesh construction methods that depend on the direct extrusion of the atrial walls from the mesh manifolds, \citep{azzolin2023:_augmenta,nagel2021:_biatrial_ssm}. In contrast, our approach is more robust, as the atrial wall volume is generated directly from the image stack in a volumetric manner. This ensures that the initial volumetric representation of the atrial myocardium is topologically sound.
	Although our approach is less prone to topological issues, mesh topology errors may still occur at later stages, such as during surface smoothing operations aimed at removing jaggedness from the voxel-based grid or during remeshing to align mesh resolution with a prescribed target. To enhance robustness, a set of topological clean-up operations is repeatedly applied after each critical meshing step. All 50 models in this study were meshed fully automatically, achieving element quality exceeding a threshold of 0.99, which has been empirically identified as critical for cardiac bidomain simulations.
	
	\subsection{Universal Atrial Coordinates} 
	
	Anatomical reference frames are key for scalable modeling studies,
	as they provide a parametric encoding of all spatial model properties
	and, thus, facilitate an unattended parameter manipulation as required for parameter sweeps \citep{bayer2018universal,roney2019uac}.
	Accurate and robust mapping based on such reference frames relies on features such as
	geometric linearity, that is, isolines are evenly distributed in space,
	and uniqueness to ensure that any parameter set encodes a unique location in space. 
	For atrial models, the first universal spatial reference system has been reported 
	for manifold models by \cite{RONEY201965}.
	This approach was modified through anatomical normalization \citep{roney2021constructing} and later extended to volumetric models \citep{roney2023:_bia_vol}. The volumetric extension relies on \glspl{uac} defined on endocardial and epicardial manifold surfaces, which serve as boundary conditions for solving a \gls{ld} problem on the volumetric mesh to achieve transmural interpolation. Additionally, by solving another \gls{ld} problem, a normalized transmural distance field is computed.
	
	While this approach is straightforward to implement, 
	the generated volumetric \glspl{uac} were of moderate quality.
	Shortcomings of this approach include manual selection of initial landmark on the \gls{la} and \gls{ra} surface, which impedes a fully automated computation of the \gls{uac}, and the construction of geodesic path to define interfaces between lateral and septal, and anterior and posterior regions of the atria, to be used as boundary conditions for \gls{ld} solves, which does not guarantee the consistent identifications of anatomical structures and reference point locations across different geometries. Moreover, \glspl{uac} are always initially built and normalized over surface biatrial models and later extended to volumetric tissue by linear projection on the transmural mesh nodes, thus impairing a sound definition of transmural interfaces and \gls{uac} distribution. Lastly, the volumetric \glspl{uac} were not tested across meshes with different resolutions.       

	Our approach overcomes these limitations by utilizing automatically selected anatomical landmarks and paths, and        
	directly implementing \glspl{uac} on volumetric meshes. We replace the anatomical paths previously used to divide the \gls{la} into posterior-anterior parts and the \gls{ra} into lateral-septal parts with well-defined volumetric interfaces that separate not only the endo- and epicardial surface regions, but also the transmural mesh nodes. The computation of these interfaces is fully automated, based on pre-labeled anatomical regions and a set of rules to define reproducible points across anatomies. The resulting interfaces therefore ensure consistent identification of anatomical structures across different anatomical models. Moreover, they represent a well-defined transmural boundary for \gls{ld} computation, thus avoiding previous ambiguity in the distribution of the transmural coordinate.

	We demonstrate the robustness of our method for computing \glspl{uac} by testing it on 50 biatrial geometries from real \gls{af} patients. Additionally, we show the applicability of our method across different mesh resolutions. The computation of our \glspl{uac} is efficient on meshes suited for the \gls{re} model, though it requires longer computational times when applied to finer grids appropriate for \gls{rd} models. The majority of the computational costs were attributed to solving the linear elasticity problem to normalize the \gls{uac} solutions. This was mainly due to the poorly deformed elements that could potentially result from the projection of the Laplace-Dirichlet solutions on the \gls{uac} subspace, and hinder the solver convergence. Enhancing the computational efficiency of \glspl{uac} for finer meshes will be a focus of future work.
	
	
	\subsection{Modeling of inter-atrial connections}
	Electrical conduction between the atria is mediated by discrete \glspl{ic}
	which play an important \gls{ep} role.
	During normal sinus rhythm, the activation of the \gls{la} is mediated by \glspl{ic}, 
	determining the points of \gls{la} earliest activation which critically shape the mid- to terminal part of the P-wave.
	Further, \glspl{ic} provide pathways for inter-atrial reentrant circuits,  
	which represents a substrate for maintaining atrial flutter tachycardia.
	
	Detailed histological studies identified five electrical main connections 
	between the otherwise electrically isolated atria \citep{sakamoto2005interatrial}:
	the \gls{bb}, the coronary sinus bundle, an inferior-posterior bundle, a superior-anterior bundle, and the muscular ring around the fossa ovalis. 
	These histo-anatomical insights are helpful for providing fundamental anatomical constraints to a model,
	but individual anatomy and \gls{ep} properties may be highly variable.
	For instance, not all \glspl{ic} may be present or electrically active in a given individual,
	the locations of \gls{ra} origin and \gls{la} insertion,
	as well as the conduction properties of the \glspl{ic} in between, may also vary to a significant extent.
	Thus, this leads to different \glspl{ic} dominating inter-atrial activation and, thus, 
	govern the site of early activation of the \gls{la}, or the \gls{ra} for a retrograde atrial activation.
	As such, all factors defining the \glspl{ic} must be considered variable,
	and, thus, flexibility represented to support their inference. 
	
	Three different approaches for incorporating the \glspl{ic} in a biatrial model have been considered previously. These were: 
	i) timed electrical stimuli on the \gls{ra} to mediate inter-atrial conduction, without including any \glspl{ic} physical representation
	ii) anatomical bridges,  explicitly represented by meshing geometries of assumed \glspl{ic} \citep{krueger2012personalization,loewe2016influence,wachter2015mesh}, 
	and, 
	iii) single ohmic resistor connecting \gls{ra} and \gls{la}        
	\citep{roney2023:_bia_vol,roney2021constructing}. 
	All these approaches are limited in terms of flexibility and generalizability. 
	Timed stimulation is flexible, but does not represent the mechanism mediating inter-atrial conduction,
	and as such does not generalize to physiological representations such as the atrial flutter.
	Explicit anatomical meshing is of limited flexibility and robustness, 
	as any change in the origin or insertion site requires remeshing.
	This approach has been employed previously in \cite{krueger2012personalization}, 
	where, in absence of \glspl{uac}, complex \emph{ad-hoc} algorithms have been used 
	to parameterize and explicitly \glspl{ic} mesh \citep{loewe2016influence,wachter2015mesh}. 
	While feasible in principle, these approaches are complex and challenging to implement robustly, as the integration of thin inter-atrial strands with 3D atrial walls is highly prone to topological mesh errors.
	To provide an example, the explicit meshing of \gls{bb} used in the study by  \cite{nagel2021:_biatrial_ssm} 
	led to degenerated elements in all models throughout the cohort.
	Finally, fixed resistive 0D coupling, as used in various studies, using bilayer manifold models \citep{roney2023:_bia_vol,roney2021constructing}
	is flexible, but the approach is not suitable for connecting volumetric atria 
	due to dimensionality mismatch between 0D and 3D elements, which may cause source-sink mismatches. 
	In principle, this problem may arise also in manifold models where 0D resistors are used 
	to couple the two 2D surfaces representing \gls{ra} and \gls{la}.
	In the absence of a physiological interpretation of the coupling resistor, 
	viable ranges for which transduction between the atria is feasible have to be found by trial and error \citep{roney2023:_bia_vol}.
	Moreover, a further limitation is the inability to impose physiological constraints on the inter-atrial conduction delays,
	as transduction is rigidly dictated by the membrane time constant and the choice of coupling resistance.

	Our novel approach to incorporating \glspl{ic} in atrial models overcomes these limitations,
	as it is highly flexible, versatile, robust, and computationally efficient,
	and adds \gls{ic} anatomy as a parameter that can be probed during model calibration. 
	An arbitrary number of \glspl{ic} can be included to connect arbitrary sites of origin and insertion 
	on \gls{ra} and \gls{la}, respectively, 
	with prescribed conduction velocities governing the activation delay across the \gls{ic}.
	The use of \glspl{uac} to define origin and insertion sites, facilitates a seamless mapping of  \glspl{ic} 
	between different anatomies.
	Further, modeling the \glspl{ic} using a cable formulation computed with the same approach of generating the \gls{hps}, overcomes 
	dimension mismatch problems between 1D and 3D model components. 

	\subsection{Calibration of atrial \gls{ep} models - using the P-wave as objective}
	Atrial \gls{ep} modeling studies aiming to achieve a patient-specific calibration 
	by 1:1 matching with direct observations
	-- \glspl{eam}, \glspl{egm} or the \gls{ecg} -- are rare.
	Rather, adjustments of conductive and cellular dynamics properties informed by literature data
	are implemented to match global metrics such as the total atrial activation time.
	For a more accurate spatial calibration, \gls{eam} activation maps have been used 
	to infer intrinsic tissue conductivities or conduction velocities \citep{lubrecht2021:_piemap}. 
	High density \gls{eam} datasets covering both the \gls{ra} and the \gls{la} 
	appear ideal for calibrating atrial \gls{ep} models, 
	as they may provide a detailed view of the overall atrial activation sequence.
	However, while such datasets can be acquired, in principle,
	they are not widely available, as their acquisition tends to prolong procedures.
	Moreover, the spatiotemporal registration uncertainty between measurements of an \gls{eam} manifold 
	and the image-derived model is significant, 
	posing major challenges for calibration procedures. 
	Using \Glspl{egm} from which \gls{eam} maps are derived is even more challenging, 
	as \glspl{egm} provide only a very local view on tissue activation. 
	These may be helpful for inferring tissue characteristics from fibrotic patches,
	e.g.\ \gls{egm} magnitude and temporal separation of fractionated complexes in \glspl{egm},
	correlate with the severity and structure of fibrosis, but do not offer benefits 
	over using \gls{eam} maps directly. 
	
	In principle, the P-wave in the \gls{ecg} appears a most natural choice for calibrating atrial \gls{ep} models,
	as it is abundantly available and can be non-invasively recorded.
	However, the P-wave provides only a global view of the atrial activation sequence, and
	inferring space-varying model parameters from such limited data is challenging.
	Owing to the ill-posed nature of this inverse problem, this may not even be feasible.
	Also challenging but, in principle, feasible, as shown recently for the ventricular activation sequence \citep{gillette2021:_framework,grandits2023:_geodesic_bp},
	is the identification of model parameters that produce atrial activation sequences  
	and replicate, with high fidelity, the P-wave in a clinical standard \gls{ecg}.
	However, there may be more than one, potentially many, parameter sets yielding the same P-wave
	which raises questions of identifiability and uniqueness \citep{grandits2024:_identify}.
	Due to these challenges, only a few attempts have been made to calibrate an atrial \gls{ep} model to the P-wave \citep{loewe2015left}, mostly limited to qualitative visual comparisons between simulated and recorded P-waves, without a quantitative assessment of the differences, as shown in this work. 
	
	Further obstacle to P-wave-based model calibration can stem from incomplete parameterization that fails to accurately produce the actual activation sequence, and insufficient computational performance that hinders a feasible calibration process, especially when relying on traditional \gls{rd} models \citep{deng2012simulation,fedele2023comprehensive,ferrer2015detailed,loewe2015left,nagel2022comparison,roney2023:_bia_vol}
	combined with high fidelity pseudo-bidomain or lower fidelity extracellular potential recovery models \citep{labarthe2014:_bilayer,nagel2021:_biatrial_ssm}, or limited fidelity in the predicted P-wave.       
	Bilayer models are popular as they mitigate performance issues by reducing the overall problem size, 
	but are still orders of magnitude slower compared to volumetric \gls{re} models as used in this study,
	and the fidelity of predicted P-waves is limited.
	Volumetric biatrial \gls{rd} models immersed in a torso model using a full- or pseudo-bidomain 
	are able to produce high-fidelity P-waves \citep{loewe2015left}.
	They offer the advantage of providing the entire potential field $\phi{\rm e}$ throughout the torso, 
	in which \gls{ecg} or more extensive body surface potential maps are embedded \citep{ferrer2015detailed, ZAPPON2024112815}. 
	However, the approach is computationally by far too expensive for calibration studies, 
	and the additional information on the torso potential field is not readily exploited 
	in the absence of body surface potential mapping data.
	
	Our analysis clearly demonstrates 
	that the \gls{relf} model, combined with an extensive parameter space -- including also the anatomy of the \glspl{ic} -- is able to produce activation sequences
	and associated P-waves with full bidomain fidelity, and at real-time performance.
	As shown previously for a ventricular \gls{ecg} model \citep{gillette2021:_framework}, and in line with atrial \gls{ecg} model \citep{nagel2022comparison},
	the \gls{relf} model yields \glspl{egm} and \glspl{ecg} that are not discernible from those 
	produced with a full fidelity \gls{rd} bidomain model, as anticipated on theoretical grounds \citep{geselowitz1989theory}. 
	As shown in Figure \ref{Fig.RD_vs_RE}, differences in activation sequence and P-wave are minor,
	and can be deemed negligible in view of the overall model uncertainties. 
	
	In our study we focused on demonstrating the flexibility and efficiency of our framework 
	in supporting comprehensive automated exploration of physical and geometrical parameter spaces,
	and refrained from attempting to accurately calibrate the atrial \gls{ep} model to the P-wave for various reasons.
	As the atrial model was integrated with torso geometry and electrode positions taken from a different subject,
	and none of the pathologies of the patient treated by AF ablation were considered,
	discrepancies can be anticipated, and the feasibility of achieving a good fit is not even guaranteed. 
	Nonetheless, we quantified differences between simulated and real \glspl{ecg} using the \gls{rmse}, 
	showing that the P-wave can be closely captured in most leads, 
	except aVL, aVR, and V1. 
	P-wave features such as positivity in lead aVL and -aVR, as well as the biphasic characteristic of lead V1 were replicated by the model. 
	
	Importantly, the fidelity and efficiency of our framework support, in principle,
	a comprehensive exploration of the parameter space spanned by atrial anatomy and \gls{ep},  
	and, thus, provides a basis for building future applications geared towards creating digital twins 
	based on P-wave calibration.
	A full single forward simulation lasted $\approx$ \SI{27}{\second} only 
	where the actual evaluation of the \gls{ep} model amounted only to $\approx$ \SI{3}{\second}.  

	\subsection{Role of left and right atrial electrophysiology in the genesis of the P-wave}
	Despite the vast acceleration in speed achieved with high fidelity \gls{ep} forward models, such as the \gls{relf},
	the use \emph{a priori} knowledge on the genesis of the P-wave is key to constrain the inference 
	by limiting the admissible parameter space and the potentially large number of different atrial activation sequences
	that produce the same P-wave.
	
	From a macroscopic perspective, the genesis of the P-wave is well understood.
	However, at a mesoscopic size scale, that is  
	the relation between the physics of depolarization wavefronts traversing individual atrial structures,
	and their relative contribution to the global P-wave, is still unclear.
	These aspects have been investigated only in a limited number of studies  \citep{ferrer2015detailed,loewe2016influence,loewe2015left}.
	In \cite{loewe2015left}, the separate contribution of \gls{ra} and \gls{la} depolarization 
	to the P-wave was investigated in two atrial and torso models derived from healthy subjects. 
	Using a pre-defined activation sequence, the study showed 
	that the \gls{ra} predominantly influences the P-wave in precordial leads V1 and V2, 
	as well as limb leads II, aVF, and III, while the \gls{la} governs the central to terminal portion of the P-wave. 
	The effect of varying the atrial activation sequence on P-wave features such as positivity or negativity in different leads was investigated in \cite{loewe2016influence} 
	in eight biatrial models of healthy subjects by initiating atrial activation at a fixed set of locations 
	representing the \gls{san} exit sites.
	The relative contributions of different atrial regions to the P-wave was moreover analyzed in \cite{ferrer2015detailed},
	for a fixed atrial activation sequence, with prescribed \gls{san} exit sites and \gls{ic} locations.
	
	
	In our study, we investigate the envelope of P-waves associated with changes in both 
	\gls{ra} and \gls{la} anatomical and physiological factors, including the location and number of early activation sites in the \gls{la}, the position and shape of the \gls{san}, and the conduction velocity in the fiber and sheet directions of the \gls{ra} body. 
	By separately analyzing the \gls{ecg} variations due to these factors, we also highlight the roles of \gls{ra} and \gls{la} activation in shaping the P-wave.
	
	The generated P-wave envelop did not fully cover the P-wave observed in the modeled patient,
	most evident in the leads aVL and aVR (refer to Figure \ref{Fig.Fig_SAN_sampling}). 
	The terminal phase of the P-wave is not well covered, most prominently visible in the more lateral precordial leads V3-V6,
	and in the leads aVR and II.
	This discrepancy may stem from various factors. Firstly, the model was built from a patient treated by AF ablation,
	with a significant fibrotic burden, which remained unaccounted for in our model, and was instead assumed to be structurally healthy. 
	As such, the dipoles at the interface between healthy atrial tissue and fibrotic patches are missing,
	leading to a nearly zero dipole once both atria were fully activated.
	Finally, torso and lead positions were taken from another subject and were not specific to the given subject.
	
	\subsubsection{Discrete interatrial conduction}
	The site of earliest activation of the \gls{la} is governed by the \gls{ra} activation sequence
	and the location and conductive properties of the \glspl{ic}.
	In this study, we only investigated the effect of varying location and size of the insertion site
	of \gls{bb} on the anterior wall of the \gls{la} upon the P-wave. 
	The use of bundles for modeling \glspl{ic} facilitates automated sweeps over these important parameters
	that govern the coupling spatio-temporal coupling of \gls{ra} and \gls{la}.
	Consistent with \cite{loewe2015left}, our simulations show 
	that \gls{la} activation influences only the mid to terminal portion of the P-wave (see Figure \ref{Fig.Tests_1BB}).
	Major differences due to varying \gls{bb} insertion in the \gls{la} were witnessed in leads aVL, I, aVR, and III, 
	while effects in other limb leads and precordial leads were minor.
	Using three \glspl{ic} to model a fan-like insertion of the \gls{bb} into the anterior \gls{la} 
	had led only to a minor spread in the P-wave envelope (see Figure \ref{Fig.Three_bb}).

	\subsubsection{Sino-atrial node}
	The \gls{san} was originally described as a crescent-shaped area 
	at the junction of the \gls{svc} and the \gls{ra} \citep{anderson1983surgical,keith1907form}, and later found to extend along the \gls{ct}, and toward the \gls{ivc} \citep{akima1978method,anderson1978development,boyett2000sinoatrial,csepe_human_2016, fedorov_conduction_2012,monfredi2010anatomy}. 
	Normal \gls{san} automaticity depends on its depolarized resting potential 
	of $\approx$ \SI{-60}{\milli\volt} relative to neighboring atrial tissue,
	with a more hyperpolarized resting potential of $\approx$ \SI{-85}{\milli\volt} \citep{csepe_human_2016}.
	To prevent hyperpolarization of \gls{san} and, thus, inhibition of automaticity,
	the \gls{san} is electrically mostly insulated to safeguard pacemaker function \citep{fedorov_conduction_2012,joyner_propagation_1986}.
	Electrical coupling of the \gls{san} with surrounding atrial myocardium 
	is limited to specific discrete sites,referred to as \gls{san} exit pathways, 
	which can number up to five \citep{fedorov_optical_2010,li_redundant_2017}.
	As such, the location of exit pathways and their activity 
	determine the earliest atrial activation site and the shape of the initial depolarization wavefront, 
	which not only affects the activation sequence of the \gls{ra}, 
	but may also influence the order of activation of the \glspl{ic}, and, in turn, alter the earliest activation site on the \gls{la} \citep{antz_electrical_1998,li_redundant_2017}.
	Consequently, \gls{san} anatomy and exit pathways are crucial parameters to consider 
	when modeling atrial activation and its reflection in the P-wave.
	
	In modeling studies, the \gls{san} is often represented as a spherical region at a fixed location  \citep{fedele2023comprehensive,ferrer2015detailed,loewe2016influence}, 
	with only a few studies accounting for the anatomical shape of the \gls{san}         \citep{gillette2022personalized,labarthe2014:_bilayer}.
	In line with previous work \citep{ferrer2015detailed,loewe2016influence}, here we demonstrate 
	a marked dependency of the P-wave upon exit sites (refer to Figure \ref{Fig.Fig_SAN_sampling}). By easy variation of \gls{san} size and location, we show that our framework may easily allow for inferring both \gls{san} anatomy and exit sites based on \gls{ecg} calibration.

	\subsubsection{Effect of the \gls{ra} endocardium on the P-wave}
	The anterior and lateral endocardial walls of the \gls{ra} are primarily composed of \glspl{pm} and \gls{ct} tissue and are attached to a thin-walled \gls{la} epicardium, consisting of only a few layers of myocytes \citep{lang2022imaging,matsuyama2004anatomical}. 
	While explicitly accounted for in bilayer models, this is often not the case in volumetric models \citep{roney2023:_bia_vol}.
	There, two different approaches have been used. Either, \gls{ct}, \glspl{pm} and the \gls{ra} portion of \gls{bb}
	were explicitly meshed, or a continuous right endocardial layer of prescribed thickness was generated
	to represent the \gls{ra} endocardium, 
	onto which the locations of the \gls{ct} and \glspl{pm} are projected either by using an anatomical atlas, or on a per-rule basis.
	Atrial \gls{ep} in the \gls{ra} endocardium is then simulated 
	by assigning faster conduction velocities in the \gls{ct} and \glspl{pm}, 
	and slower conduction velocities to the remaining endocardial tissue between the protruding structures \citep{azzolin2023:_augmenta,ferrer2015detailed,loewe2015left,roney2023:_bia_vol}. 
	As such, the absence of tissue between the \gls{ct} and \glspl{pm} is typically not accounted for 
	which can be interpreted as an increase in the effective epicardial wall width and the associated source strength. 
	
	We investigate the effect of modeling the \gls{ra} endocardial tissue as an electrically active layer 
	upon atrial sources and the P-wave. 
	Treating the \gls{ra} endocardial layer in between \gls{ct} and \glspl{pm} as blood pool 
	slowed down the activation of the posterior-lateral wall of the \gls{ra} and the \gls{raa}. 
	While this had a minor impact on the total activation time of the \gls{ra}, 
	a noticeable reduction in peak P-wave amplitude by $\approx$1\% of was witnessed, with no change in P-wave morphology. 
	The observed deviations in P-wave can be deemed as minor relative to the overall uncertainty.
	This may be not the case when simulating endocardial recorded \glspl{egm},
	when electrodes are located near the atrial sources.

	\section{Limitations}
	\label{Sec:limitation}
	While our comprehensive end-to-end workflow for generating biatrial anatomies combined with a fast forward \gls{relf} model 
	is able to deliver anatomically accurate models for predicting \gls{ecg} and \glspl{egm} with full physical fidelity,
	numerous notable limitations exist.
	
	Our study introduces a technological framework for the calibration of atrial \gls{ep} models 
	using observable extracellular recordings such as the \glspl{egm} and the \gls{ecg} -- specifically, the P-wave -- 
	as objective to infer the high dimensional parameter space of such models. 
	However, in its current form, an immediate inference to match a P-wave is not readily feasible. 
	Despite the computational efficiency, a brute force sampling approach may still be too expensive \citep{gillette2021:_framework}.
	Similar to previous work on identifying the ventricular activation sequence, a dedicated approach 
	combining \emph{a priori} knowledge to constrain the search space in combination with an optimization approach 
	that facilitates a gradient-based search for parameters is likely to be more accurate and efficient \citep{grandits2021:_geasi,grandits2023:_geodesic_bp}.
	Sampling narrow, physiologically constrained parameter corridors, as done here, may provide a better approximation than current standard methods. However, this improvement cannot be guaranteed for two key reasons.        
	First, as the atrial activation sequence cannot be measured accurately, not even invasively with most advanced \gls{eam} methodology,
	a solid ground truth reference is missing.
	Secondly, it is highly likely that a large number of parameter choices exist 
	which yield identical perfect matches of the \glspl{ecg}.
	Thus, uniqueness and identifiability is a serious concern that must be adequately addressed 
	to narrow down the number of different sampled parameters.
	The development of suitable strategies and their evaluation was beyond the scope of this study,  and substantial further research efforts are required.
	
	An important factor limiting the anatomical accuracy of models generated by our end-to-end workflow 
	is the assumption of uniform wall thicknesses. 
	As the accurate segmentation of the atrial walls from current clinical images is challenging,
	in our volumetric wall reconstruction, the atrial wall widths were prescribed and held constant over individual anatomical regions. 
	However, atrial wall widths in a patient's heart may be different and more heterogeneous than represented in the model.
	As a consequence, the local variation in local conduction velocity mediated by altered source-sink relation may not be accurately represented.
	
	A further limitation is the reduced computational efficiency in generating finer resolutions for \gls{rd} models.
	This is largely caused by the computation of the \glspl{uac}, and not due to the extra cost of generating finer meshes per se. Specifically, projecting the biatrial mesh onto the \gls{uac} space to solve the linear elasticity problem 
	can induce mesh deformations and, consequently, iterative solver convergence issues.
	More advanced preconditioners will be incorporated \citep{augustin2016anatomically} in the future to enhance computational efficiency. 
	Due to the nearly identical anatomical meshes, a faster alternative may be 
	the Euclidian interpolation of \glspl{uac} from the lower resolution \gls{re} to the higher resolution \gls{rd} mesh, but this has not been investigated here.
	
	While our novel flexible and rapid approach to modeling \glspl{ic} is advantageous for model calibration, 
	the used of single cables only approximates the effects of conduction between the atria, 
	but not the actual sources contributing to the extracellular potential fields.
	The actual anatomical structures underlying the \glspl{ic} are rather muscular bands or 3D strands. 
	As such, the volumetric source density of the cable approach underestimates the contribution to the potential field
	in their vicinity.
	While this does not constitute a problem for simulating activation sequences and the \gls{ecg}, 
	simulating the sensing of \gls{egm} will be inaccurate and also methodologically challenging to implement, 
	owing to the non-conformity of \gls{ic} grid with the atrial or the volume conductor mesh.
	This problem is readily mitigated by using bundles instead of single cables,
	which may be straightforwardly incorporated, as shown for modeling a fan-like insertion of the \gls{bb},
	where a bundle of three \glspl{ic} was used (see Figure \ref{Fig.Three_bb}). 
	Alternatively, the cable-based representation can be replaced with an explicit mesh of the \glspl{ic}
	once the model has been calibrated.
	
	Further, the cables representing the \glspl{ic} are automatically generated, 
	following the geodesic paths connecting two mesh nodes on \gls{ra} and \gls{ra}, respectively. 
	This might result in a mismatch in trajectory and length between modeled and real \gls{ic} and, thus, physiologically non-feasible conduction velocities may be required to meet prescribed inter-atrial conduction delays.
	Specifically, this may complicate the transfer of parameters from \gls{re} to \gls{rd} models 
	where conduction velocity ranges are limited.
	
	Further, we have not included the segmentation of the torso along with the corresponding atrial geometries. 
	Large clinical data cohorts often lack \gls{ct} or \gls{mri} scans of the torso. 
	Nevertheless, our workflow supports the meshing of a generic torso surface around the generated atrial geometries, 
	making the integration of automatic torso segmentation and corresponding meshing into our workflow a straightforward enhancement.
	
	Finally, we have not yet incorporated considerations for atrial disorders and pathologies, 
	which would be essential for clinical applications. While idealized damage and non-conductive substrates can be easily modeled using our \gls{ep} simulator—along with pacing protocols to simulate conditions such as \gls{af}, our current \gls{relf} model does not support \gls{ecg} simulations for rhythm disorders with reentrant phenotypes, as the underlying Eikonal model driving activation does not accommodate reentrant patterns. 
	Advancements in altering atrial depolarization profiles may be necessary. Furthermore, real damaged tissue data from clinical sources is not yet integrated into our workflow. However, our reference framework is designed to easily accommodate new functionalities and define complex parameter variations, paving the way for the automatic and near-real-time creation of atrial digital twins.
	
	\section{Conclusions}    
	\label{Sec:conclusion}
	In this study, we developed an automated, efficient, and flexible multi-scale workflow for creating anatomically accurate biatrial models with high-quality meshes, suitable for various forward \gls{ep} models, and extended it to a function twinning stage, where efficient model calibration based on \gls{ecg} data is possible. The workflow utilizes a segmented blood pool and a deep learning-based segmentation network to automatically identify and label principal atrial structures. It moreover generates endocardial and epicardial walls based on predefined rules, ensuring robust model openings and high mesh quality. Detailed anatomical and functional structures, including fiber architecture and space-varying parameters, are incorporated into the models. 
	We further introduced a flexible method for representing \glspl{ic} using auto-generated cables, allowing for realistic and customizable conduction properties. Our framework finally includes a clinically compatible \gls{ecg} forward generation system, integrating both \gls{rd} and \gls{re} models with Lead Field methods for accurate \gls{ecg} trace generation, and possibly allows for the integration of the Pseudo-bidomain approach or $\phi_e$-recovery. Among these, we show that the \gls{relf} method proved effective for real-time \gls{ecg} computation and reliable source approximation in sinus rhythm. 
	
	Overall, our workflow facilitates efficient exploration of the \gls{ep} parameter space for calibration of biatrial models based on \gls{ecg} data, making it a valuable tool for advancing cardiac \gls{ep} research.
	
	\FloatBarrier
	
	\section*{Acknowledgments}
	
	This research received support from the Austrian Science Fund (FWF) grant no. 10.55776/I6476, and from the European Union’s Horizon 2020 research and innovation program under the Marie Sk\l{}odowska-Curie grant TwinCare-AF (grant agreement no. 101148636). The National Institutes of Health supported this work under grant 1R01HL158667 and by the Austrian Science Fund (FWF) under grant 10.55776/P37063.
	E. Zappon acknowledges her membership to INdAM GNCS - Gruppo Nazionale per il Calcolo Scientifico (National Group for Scientific Computing, Italy).
	\begin{figure}[!h]
		\centering
		\includegraphics[width = 0.4\textwidth]{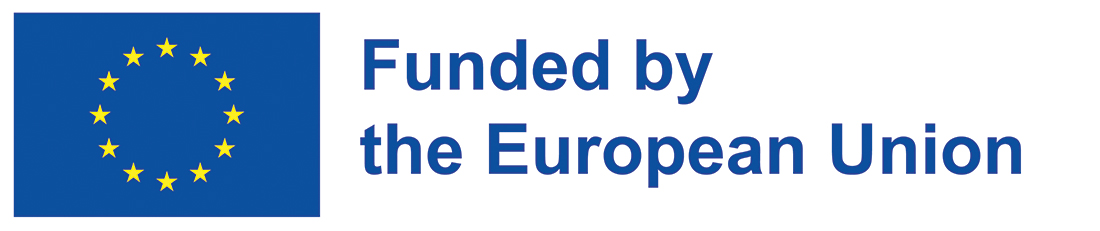}
	\end{figure}

	\bibliographystyle{model2-names}

	\bibliography{References}

\begin{thebibliography}{118}
\expandafter\ifx\csname natexlab\endcsname\relax\def\natexlab#1{#1}\fi
\providecommand{\url}[1]{\texttt{#1}}
\providecommand{\href}[2]{#2}
\providecommand{\path}[1]{#1}
\providecommand{\DOIprefix}{doi:}
\providecommand{\ArXivprefix}{arXiv:}
\providecommand{\URLprefix}{URL: }
\providecommand{\Pubmedprefix}{pmid:}
\providecommand{\doi}[1]{\href{http://dx.doi.org/#1}{\path{#1}}}
\providecommand{\Pubmed}[1]{\href{pmid:#1}{\path{#1}}}
\providecommand{\bibinfo}[2]{#2}
\ifx\xfnm\relax \def\xfnm[#1]{\unskip,\space#1}\fi
\bibitem[{Akima(1978)}]{akima1978method}
\bibinfo{author}{Akima, H.}, \bibinfo{year}{1978}.
\newblock \bibinfo{title}{A method of bivariate interpolation and smooth
  surface fitting for irregularly distributed data points}.
\newblock \bibinfo{journal}{ACM Transactions on Mathematical Software (TOMS)}
  \bibinfo{volume}{4}, \bibinfo{pages}{148--159}.
\bibitem[{Anderson et~al.(1983)Anderson, Ho and Becker}]{anderson1983surgical}
\bibinfo{author}{Anderson, R.H.}, \bibinfo{author}{Ho, S.Y.},
  \bibinfo{author}{Becker, A.E.}, \bibinfo{year}{1983}.
\newblock \bibinfo{title}{The surgical anatomy of the conduction tissues.}
\newblock \bibinfo{journal}{Thorax} \bibinfo{volume}{38},
  \bibinfo{pages}{408--420}.
\bibitem[{Anderson et~al.(1978)Anderson, Siew~Yen, Becker and
  Gosling}]{anderson1978development}
\bibinfo{author}{Anderson, R.H.}, \bibinfo{author}{Siew~Yen, H.},
  \bibinfo{author}{Becker, A.E.}, \bibinfo{author}{Gosling, J.A.},
  \bibinfo{year}{1978}.
\newblock \bibinfo{title}{The development of the sinoatrial node}, in:
  \bibinfo{booktitle}{The Sinus Node: Structure, Function, and Clinical
  Relevance}. \bibinfo{publisher}{Springer}, pp. \bibinfo{pages}{166--182}.
\bibitem[{Antz et~al.(1998a)Antz, Otomo, Arruda, Scherlag, Pitha, Tondo,
  Lazzara and Jackman}]{antz1998electrical}
\bibinfo{author}{Antz, M.}, \bibinfo{author}{Otomo, K.},
  \bibinfo{author}{Arruda, M.}, \bibinfo{author}{Scherlag, B.J.},
  \bibinfo{author}{Pitha, J.}, \bibinfo{author}{Tondo, C.},
  \bibinfo{author}{Lazzara, R.}, \bibinfo{author}{Jackman, W.M.},
  \bibinfo{year}{1998}a.
\newblock \bibinfo{title}{Electrical conduction between the right atrium and
  the left atrium via the musculature of the coronary sinus}.
\newblock \bibinfo{journal}{Circulation} \bibinfo{volume}{98},
  \bibinfo{pages}{1790--1795}.
\bibitem[{Antz et~al.(1998b)Antz, Otomo, Arruda, Scherlag, Pitha, Tondo,
  Lazzara and Jackman}]{antz_electrical_1998}
\bibinfo{author}{Antz, M.}, \bibinfo{author}{Otomo, K.},
  \bibinfo{author}{Arruda, M.}, \bibinfo{author}{Scherlag, B.J.},
  \bibinfo{author}{Pitha, J.}, \bibinfo{author}{Tondo, C.},
  \bibinfo{author}{Lazzara, R.}, \bibinfo{author}{Jackman, W.M.},
  \bibinfo{year}{1998}b.
\newblock \bibinfo{title}{Electrical {Conduction} {Between} the {Right}
  {Atrium} and the {Left} {Atrium} via the {Musculature} of the {Coronary}
  {Sinus}}.
\newblock \bibinfo{journal}{Circulation} \bibinfo{volume}{98},
  \bibinfo{pages}{1790--1795}.
\newblock \URLprefix
  \url{https://www.ahajournals.org/doi/10.1161/01.CIR.98.17.1790},
  \DOIprefix\doi{10.1161/01.CIR.98.17.1790}.
\bibitem[{Augustin et~al.(2016)Augustin, Neic, Liebmann, Prassl, Niederer,
  Haase and Plank}]{augustin2016anatomically}
\bibinfo{author}{Augustin, C.M.}, \bibinfo{author}{Neic, A.},
  \bibinfo{author}{Liebmann, M.}, \bibinfo{author}{Prassl, A.J.},
  \bibinfo{author}{Niederer, S.A.}, \bibinfo{author}{Haase, G.},
  \bibinfo{author}{Plank, G.}, \bibinfo{year}{2016}.
\newblock \bibinfo{title}{Anatomically accurate high resolution modeling of
  human whole heart electromechanics: a strongly scalable algebraic multigrid
  solver method for nonlinear deformation}.
\newblock \bibinfo{journal}{Journal of computational physics}
  \bibinfo{volume}{305}, \bibinfo{pages}{622--646}.
\bibitem[{Azzolin et~al.(2023)Azzolin, Eichenlaub, Nagel, Nairn, Sánchez,
  Unger, Arentz, Westermann, Dössel, Jadidi and Loewe}]{azzolin2023:_augmenta}
\bibinfo{author}{Azzolin, L.}, \bibinfo{author}{Eichenlaub, M.},
  \bibinfo{author}{Nagel, C.}, \bibinfo{author}{Nairn, D.},
  \bibinfo{author}{Sánchez, J.}, \bibinfo{author}{Unger, L.},
  \bibinfo{author}{Arentz, T.}, \bibinfo{author}{Westermann, D.},
  \bibinfo{author}{Dössel, O.}, \bibinfo{author}{Jadidi, A.},
  \bibinfo{author}{Loewe, A.}, \bibinfo{year}{2023}.
\newblock \bibinfo{title}{{AugmentA}: {Patient}-specific augmented atrial model
  generation tool}.
\newblock \bibinfo{journal}{Computerized Medical Imaging and Graphics}
  \bibinfo{volume}{108}, \bibinfo{pages}{102265}.
\newblock \URLprefix
  \url{https://linkinghub.elsevier.com/retrieve/pii/S0895611123000836},
  \DOIprefix\doi{10.1016/j.compmedimag.2023.102265}.
\bibitem[{Azzolin et~al.(2020)Azzolin, Luongo, Ventura, Saiz, D{\"o}sse and
  Loewe}]{azzolin2020influence}
\bibinfo{author}{Azzolin, L.}, \bibinfo{author}{Luongo, G.},
  \bibinfo{author}{Ventura, S.R.}, \bibinfo{author}{Saiz, J.},
  \bibinfo{author}{D{\"o}sse, O.}, \bibinfo{author}{Loewe, A.},
  \bibinfo{year}{2020}.
\newblock \bibinfo{title}{Influence of gradient and smoothness of atrial wall
  thickness on initiation and maintenance of atrial fibrillation}, in:
  \bibinfo{booktitle}{2020 Computing in Cardiology},
  \bibinfo{organization}{IEEE}. pp. \bibinfo{pages}{1--4}.
\bibitem[{Bachmann(1916)}]{Bachmann1916}
\bibinfo{author}{Bachmann, G.}, \bibinfo{year}{1916}.
\newblock \bibinfo{title}{The inter-auricular time interval}.
\newblock \bibinfo{journal}{American Journal of Physiology-Legacy Content}
  \bibinfo{volume}{41}, \bibinfo{pages}{309--320}.
\newblock \DOIprefix\doi{10.1152/ajplegacy.1916.41.3.309}.
\bibitem[{Batdorf et~al.(1997)Batdorf, Freitag, Ollivier-Gooch, Batdorf,
  Freitag and Ollivier-Gooch}]{batdorf1997computational}
\bibinfo{author}{Batdorf, M.}, \bibinfo{author}{Freitag, L.},
  \bibinfo{author}{Ollivier-Gooch, C.}, \bibinfo{author}{Batdorf, M.},
  \bibinfo{author}{Freitag, L.}, \bibinfo{author}{Ollivier-Gooch, C.},
  \bibinfo{year}{1997}.
\newblock \bibinfo{title}{Computational study of the effect of unstructured
  mesh quality on solution efficiency}, in: \bibinfo{booktitle}{13th
  Computational Fluid Dynamics Conference}, p. \bibinfo{pages}{1888}.
\bibitem[{Bayer et~al.(2018)Bayer, Prassl, Pashaei, Gomez, Frontera, Neic,
  Plank and Vigmond}]{bayer2018universal}
\bibinfo{author}{Bayer, J.}, \bibinfo{author}{Prassl, A.J.},
  \bibinfo{author}{Pashaei, A.}, \bibinfo{author}{Gomez, J.F.},
  \bibinfo{author}{Frontera, A.}, \bibinfo{author}{Neic, A.},
  \bibinfo{author}{Plank, G.}, \bibinfo{author}{Vigmond, E.J.},
  \bibinfo{year}{2018}.
\newblock \bibinfo{title}{Universal ventricular coordinates: A generic
  framework for describing position within the heart and transferring data}.
\newblock \bibinfo{journal}{Medical Image Analysis} \bibinfo{volume}{45},
  \bibinfo{pages}{83--93}.
\bibitem[{Beinart et~al.(2011)Beinart, Abbara, Blum, Ferencik, Heist, Ruskin
  and Mansour}]{beinart2011left}
\bibinfo{author}{Beinart, R.}, \bibinfo{author}{Abbara, S.},
  \bibinfo{author}{Blum, A.}, \bibinfo{author}{Ferencik, M.},
  \bibinfo{author}{Heist, K.}, \bibinfo{author}{Ruskin, J.},
  \bibinfo{author}{Mansour, M.}, \bibinfo{year}{2011}.
\newblock \bibinfo{title}{Left atrial wall thickness variability measured by ct
  scans in patients undergoing pulmonary vein isolation}.
\newblock \bibinfo{journal}{Journal of cardiovascular electrophysiology}
  \bibinfo{volume}{22}, \bibinfo{pages}{1232--1236}.
\bibitem[{Bhagirath et~al.(2024)Bhagirath, Strocchi, Bishop, Boyle and
  Plank}]{bhagirath2024:_bits}
\bibinfo{author}{Bhagirath, P.}, \bibinfo{author}{Strocchi, M.},
  \bibinfo{author}{Bishop, M.J.}, \bibinfo{author}{Boyle, P.},
  \bibinfo{author}{Plank, G.}, \bibinfo{year}{2024}.
\newblock \bibinfo{title}{From bits to bedside: Entering the age of digital
  twins in cardiac electrophysiology}.
\newblock \bibinfo{journal}{Europace} .
\bibitem[{Bishop and Plank(2025)}]{bishop24:_devarp_trap}
\bibinfo{author}{Bishop, M.}, \bibinfo{author}{Plank, G.},
  \bibinfo{year}{2025}.
\newblock \bibinfo{title}{Stochastic behaviour of arrhythmia induction in
  virtual heart models suggests caution in offering mechanistic insights.}
\newblock \bibinfo{journal}{Nature Cardiovascular Research} \bibinfo{note}{In
  press}.
\bibitem[{Bishop et~al.(2016)Bishop, Rajani, Plank, Gaddum, Carr-White, Wright,
  O'Neill and Niederer}]{bishop2016three}
\bibinfo{author}{Bishop, M.}, \bibinfo{author}{Rajani, R.},
  \bibinfo{author}{Plank, G.}, \bibinfo{author}{Gaddum, N.},
  \bibinfo{author}{Carr-White, G.}, \bibinfo{author}{Wright, M.},
  \bibinfo{author}{O'Neill, M.}, \bibinfo{author}{Niederer, S.},
  \bibinfo{year}{2016}.
\newblock \bibinfo{title}{Three-dimensional atrial wall thickness maps to
  inform catheter ablation procedures for atrial fibrillation}.
\newblock \bibinfo{journal}{Europace} \bibinfo{volume}{18},
  \bibinfo{pages}{376--383}.
\bibitem[{Bishop and Plank(2011a)}]{bishop2011bidomain}
\bibinfo{author}{Bishop, M.J.}, \bibinfo{author}{Plank, G.},
  \bibinfo{year}{2011}a.
\newblock \bibinfo{title}{Bidomain {ECG} simulations using an augmented
  monodomain model for the cardiac source}.
\newblock \bibinfo{journal}{IEEE transactions on biomedical engineering}
  \bibinfo{volume}{58}, \bibinfo{pages}{2297--2307}.
\bibitem[{Bishop and Plank(2011b)}]{bishop2011bathloading}
\bibinfo{author}{Bishop, M.J.}, \bibinfo{author}{Plank, G.},
  \bibinfo{year}{2011}b.
\newblock \bibinfo{title}{Representing cardiac bidomain bath-loading effects by
  an augmented monodomain approach: application to complex ventricular models.}
\newblock \bibinfo{journal}{IEEE transactions on bio-medical engineering}
  \bibinfo{volume}{58}, \bibinfo{pages}{1066--1075}.
\newblock \DOIprefix\doi{10.1109/TBME.2010.2096425}.
\bibitem[{Boyett et~al.(2000)Boyett, Honjo and Kodama}]{boyett2000sinoatrial}
\bibinfo{author}{Boyett, M.R.}, \bibinfo{author}{Honjo, H.},
  \bibinfo{author}{Kodama, I.}, \bibinfo{year}{2000}.
\newblock \bibinfo{title}{The sinoatrial node, a heterogeneous pacemaker
  structure}.
\newblock \bibinfo{journal}{Cardiovascular research} \bibinfo{volume}{47},
  \bibinfo{pages}{658--687}.
\bibitem[{Boyle et~al.(2019)Boyle, Zghaib, Zahid, Ali, Deng, Franceschi, Hakim,
  Murphy, Prakosa, Zimmerman et~al.}]{boyle2019computationally}
\bibinfo{author}{Boyle, P.M.}, \bibinfo{author}{Zghaib, T.},
  \bibinfo{author}{Zahid, S.}, \bibinfo{author}{Ali, R.L.},
  \bibinfo{author}{Deng, D.}, \bibinfo{author}{Franceschi, W.H.},
  \bibinfo{author}{Hakim, J.B.}, \bibinfo{author}{Murphy, M.J.},
  \bibinfo{author}{Prakosa, A.}, \bibinfo{author}{Zimmerman, S.L.}, et~al.,
  \bibinfo{year}{2019}.
\newblock \bibinfo{title}{Computationally guided personalized targeted ablation
  of persistent atrial fibrillation}.
\newblock \bibinfo{journal}{Nature biomedical engineering} \bibinfo{volume}{3},
  \bibinfo{pages}{870--879}.
\bibitem[{van Campenhout et~al.(2013)van Campenhout, Yaksh, Kik, de~Jaegere,
  Ho, Allessie and de~Groot}]{Margo2013}
\bibinfo{author}{van Campenhout, M.J.}, \bibinfo{author}{Yaksh, A.},
  \bibinfo{author}{Kik, C.}, \bibinfo{author}{de~Jaegere, P.P.},
  \bibinfo{author}{Ho, S.Y.}, \bibinfo{author}{Allessie, M.A.},
  \bibinfo{author}{de~Groot, N.M.}, \bibinfo{year}{2013}.
\newblock \bibinfo{title}{Bachmann’s bundle}.
\newblock \bibinfo{journal}{Circulation: Arrhythmia and Electrophysiology}
  \bibinfo{volume}{6}, \bibinfo{pages}{1041--1046}.
\newblock \DOIprefix\doi{10.1161/CIRCEP.113.000758}.
\bibitem[{Chauvin et~al.(2000)Chauvin, Shah, Ha\"{i}ssaguerre, Marcellin and
  Brechenmacher}]{chauvin2000anatomic}
\bibinfo{author}{Chauvin, M.}, \bibinfo{author}{Shah, D.C.},
  \bibinfo{author}{Ha\"{i}ssaguerre, M.}, \bibinfo{author}{Marcellin, L.},
  \bibinfo{author}{Brechenmacher, C.}, \bibinfo{year}{2000}.
\newblock \bibinfo{title}{The anatomic basis of connections between the
  coronary sinus musculature and the left atrium in humans}.
\newblock \bibinfo{journal}{Circulation} \bibinfo{volume}{101},
  \bibinfo{pages}{647--652}.
\bibitem[{Conley et~al.(2016)Conley, Delaney and Jiao}]{conley2016overcoming}
\bibinfo{author}{Conley, R.}, \bibinfo{author}{Delaney, T.J.},
  \bibinfo{author}{Jiao, X.}, \bibinfo{year}{2016}.
\newblock \bibinfo{title}{Overcoming element quality dependence of finite
  elements with adaptive extended stencil fem (aes-fem)}.
\newblock \bibinfo{journal}{International Journal for Numerical Methods in
  Engineering} \bibinfo{volume}{108}, \bibinfo{pages}{1054--1085}.
\bibitem[{Corradi et~al.(2011)Corradi, Maestri, Macchi and
  Callegari}]{corradi2011atria}
\bibinfo{author}{Corradi, D.}, \bibinfo{author}{Maestri, R.},
  \bibinfo{author}{Macchi, E.}, \bibinfo{author}{Callegari, S.},
  \bibinfo{year}{2011}.
\newblock \bibinfo{title}{The atria: from morphology to function}.
\newblock \bibinfo{journal}{Journal of Cardiovascular Electrophysiology}
  \bibinfo{volume}{22}, \bibinfo{pages}{223--235}.
\bibitem[{Corral-Acero et~al.(2020)Corral-Acero, Margara, Marciniak, Rodero,
  Loncaric, Feng, Gilbert, Fernandes, Bukhari, Wajdan, Martinez, Santos,
  Shamohammdi, Luo, Westphal, Leeson, DiAchille, Gurev, Mayr, Geris,
  Pathmanathan, Morrison, Cornelussen, Prinzen, Delhaas, Doltra, Sitges,
  Vigmond, Zacur, Grau, Rodriguez, Remme, Niederer, Mortier, McLeod, Potse,
  Pueyo, Bueno-Orovio and Lamata}]{corral2020digital}
\bibinfo{author}{Corral-Acero, J.}, \bibinfo{author}{Margara, F.},
  \bibinfo{author}{Marciniak, M.}, \bibinfo{author}{Rodero, C.},
  \bibinfo{author}{Loncaric, F.}, \bibinfo{author}{Feng, Y.},
  \bibinfo{author}{Gilbert, A.}, \bibinfo{author}{Fernandes, J.F.},
  \bibinfo{author}{Bukhari, H.A.}, \bibinfo{author}{Wajdan, A.},
  \bibinfo{author}{Martinez, M.V.}, \bibinfo{author}{Santos, M.S.},
  \bibinfo{author}{Shamohammdi, M.}, \bibinfo{author}{Luo, H.},
  \bibinfo{author}{Westphal, P.}, \bibinfo{author}{Leeson, P.},
  \bibinfo{author}{DiAchille, P.}, \bibinfo{author}{Gurev, V.},
  \bibinfo{author}{Mayr, M.}, \bibinfo{author}{Geris, L.},
  \bibinfo{author}{Pathmanathan, P.}, \bibinfo{author}{Morrison, T.},
  \bibinfo{author}{Cornelussen, R.}, \bibinfo{author}{Prinzen, F.},
  \bibinfo{author}{Delhaas, T.}, \bibinfo{author}{Doltra, A.},
  \bibinfo{author}{Sitges, M.}, \bibinfo{author}{Vigmond, E.J.},
  \bibinfo{author}{Zacur, E.}, \bibinfo{author}{Grau, V.},
  \bibinfo{author}{Rodriguez, B.}, \bibinfo{author}{Remme, E.W.},
  \bibinfo{author}{Niederer, S.}, \bibinfo{author}{Mortier, P.},
  \bibinfo{author}{McLeod, K.}, \bibinfo{author}{Potse, M.},
  \bibinfo{author}{Pueyo, E.}, \bibinfo{author}{Bueno-Orovio, A.},
  \bibinfo{author}{Lamata, P.}, \bibinfo{year}{2020}.
\newblock \bibinfo{title}{{The 'Digital Twin' to enable the vision of precision
  cardiology.}}
\newblock \bibinfo{journal}{Eur. Heart J.} ,
  \bibinfo{pages}{1--11}\DOIprefix\doi{10.1093/eurheartj/ehaa159}.
\bibitem[{Costa et~al.(2014)Costa, Campos, Prassl, Dos~Santos,
  Sanchez-Quintana, Ahammer, Hofer and Plank}]{costa2014:_fibrotic}
\bibinfo{author}{Costa, C.M.}, \bibinfo{author}{Campos, F.O.},
  \bibinfo{author}{Prassl, A.J.}, \bibinfo{author}{Dos~Santos, R.W.},
  \bibinfo{author}{Sanchez-Quintana, D.}, \bibinfo{author}{Ahammer, H.},
  \bibinfo{author}{Hofer, E.}, \bibinfo{author}{Plank, G.},
  \bibinfo{year}{2014}.
\newblock \bibinfo{title}{An efficient finite element approach for modeling
  fibrotic clefts in the heart}.
\newblock \bibinfo{journal}{IEEE Transactions on Biomedical Engineering}
  \bibinfo{volume}{61}, \bibinfo{pages}{900--910}.
\newblock \DOIprefix\doi{10.1109/TBME.2013.2292320}.
\bibitem[{Courtemanche et~al.(1998)Courtemanche, Ramirez and
  Nattel}]{courtemanche1998ionic}
\bibinfo{author}{Courtemanche, M.}, \bibinfo{author}{Ramirez, R.J.},
  \bibinfo{author}{Nattel, S.}, \bibinfo{year}{1998}.
\newblock \bibinfo{title}{Ionic mechanisms underlying human atrial action
  potential properties: insights from a mathematical model}.
\newblock \bibinfo{journal}{American Journal of Physiology-Heart and
  Circulatory Physiology} \bibinfo{volume}{275}, \bibinfo{pages}{H301--H321}.
\bibitem[{Coveney et~al.(2022)Coveney, Roney, Corrado, Wilkinson, Oakley,
  Niederer and Clayton}]{coveney2022calibrating}
\bibinfo{author}{Coveney, S.}, \bibinfo{author}{Roney, C.H.},
  \bibinfo{author}{Corrado, C.}, \bibinfo{author}{Wilkinson, R.D.},
  \bibinfo{author}{Oakley, J.E.}, \bibinfo{author}{Niederer, S.A.},
  \bibinfo{author}{Clayton, R.H.}, \bibinfo{year}{2022}.
\newblock \bibinfo{title}{Calibrating cardiac electrophysiology models using
  latent gaussian processes on atrial manifolds}.
\newblock \bibinfo{journal}{Scientific Reports} \bibinfo{volume}{12},
  \bibinfo{pages}{16572}.
\bibitem[{Crozier et~al.(2016)Crozier, Augustin, Neic, Prassl, Holler, Fastl,
  Hennemuth, Bredies, Kuehne, Bishop et~al.}]{crozier2016image}
\bibinfo{author}{Crozier, A.}, \bibinfo{author}{Augustin, C.},
  \bibinfo{author}{Neic, A.}, \bibinfo{author}{Prassl, A.},
  \bibinfo{author}{Holler, M.}, \bibinfo{author}{Fastl, T.},
  \bibinfo{author}{Hennemuth, A.}, \bibinfo{author}{Bredies, K.},
  \bibinfo{author}{Kuehne, T.}, \bibinfo{author}{Bishop, M.}, et~al.,
  \bibinfo{year}{2016}.
\newblock \bibinfo{title}{Image-based personalization of cardiac anatomy for
  coupled electromechanical modeling}.
\newblock \bibinfo{journal}{Annals of biomedical engineering}
  \bibinfo{volume}{44}, \bibinfo{pages}{58--70}.
\bibitem[{Csepe et~al.(2016)Csepe, Zhao, Hansen, Li, Sul, Lim, Wang, Simonetti,
  Kilic, Mohler, Janssen and Fedorov}]{csepe_human_2016}
\bibinfo{author}{Csepe, T.A.}, \bibinfo{author}{Zhao, J.},
  \bibinfo{author}{Hansen, B.J.}, \bibinfo{author}{Li, N.},
  \bibinfo{author}{Sul, L.V.}, \bibinfo{author}{Lim, P.},
  \bibinfo{author}{Wang, Y.}, \bibinfo{author}{Simonetti, O.P.},
  \bibinfo{author}{Kilic, A.}, \bibinfo{author}{Mohler, P.J.},
  \bibinfo{author}{Janssen, P.M.}, \bibinfo{author}{Fedorov, V.V.},
  \bibinfo{year}{2016}.
\newblock \bibinfo{title}{Human sinoatrial node structure: {3D} microanatomy of
  sinoatrial conduction pathways}.
\newblock \bibinfo{journal}{Progress in Biophysics and Molecular Biology}
  \bibinfo{volume}{120}, \bibinfo{pages}{164--178}.
\newblock \URLprefix
  \url{https://linkinghub.elsevier.com/retrieve/pii/S0079610715002606},
  \DOIprefix\doi{10.1016/j.pbiomolbio.2015.12.011}.
\bibitem[{David M.~Harrild(2000)}]{harrild2000:_atria}
\bibinfo{author}{David M.~Harrild, C.S.H..}, \bibinfo{year}{2000}.
\newblock \bibinfo{title}{A {Computer} {Model} of {Normal} {Conduction} in the
  {Human} {Atria}}.
\newblock \bibinfo{journal}{Circulation Research} \bibinfo{volume}{87},
  \bibinfo{pages}{e25--e36}.
\newblock \URLprefix
  \url{https://www.ahajournals.org/doi/10.1161/01.res.87.7.e25},
  \DOIprefix\doi{10.1161/01.RES.87.7.e25}. \bibinfo{note}{publisher: American
  Heart Association}.
\bibitem[{Deng et~al.(2012)Deng, Gong, Shou, Jiao, Zhang, Ye and
  Xia}]{deng2012simulation}
\bibinfo{author}{Deng, D.d.}, \bibinfo{author}{Gong, Y.l.},
  \bibinfo{author}{Shou, G.f.}, \bibinfo{author}{Jiao, P.f.},
  \bibinfo{author}{Zhang, H.g.}, \bibinfo{author}{Ye, X.s.},
  \bibinfo{author}{Xia, L.}, \bibinfo{year}{2012}.
\newblock \bibinfo{title}{Simulation of biatrial conduction via different
  pathways during sinus rhythm with a detailed human atrial model}.
\newblock \bibinfo{journal}{Journal of Zhejiang University Science B}
  \bibinfo{volume}{13}, \bibinfo{pages}{676--694}.
\bibitem[{Dewland et~al.(2013)Dewland, Wintermark, Vaysman, Smith, Tong,
  Vittinghoff and Marcus}]{dewland2013use}
\bibinfo{author}{Dewland, T.A.}, \bibinfo{author}{Wintermark, M.},
  \bibinfo{author}{Vaysman, A.}, \bibinfo{author}{Smith, L.M.},
  \bibinfo{author}{Tong, E.}, \bibinfo{author}{Vittinghoff, E.},
  \bibinfo{author}{Marcus, G.M.}, \bibinfo{year}{2013}.
\newblock \bibinfo{title}{Use of computed tomography to identify atrial
  fibrillation associated differences in left atrial wall thickness and
  density}.
\newblock \bibinfo{journal}{Pacing and Clinical Electrophysiology}
  \bibinfo{volume}{36}, \bibinfo{pages}{55--62}.
\bibitem[{Dobrzynski et~al.(2005)Dobrzynski, Li, Tellez, Greener, Nikolski,
  Wright, Parson, Jones, Lancaster, Yamamoto et~al.}]{dobrzynski2005computer}
\bibinfo{author}{Dobrzynski, H.}, \bibinfo{author}{Li, J.},
  \bibinfo{author}{Tellez, J.}, \bibinfo{author}{Greener, I.},
  \bibinfo{author}{Nikolski, V.}, \bibinfo{author}{Wright, S.},
  \bibinfo{author}{Parson, S.}, \bibinfo{author}{Jones, S.},
  \bibinfo{author}{Lancaster, M.}, \bibinfo{author}{Yamamoto, M.}, et~al.,
  \bibinfo{year}{2005}.
\newblock \bibinfo{title}{Computer three-dimensional reconstruction of the
  sinoatrial node}.
\newblock \bibinfo{journal}{Circulation} \bibinfo{volume}{111},
  \bibinfo{pages}{846--854}.
\bibitem[{D{\"o}ssel et~al.(2012)D{\"o}ssel, Krueger, Weber, Wilhelms and
  Seemann}]{dossel2012computational}
\bibinfo{author}{D{\"o}ssel, O.}, \bibinfo{author}{Krueger, M.W.},
  \bibinfo{author}{Weber, F.M.}, \bibinfo{author}{Wilhelms, M.},
  \bibinfo{author}{Seemann, G.}, \bibinfo{year}{2012}.
\newblock \bibinfo{title}{Computational modeling of the human atrial anatomy
  and electrophysiology}.
\newblock \bibinfo{journal}{Medical \& biological engineering \& computing}
  \bibinfo{volume}{50}, \bibinfo{pages}{773--799}.
\bibitem[{Fedele et~al.(2023)Fedele, Piersanti, Regazzoni, Salvador, Africa,
  Bucelli, Zingaro, Quarteroni et~al.}]{fedele2023comprehensive}
\bibinfo{author}{Fedele, M.}, \bibinfo{author}{Piersanti, R.},
  \bibinfo{author}{Regazzoni, F.}, \bibinfo{author}{Salvador, M.},
  \bibinfo{author}{Africa, P.C.}, \bibinfo{author}{Bucelli, M.},
  \bibinfo{author}{Zingaro, A.}, \bibinfo{author}{Quarteroni, A.}, et~al.,
  \bibinfo{year}{2023}.
\newblock \bibinfo{title}{A comprehensive and biophysically detailed
  computational model of the whole human heart electromechanics}.
\newblock \bibinfo{journal}{Computer Methods in Applied Mechanics and
  Engineering} \bibinfo{volume}{410}, \bibinfo{pages}{115983}.
\bibitem[{Fedorov et~al.(2012)Fedorov, Glukhov and
  Chang}]{fedorov_conduction_2012}
\bibinfo{author}{Fedorov, V.V.}, \bibinfo{author}{Glukhov, A.V.},
  \bibinfo{author}{Chang, R.}, \bibinfo{year}{2012}.
\newblock \bibinfo{title}{Conduction barriers and pathways of the sinoatrial
  pacemaker complex: their role in normal rhythm and atrial arrhythmias}.
\newblock \bibinfo{journal}{American Journal of Physiology-Heart and
  Circulatory Physiology} \bibinfo{volume}{302}, \bibinfo{pages}{H1773--H1783}.
\newblock \URLprefix
  \url{https://www.physiology.org/doi/10.1152/ajpheart.00892.2011},
  \DOIprefix\doi{10.1152/ajpheart.00892.2011}.
\bibitem[{Fedorov et~al.(2010)Fedorov, Glukhov, Chang, Kostecki, Aferol,
  Hucker, Wuskell, Loew, Schuessler, Moazami and Efimov}]{fedorov_optical_2010}
\bibinfo{author}{Fedorov, V.V.}, \bibinfo{author}{Glukhov, A.V.},
  \bibinfo{author}{Chang, R.}, \bibinfo{author}{Kostecki, G.},
  \bibinfo{author}{Aferol, H.}, \bibinfo{author}{Hucker, W.J.},
  \bibinfo{author}{Wuskell, J.P.}, \bibinfo{author}{Loew, L.M.},
  \bibinfo{author}{Schuessler, R.B.}, \bibinfo{author}{Moazami, N.},
  \bibinfo{author}{Efimov, I.R.}, \bibinfo{year}{2010}.
\newblock \bibinfo{title}{Optical {Mapping} of the {Isolated}
  {Coronary}-{Perfused} {Human} {Sinus} {Node}}.
\newblock \bibinfo{journal}{Journal of the American College of Cardiology}
  \bibinfo{volume}{56}, \bibinfo{pages}{1386--1394}.
\newblock \URLprefix
  \url{https://linkinghub.elsevier.com/retrieve/pii/S0735109710034571},
  \DOIprefix\doi{10.1016/j.jacc.2010.03.098}.
\bibitem[{Ferrer et~al.(2015)Ferrer, Sebasti{\'a}n, S{\'a}nchez-Quintana,
  Rodriguez, Godoy, Martinez and Saiz}]{ferrer2015detailed}
\bibinfo{author}{Ferrer, A.}, \bibinfo{author}{Sebasti{\'a}n, R.},
  \bibinfo{author}{S{\'a}nchez-Quintana, D.}, \bibinfo{author}{Rodriguez,
  J.F.}, \bibinfo{author}{Godoy, E.J.}, \bibinfo{author}{Martinez, L.},
  \bibinfo{author}{Saiz, J.}, \bibinfo{year}{2015}.
\newblock \bibinfo{title}{Detailed anatomical and electrophysiological models
  of human atria and torso for the simulation of atrial activation}.
\newblock \bibinfo{journal}{PloS one} \bibinfo{volume}{10},
  \bibinfo{pages}{e0141573}.
\bibitem[{Geselowitz(1989)}]{geselowitz1989theory}
\bibinfo{author}{Geselowitz, D.B.}, \bibinfo{year}{1989}.
\newblock \bibinfo{title}{On the theory of the electrocardiogram}.
\newblock \bibinfo{journal}{Proceedings of the IEEE} \bibinfo{volume}{77},
  \bibinfo{pages}{857--876}.
\bibitem[{Gillette et~al.(2021a)Gillette, Gsell, Bouyssier, Prassl, Neic,
  Vigmond and Plank}]{gillette2021automated}
\bibinfo{author}{Gillette, K.}, \bibinfo{author}{Gsell, M.A.},
  \bibinfo{author}{Bouyssier, J.}, \bibinfo{author}{Prassl, A.J.},
  \bibinfo{author}{Neic, A.}, \bibinfo{author}{Vigmond, E.J.},
  \bibinfo{author}{Plank, G.}, \bibinfo{year}{2021}a.
\newblock \bibinfo{title}{Automated framework for the inclusion of a
  his--purkinje system in cardiac digital twins of ventricular
  electrophysiology}.
\newblock \bibinfo{journal}{Annals of biomedical engineering}
  \bibinfo{volume}{49}, \bibinfo{pages}{3143--3153}.
\bibitem[{Gillette et~al.(2022)Gillette, Gsell, Strocchi, Grandits, Neic,
  Manninger, Scherr, Roney, Prassl, Augustin et~al.}]{gillette2022personalized}
\bibinfo{author}{Gillette, K.}, \bibinfo{author}{Gsell, M.A.},
  \bibinfo{author}{Strocchi, M.}, \bibinfo{author}{Grandits, T.},
  \bibinfo{author}{Neic, A.}, \bibinfo{author}{Manninger, M.},
  \bibinfo{author}{Scherr, D.}, \bibinfo{author}{Roney, C.H.},
  \bibinfo{author}{Prassl, A.J.}, \bibinfo{author}{Augustin, C.M.}, et~al.,
  \bibinfo{year}{2022}.
\newblock \bibinfo{title}{A personalized real-time virtual model of whole heart
  electrophysiology}.
\newblock \bibinfo{journal}{Frontiers in Physiology} \bibinfo{volume}{13},
  \bibinfo{pages}{907190}.
\bibitem[{Gillette et~al.(2021b)Gillette, Gsell, Prassl, Karabelas, Reiter,
  Reiter, Grandits, Payer, Štern, Urschler, Bayer, Augustin, Neic, Pock,
  Vigmond and Plank}]{gillette2021:_framework}
\bibinfo{author}{Gillette, K.}, \bibinfo{author}{Gsell, M.A.F.},
  \bibinfo{author}{Prassl, A.J.}, \bibinfo{author}{Karabelas, E.},
  \bibinfo{author}{Reiter, U.}, \bibinfo{author}{Reiter, G.},
  \bibinfo{author}{Grandits, T.}, \bibinfo{author}{Payer, C.},
  \bibinfo{author}{Štern, D.}, \bibinfo{author}{Urschler, M.},
  \bibinfo{author}{Bayer, J.D.}, \bibinfo{author}{Augustin, C.M.},
  \bibinfo{author}{Neic, A.}, \bibinfo{author}{Pock, T.},
  \bibinfo{author}{Vigmond, E.J.}, \bibinfo{author}{Plank, G.},
  \bibinfo{year}{2021}b.
\newblock \bibinfo{title}{A {Framework} for the generation of digital twins of
  cardiac electrophysiology from clinical 12-leads {ECGs}.}
\newblock \bibinfo{journal}{Medical image analysis} \bibinfo{volume}{71},
  \bibinfo{pages}{102080}.
\newblock \URLprefix \url{https://doi.org/10.1016/j.media.2021.102080},
  \DOIprefix\doi{10.1016/j.media.2021.102080}. \bibinfo{note}{publisher:
  Elsevier B.V.}
\bibitem[{Grandits et~al.(2021)Grandits, Effland, Pock, Krause, Plank and
  Pezzuto}]{grandits2021:_geasi}
\bibinfo{author}{Grandits, T.}, \bibinfo{author}{Effland, A.},
  \bibinfo{author}{Pock, T.}, \bibinfo{author}{Krause, R.},
  \bibinfo{author}{Plank, G.}, \bibinfo{author}{Pezzuto, S.},
  \bibinfo{year}{2021}.
\newblock \bibinfo{title}{{GEASI}: {Geodesic}-based earliest activation sites
  identification in cardiac models.}
\newblock \bibinfo{journal}{International journal for numerical methods in
  biomedical engineering} \bibinfo{volume}{37}, \bibinfo{pages}{e3505}.
\newblock \URLprefix \url{http://arxiv.org/abs/2102.09962},
  \DOIprefix\doi{10.1002/cnm.3505}. \bibinfo{note}{arXiv: 2102.09962}.
\bibitem[{Grandits et~al.(2024a)Grandits, Gillette, Plank and
  Pezzuto}]{grandits2024:_identify}
\bibinfo{author}{Grandits, T.}, \bibinfo{author}{Gillette, K.},
  \bibinfo{author}{Plank, G.}, \bibinfo{author}{Pezzuto, S.},
  \bibinfo{year}{2024}a.
\newblock \bibinfo{title}{Accurate and {Efficient} {Cardiac} {Digital} {Twin}
  from surface {ECGs}: {Insights} into {Identifiability} of {Ventricular}
  {Conduction} {System}}.
\newblock \URLprefix \url{http://arxiv.org/abs/2411.00165},
  \DOIprefix\doi{10.48550/arXiv.2411.00165}. \bibinfo{note}{arXiv:2411.00165}.
\bibitem[{Grandits et~al.(2024b)Grandits, Verhülsdonk, Haase, Effland and
  Pezzuto}]{grandits2023:_geodesic_bp}
\bibinfo{author}{Grandits, T.}, \bibinfo{author}{Verhülsdonk, J.},
  \bibinfo{author}{Haase, G.}, \bibinfo{author}{Effland, A.},
  \bibinfo{author}{Pezzuto, S.}, \bibinfo{year}{2024}b.
\newblock \bibinfo{title}{Digital {Twinning} of {Cardiac} {Electrophysiology}
  {Models} {From} the {Surface} {ECG}: {A} {Geodesic} {Backpropagation}
  {Approach}}.
\newblock \bibinfo{journal}{IEEE Transactions on Biomedical Engineering}
  \bibinfo{volume}{71}, \bibinfo{pages}{1281--1288}.
\newblock \URLprefix \url{https://ieeexplore.ieee.org/document/10339854},
  \DOIprefix\doi{10.1109/TBME.2023.3331876}. \bibinfo{note}{conference Name:
  IEEE Transactions on Biomedical Engineering}.
\bibitem[{Gray et~al.(1996)Gray, Pertsov and Jalife}]{gray1996incomplete}
\bibinfo{author}{Gray, R.A.}, \bibinfo{author}{Pertsov, A.M.},
  \bibinfo{author}{Jalife, J.}, \bibinfo{year}{1996}.
\newblock \bibinfo{title}{Incomplete reentry and epicardial breakthrough
  patterns during atrial fibrillation in the sheep heart}.
\newblock \bibinfo{journal}{Circulation} \bibinfo{volume}{94},
  \bibinfo{pages}{2649--2661}.
\bibitem[{Gsell et~al.(2024)Gsell, Neic, Bishop, Gillette, Prassl, Augustin,
  Vigmond and Plank}]{gsell2024:_forcepss}
\bibinfo{author}{Gsell, M.A.}, \bibinfo{author}{Neic, A.},
  \bibinfo{author}{Bishop, M.J.}, \bibinfo{author}{Gillette, K.},
  \bibinfo{author}{Prassl, A.J.}, \bibinfo{author}{Augustin, C.M.},
  \bibinfo{author}{Vigmond, E.J.}, \bibinfo{author}{Plank, G.},
  \bibinfo{year}{2024}.
\newblock \bibinfo{title}{{ForCEPSS}—{A} framework for cardiac
  electrophysiology simulations standardization}.
\newblock \bibinfo{journal}{Computer Methods and Programs in Biomedicine}
  \bibinfo{volume}{251}, \bibinfo{pages}{108189}.
\newblock \URLprefix
  \url{https://linkinghub.elsevier.com/retrieve/pii/S0169260724001858},
  \DOIprefix\doi{10.1016/j.cmpb.2024.108189}.
\bibitem[{Hansson et~al.(1998)Hansson, Holm, Blomstr{\"o}m, Johansson,
  L{\"u}hrs, Brandt and Olsson}]{hansson1998right}
\bibinfo{author}{Hansson, A.}, \bibinfo{author}{Holm, M.},
  \bibinfo{author}{Blomstr{\"o}m, P.}, \bibinfo{author}{Johansson, R.},
  \bibinfo{author}{L{\"u}hrs, C.}, \bibinfo{author}{Brandt, J.},
  \bibinfo{author}{Olsson, S.}, \bibinfo{year}{1998}.
\newblock \bibinfo{title}{Right atrial free wall conduction velocity and degree
  of anisotropy in patients with stable sinus rhythm studied during open heart
  surgery}.
\newblock \bibinfo{journal}{European heart journal} \bibinfo{volume}{19},
  \bibinfo{pages}{293--300}.
\bibitem[{Harrild and Henriquez(2000)}]{harrild2000computer}
\bibinfo{author}{Harrild, D.M.}, \bibinfo{author}{Henriquez, C.S.},
  \bibinfo{year}{2000}.
\newblock \bibinfo{title}{A computer model of normal conduction in the human
  atria}.
\newblock \bibinfo{journal}{Circulation research} \bibinfo{volume}{87},
  \bibinfo{pages}{e25--e36}.
\bibitem[{Ho and S{\'a}nchez-Quintana(2009)}]{ho2009importance}
\bibinfo{author}{Ho, S.}, \bibinfo{author}{S{\'a}nchez-Quintana, D.},
  \bibinfo{year}{2009}.
\newblock \bibinfo{title}{The importance of atrial structure and fibers}.
\newblock \bibinfo{journal}{Clinical Anatomy: The Official Journal of the
  American Association of Clinical Anatomists and the British Association of
  Clinical Anatomists} \bibinfo{volume}{22}, \bibinfo{pages}{52--63}.
\bibitem[{Honjo et~al.(1996)Honjo, Boyett, Kodama and
  Toyama}]{honjo_correlation_1996}
\bibinfo{author}{Honjo, H.}, \bibinfo{author}{Boyett, M.},
  \bibinfo{author}{Kodama, I.}, \bibinfo{author}{Toyama, J.},
  \bibinfo{year}{1996}.
\newblock \bibinfo{title}{Correlation between electrical activity and the size
  of rabbit sino-atrial node cells.}
\newblock \bibinfo{journal}{The Journal of physiology} \bibinfo{volume}{496},
  \bibinfo{pages}{795--808}.
\bibitem[{Hopman et~al.(2023)Hopman, Visch, Bhagirath, van~der Laan, Mulder,
  Razeghi, Kemme, Niederer, Allaart and G{\"o}tte}]{hopman2023right}
\bibinfo{author}{Hopman, L.H.}, \bibinfo{author}{Visch, J.E.},
  \bibinfo{author}{Bhagirath, P.}, \bibinfo{author}{van~der Laan, A.M.},
  \bibinfo{author}{Mulder, M.J.}, \bibinfo{author}{Razeghi, O.},
  \bibinfo{author}{Kemme, M.J.}, \bibinfo{author}{Niederer, S.A.},
  \bibinfo{author}{Allaart, C.P.}, \bibinfo{author}{G{\"o}tte, M.J.},
  \bibinfo{year}{2023}.
\newblock \bibinfo{title}{Right atrial function and fibrosis in relation to
  successful atrial fibrillation ablation}.
\newblock \bibinfo{journal}{European Heart Journal-Cardiovascular Imaging}
  \bibinfo{volume}{24}, \bibinfo{pages}{336--345}.
\bibitem[{Joyner and Van~Capelle(1986)}]{joyner_propagation_1986}
\bibinfo{author}{Joyner, R.}, \bibinfo{author}{Van~Capelle, F.},
  \bibinfo{year}{1986}.
\newblock \bibinfo{title}{Propagation through electrically coupled cells. {How}
  a small {SA} node drives a large atrium}.
\newblock \bibinfo{journal}{Biophysical Journal} \bibinfo{volume}{50},
  \bibinfo{pages}{1157--1164}.
\newblock \URLprefix
  \url{https://linkinghub.elsevier.com/retrieve/pii/S0006349586835597},
  \DOIprefix\doi{10.1016/S0006-3495(86)83559-7}.
\bibitem[{Kanchi and Masud(2007)}]{kanchi20073d}
\bibinfo{author}{Kanchi, H.}, \bibinfo{author}{Masud, A.},
  \bibinfo{year}{2007}.
\newblock \bibinfo{title}{A 3d adaptive mesh moving scheme}.
\newblock \bibinfo{journal}{International Journal for Numerical Methods in
  Fluids} \bibinfo{volume}{54}, \bibinfo{pages}{923--944}.
\bibitem[{Karabelas et~al.(2018)Karabelas, Gsell, Augustin, Marx, Neic, Prassl,
  Goubergrits, Kuehne and Plank}]{karabelas2018:_towards}
\bibinfo{author}{Karabelas, E.}, \bibinfo{author}{Gsell, M.A.F.},
  \bibinfo{author}{Augustin, C.M.}, \bibinfo{author}{Marx, L.},
  \bibinfo{author}{Neic, A.}, \bibinfo{author}{Prassl, A.J.},
  \bibinfo{author}{Goubergrits, L.}, \bibinfo{author}{Kuehne, T.},
  \bibinfo{author}{Plank, G.}, \bibinfo{year}{2018}.
\newblock \bibinfo{title}{Towards a {Computational} {Framework} for {Modeling}
  the {Impact} of {Aortic} {Coarctations} {Upon} {Left} {Ventricular} {Load}.}
\newblock \bibinfo{journal}{Frontiers in physiology} \bibinfo{volume}{9},
  \bibinfo{pages}{538}.
\newblock \URLprefix \url{http://www.ncbi.nlm.nih.gov/pubmed/29892227},
  \DOIprefix\doi{10.3389/fphys.2018.00538}.
\bibitem[{Keith and Flack(1907)}]{keith1907form}
\bibinfo{author}{Keith, A.}, \bibinfo{author}{Flack, M.}, \bibinfo{year}{1907}.
\newblock \bibinfo{title}{The form and nature of the muscular connections
  between the primary divisions of the vertebrate heart}.
\newblock \bibinfo{journal}{Journal of anatomy and physiology}
  \bibinfo{volume}{41}, \bibinfo{pages}{172}.
\bibitem[{Keller et~al.(2010)Keller, Weber, Seemann and
  Dossel}]{keller2010ranking}
\bibinfo{author}{Keller, D.U.}, \bibinfo{author}{Weber, F.M.},
  \bibinfo{author}{Seemann, G.}, \bibinfo{author}{Dossel, O.},
  \bibinfo{year}{2010}.
\newblock \bibinfo{title}{Ranking the influence of tissue conductivities on
  forward-calculated ecgs}.
\newblock \bibinfo{journal}{IEEE Transactions on Biomedical Engineering}
  \bibinfo{volume}{57}, \bibinfo{pages}{1568--1576}.
\bibitem[{Kharbanda et~al.(2019)Kharbanda, {\"O}zdemir, Taverne, Kik, Bogers
  and de~Groot}]{kharbanda2019current}
\bibinfo{author}{Kharbanda, R.K.}, \bibinfo{author}{{\"O}zdemir, E.H.},
  \bibinfo{author}{Taverne, Y.J.}, \bibinfo{author}{Kik, C.},
  \bibinfo{author}{Bogers, A.J.}, \bibinfo{author}{de~Groot, N.M.},
  \bibinfo{year}{2019}.
\newblock \bibinfo{title}{Current concepts of anatomy, electrophysiology, and
  therapeutic implications of the interatrial septum}.
\newblock \bibinfo{journal}{JACC: Clinical Electrophysiology}
  \bibinfo{volume}{5}, \bibinfo{pages}{647--656}.
\bibitem[{Knol et~al.(2019)Knol, Teuwen, Kleinrensink, Bogers, de~Groot and
  Taverne}]{knol2019bachmann}
\bibinfo{author}{Knol, W.G.}, \bibinfo{author}{Teuwen, C.P.},
  \bibinfo{author}{Kleinrensink, G.J.}, \bibinfo{author}{Bogers, A.J.},
  \bibinfo{author}{de~Groot, N.M.}, \bibinfo{author}{Taverne, Y.J.},
  \bibinfo{year}{2019}.
\newblock \bibinfo{title}{The bachmann bundle and interatrial conduction:
  comparing atrial morphology to electrical activity}.
\newblock \bibinfo{journal}{Heart Rhythm} \bibinfo{volume}{16},
  \bibinfo{pages}{606--614}.
\bibitem[{Knupp(2022)}]{knupp2022worst}
\bibinfo{author}{Knupp, P.}, \bibinfo{year}{2022}.
\newblock \bibinfo{title}{Worst case mesh quality in the target matrix
  optimization paradigm}.
\newblock \bibinfo{journal}{Engineering with Computers} \bibinfo{volume}{38},
  \bibinfo{pages}{5695--5711}.
\bibitem[{Krueger et~al.(2012)Krueger, Seemann, Rhode, Keller, Schilling,
  Arujuna, Gill, O'Neill, Razavi and Dossel}]{krueger2012personalization}
\bibinfo{author}{Krueger, M.W.}, \bibinfo{author}{Seemann, G.},
  \bibinfo{author}{Rhode, K.}, \bibinfo{author}{Keller, D.U.},
  \bibinfo{author}{Schilling, C.}, \bibinfo{author}{Arujuna, A.},
  \bibinfo{author}{Gill, J.}, \bibinfo{author}{O'Neill, M.D.},
  \bibinfo{author}{Razavi, R.}, \bibinfo{author}{Dossel, O.},
  \bibinfo{year}{2012}.
\newblock \bibinfo{title}{Personalization of atrial anatomy and
  electrophysiology as a basis for clinical modeling of radio-frequency
  ablation of atrial fibrillation}.
\newblock \bibinfo{journal}{IEEE transactions on medical imaging}
  \bibinfo{volume}{32}, \bibinfo{pages}{73--84}.
\bibitem[{Labarthe et~al.(2014)Labarthe, Bayer, Coudiere, Henry, Cochet, Jais
  and Vigmond}]{labarthe2014:_bilayer}
\bibinfo{author}{Labarthe, S.}, \bibinfo{author}{Bayer, J.},
  \bibinfo{author}{Coudiere, Y.}, \bibinfo{author}{Henry, J.},
  \bibinfo{author}{Cochet, H.}, \bibinfo{author}{Jais, P.},
  \bibinfo{author}{Vigmond, E.}, \bibinfo{year}{2014}.
\newblock \bibinfo{title}{A bilayer model of human atria: mathematical
  background, construction, and assessment}.
\newblock \bibinfo{journal}{Europace} \bibinfo{volume}{16},
  \bibinfo{pages}{iv21--iv29}.
\newblock \URLprefix
  \url{https://academic.oup.com/europace/article-lookup/doi/10.1093/europace/euu256},
  \DOIprefix\doi{10.1093/europace/euu256}.
\bibitem[{Lang et~al.(2022)Lang, Cameli, Sade, Faletra, Fortuni, Rossi and
  Soulat-Dufour}]{lang2022imaging}
\bibinfo{author}{Lang, R.M.}, \bibinfo{author}{Cameli, M.},
  \bibinfo{author}{Sade, L.E.}, \bibinfo{author}{Faletra, F.F.},
  \bibinfo{author}{Fortuni, F.}, \bibinfo{author}{Rossi, A.},
  \bibinfo{author}{Soulat-Dufour, L.}, \bibinfo{year}{2022}.
\newblock \bibinfo{title}{Imaging assessment of the right atrium: anatomy and
  function}.
\newblock \bibinfo{journal}{European Heart Journal-Cardiovascular Imaging}
  \bibinfo{volume}{23}, \bibinfo{pages}{867--884}.
\bibitem[{Lemery et~al.(2007)Lemery, Birnie, Tang, Green, Gollob, Hendry and
  Lau}]{lemery2007normal}
\bibinfo{author}{Lemery, R.}, \bibinfo{author}{Birnie, D.},
  \bibinfo{author}{Tang, A.S.}, \bibinfo{author}{Green, M.},
  \bibinfo{author}{Gollob, M.}, \bibinfo{author}{Hendry, M.},
  \bibinfo{author}{Lau, E.}, \bibinfo{year}{2007}.
\newblock \bibinfo{title}{Normal atrial activation and voltage during sinus
  rhythm in the human heart: an endocardial and epicardial mapping study in
  patients with a history of atrial fibrillation}.
\newblock \bibinfo{journal}{Journal of cardiovascular electrophysiology}
  \bibinfo{volume}{18}, \bibinfo{pages}{402--408}.
\bibitem[{Lemery et~al.(2003)Lemery, Guiraudon and Veinot}]{lemery2003anatomic}
\bibinfo{author}{Lemery, R.}, \bibinfo{author}{Guiraudon, G.},
  \bibinfo{author}{Veinot, J.P.}, \bibinfo{year}{2003}.
\newblock \bibinfo{title}{Anatomic description of bachmann's bundle and its
  relation to the atrial septum}.
\newblock \bibinfo{journal}{The American journal of cardiology}
  \bibinfo{volume}{91}, \bibinfo{pages}{1482--1482}.
\bibitem[{Li et~al.(2017)Li, Hansen, Csepe, Zhao, Ignozzi, Sul, Zakharkin,
  Kalyanasundaram, Davis, Biesiadecki, Kilic, Janssen, Mohler, Weiss, Hummel
  and Fedorov}]{li_redundant_2017}
\bibinfo{author}{Li, N.}, \bibinfo{author}{Hansen, B.J.},
  \bibinfo{author}{Csepe, T.A.}, \bibinfo{author}{Zhao, J.},
  \bibinfo{author}{Ignozzi, A.J.}, \bibinfo{author}{Sul, L.V.},
  \bibinfo{author}{Zakharkin, S.O.}, \bibinfo{author}{Kalyanasundaram, A.},
  \bibinfo{author}{Davis, J.P.}, \bibinfo{author}{Biesiadecki, B.J.},
  \bibinfo{author}{Kilic, A.}, \bibinfo{author}{Janssen, P.M.L.},
  \bibinfo{author}{Mohler, P.J.}, \bibinfo{author}{Weiss, R.},
  \bibinfo{author}{Hummel, J.D.}, \bibinfo{author}{Fedorov, V.V.},
  \bibinfo{year}{2017}.
\newblock \bibinfo{title}{Redundant and diverse intranodal pacemakers and
  conduction pathways protect the human sinoatrial node from failure}.
\newblock \bibinfo{journal}{Science Translational Medicine}
  \bibinfo{volume}{9}, \bibinfo{pages}{eaam5607}.
\newblock \URLprefix
  \url{https://www.science.org/doi/10.1126/scitranslmed.aam5607},
  \DOIprefix\doi{10.1126/scitranslmed.aam5607}.
\bibitem[{Loewe et~al.(2016)Loewe, Krueger, Holmqvist, D{\"o}ssel, Seemann and
  Platonov}]{loewe2016influence}
\bibinfo{author}{Loewe, A.}, \bibinfo{author}{Krueger, M.W.},
  \bibinfo{author}{Holmqvist, F.}, \bibinfo{author}{D{\"o}ssel, O.},
  \bibinfo{author}{Seemann, G.}, \bibinfo{author}{Platonov, P.G.},
  \bibinfo{year}{2016}.
\newblock \bibinfo{title}{Influence of the earliest right atrial activation
  site and its proximity to interatrial connections on p-wave morphology}.
\newblock \bibinfo{journal}{EP Europace} \bibinfo{volume}{18},
  \bibinfo{pages}{iv35--iv43}.
\bibitem[{Loewe et~al.(2015)Loewe, Krueger, Platonov, Holmqvist, D{\"o}ssel and
  Seemann}]{loewe2015left}
\bibinfo{author}{Loewe, A.}, \bibinfo{author}{Krueger, M.W.},
  \bibinfo{author}{Platonov, P.G.}, \bibinfo{author}{Holmqvist, F.},
  \bibinfo{author}{D{\"o}ssel, O.}, \bibinfo{author}{Seemann, G.},
  \bibinfo{year}{2015}.
\newblock \bibinfo{title}{Left and right atrial contribution to the p-wave in
  realistic computational models}, in: \bibinfo{booktitle}{Functional Imaging
  and Modeling of the Heart: 8th International Conference, FIMH 2015,
  Maastricht, The Netherlands, June 25-27, 2015. Proceedings 8},
  \bibinfo{organization}{Springer}. pp. \bibinfo{pages}{439--447}.
\bibitem[{Lubrecht et~al.(2021)Lubrecht, Grandits, Gharaviri, Schotten, Pock,
  Plank, Krause, Auricchio, Conte and Pezzuto}]{lubrecht2021:_piemap}
\bibinfo{author}{Lubrecht, J.M.}, \bibinfo{author}{Grandits, T.},
  \bibinfo{author}{Gharaviri, A.}, \bibinfo{author}{Schotten, U.},
  \bibinfo{author}{Pock, T.}, \bibinfo{author}{Plank, G.},
  \bibinfo{author}{Krause, R.}, \bibinfo{author}{Auricchio, A.},
  \bibinfo{author}{Conte, G.}, \bibinfo{author}{Pezzuto, S.},
  \bibinfo{year}{2021}.
\newblock \bibinfo{title}{Automatic reconstruction of the left atrium
  activation from sparse intracardiac contact recordings by inverse estimate of
  fibre structure and anisotropic conduction in a patient-specific model}.
\newblock \bibinfo{journal}{EP Europace} \bibinfo{volume}{23},
  \bibinfo{pages}{i63--i70}.
\newblock \URLprefix
  \url{https://academic.oup.com/europace/article/23/Supplement_1/i63/6158558},
  \DOIprefix\doi{10.1093/europace/euaa392}.
\bibitem[{Matsuyama et~al.(2004)Matsuyama, Inoue, Kobayashi, Sakai, Saito,
  Katagiri and Ota}]{matsuyama2004anatomical}
\bibinfo{author}{Matsuyama, T.a.}, \bibinfo{author}{Inoue, S.},
  \bibinfo{author}{Kobayashi, Y.}, \bibinfo{author}{Sakai, T.},
  \bibinfo{author}{Saito, T.}, \bibinfo{author}{Katagiri, T.},
  \bibinfo{author}{Ota, H.}, \bibinfo{year}{2004}.
\newblock \bibinfo{title}{Anatomical diversity and age-related histological
  changes in the human right atrial posterolateral wall}.
\newblock \bibinfo{journal}{EP Europace} \bibinfo{volume}{6},
  \bibinfo{pages}{307--315}.
\bibitem[{Monfredi et~al.(2010)Monfredi, Dobrzynski, Mondal, Boyett and
  Morris}]{monfredi2010anatomy}
\bibinfo{author}{Monfredi, O.}, \bibinfo{author}{Dobrzynski, H.},
  \bibinfo{author}{Mondal, T.}, \bibinfo{author}{Boyett, M.R.},
  \bibinfo{author}{Morris, G.M.}, \bibinfo{year}{2010}.
\newblock \bibinfo{title}{The anatomy and physiology of the sinoatrial node—a
  contemporary review}.
\newblock \bibinfo{journal}{Pacing and clinical electrophysiology}
  \bibinfo{volume}{33}, \bibinfo{pages}{1392--1406}.
\bibitem[{Multerer and Pezzuto(2021)}]{multerer2021uncertainty}
\bibinfo{author}{Multerer, M.}, \bibinfo{author}{Pezzuto, S.},
  \bibinfo{year}{2021}.
\newblock \bibinfo{title}{Uncertainty quantification for the 12-lead ecg: a
  lead field approach}.
\newblock \bibinfo{journal}{arXiv preprint arXiv:2102.09960} .
\bibitem[{Muñoz et~al.(2011)Muñoz, Kaur and Vigmond}]{munoz_onset_2011}
\bibinfo{author}{Muñoz, M.A.}, \bibinfo{author}{Kaur, J.},
  \bibinfo{author}{Vigmond, E.J.}, \bibinfo{year}{2011}.
\newblock \bibinfo{title}{Onset of atrial arrhythmias elicited by autonomic
  modulation of rabbit sinoatrial node activity: a modeling study}.
\newblock \bibinfo{journal}{American Journal of Physiology-Heart and
  Circulatory Physiology} \bibinfo{volume}{301}, \bibinfo{pages}{H1974--H1983}.
\newblock \DOIprefix\doi{10.1152/ajpheart.00059.2011}.
\bibitem[{Nagel et~al.(2022)Nagel, Espinosa, Gillette, Gsell, S{\'a}nchez,
  Plank, D{\"o}ssel and Loewe}]{nagel2022comparison}
\bibinfo{author}{Nagel, C.}, \bibinfo{author}{Espinosa, C.B.},
  \bibinfo{author}{Gillette, K.}, \bibinfo{author}{Gsell, M.A.},
  \bibinfo{author}{S{\'a}nchez, J.}, \bibinfo{author}{Plank, G.},
  \bibinfo{author}{D{\"o}ssel, O.}, \bibinfo{author}{Loewe, A.},
  \bibinfo{year}{2022}.
\newblock \bibinfo{title}{Comparison of propagation models and forward
  calculation methods on cellular, tissue and organ scale atrial
  electrophysiology}.
\newblock \bibinfo{journal}{IEEE Transactions on Biomedical Engineering}
  \bibinfo{volume}{70}, \bibinfo{pages}{511--522}.
\bibitem[{Nagel et~al.(2021)Nagel, Schuler, Dössel and
  Loewe}]{nagel2021:_biatrial_ssm}
\bibinfo{author}{Nagel, C.}, \bibinfo{author}{Schuler, S.},
  \bibinfo{author}{Dössel, O.}, \bibinfo{author}{Loewe, A.},
  \bibinfo{year}{2021}.
\newblock \bibinfo{title}{A bi-atrial statistical shape model for large-scale
  in silico studies of human atria: {Model} development and application to
  {ECG} simulations}.
\newblock \bibinfo{journal}{Medical Image Analysis} \bibinfo{volume}{74},
  \bibinfo{pages}{102210}.
\newblock \URLprefix
  \url{https://linkinghub.elsevier.com/retrieve/pii/S1361841521002553},
  \DOIprefix\doi{10.1016/j.media.2021.102210}.
\bibitem[{Neic et~al.(2017a)Neic, Campos, Prassl, Neiderer, Bishop, Vigmond and
  Plank}]{neic2017eikonal}
\bibinfo{author}{Neic, A.}, \bibinfo{author}{Campos, F.O.},
  \bibinfo{author}{Prassl, A.J.}, \bibinfo{author}{Neiderer, S.A.},
  \bibinfo{author}{Bishop, M.J.}, \bibinfo{author}{Vigmond, E.J.},
  \bibinfo{author}{Plank, G.}, \bibinfo{year}{2017}a.
\newblock \bibinfo{title}{Efficient computation of electrograms and ecgs in
  human whole heart simulations using a reaction-eikonal model}.
\bibitem[{Neic et~al.(2017b)Neic, Campos, Prassl, Niederer, Bishop, Vigmond and
  Plank}]{neic17:_reaction_eikonal}
\bibinfo{author}{Neic, A.}, \bibinfo{author}{Campos, F.O.},
  \bibinfo{author}{Prassl, A.J.}, \bibinfo{author}{Niederer, S.A.},
  \bibinfo{author}{Bishop, M.J.}, \bibinfo{author}{Vigmond, E.J.},
  \bibinfo{author}{Plank, G.}, \bibinfo{year}{2017}b.
\newblock \bibinfo{title}{Efficient computation of electrograms and ecgs in
  human whole heart simulations using a reaction-eikonal model.}
\newblock \bibinfo{journal}{Journal of computational physics}
  \bibinfo{volume}{346}, \bibinfo{pages}{191--211}.
\newblock \DOIprefix\doi{10.1016/j.jcp.2017.06.020}.
\bibitem[{Neic et~al.(under review)Neic, Gsell, Karabelas, Prassl and
  Plank}]{neic19:_meshtool}
\bibinfo{author}{Neic, A.}, \bibinfo{author}{Gsell, M.},
  \bibinfo{author}{Karabelas, E.}, \bibinfo{author}{Prassl, A.},
  \bibinfo{author}{Plank, G.}, \bibinfo{year}{under review}.
\newblock \bibinfo{title}{Automating image-based mesh generation and
  manipulation tasks in cardiac modeling workflows using meshtool}.
\newblock \bibinfo{journal}{Software X} .
\bibitem[{Neic et~al.(2020)Neic, Gsell, Karabelas, Prassl and
  Plank}]{neic2020automating}
\bibinfo{author}{Neic, A.}, \bibinfo{author}{Gsell, M.A.},
  \bibinfo{author}{Karabelas, E.}, \bibinfo{author}{Prassl, A.J.},
  \bibinfo{author}{Plank, G.}, \bibinfo{year}{2020}.
\newblock \bibinfo{title}{Automating image-based mesh generation and
  manipulation tasks in cardiac modeling workflows using meshtool}.
\newblock \bibinfo{journal}{SoftwareX} \bibinfo{volume}{11},
  \bibinfo{pages}{100454}.
\bibitem[{Niederer et~al.(2020)Niederer, Aboelkassem, Cantwell, Corrado,
  Coveney, Cherry, Delhaas, Fenton, Panfilov, Pathmanathan, Plank, Riabiz,
  Roney, {Dos Santos} and Wang}]{niederer2020creation}
\bibinfo{author}{Niederer, S.A.}, \bibinfo{author}{Aboelkassem, Y.},
  \bibinfo{author}{Cantwell, C.D.}, \bibinfo{author}{Corrado, C.},
  \bibinfo{author}{Coveney, S.}, \bibinfo{author}{Cherry, E.M.},
  \bibinfo{author}{Delhaas, T.}, \bibinfo{author}{Fenton, F.H.},
  \bibinfo{author}{Panfilov, A.V.}, \bibinfo{author}{Pathmanathan, P.},
  \bibinfo{author}{Plank, G.}, \bibinfo{author}{Riabiz, M.},
  \bibinfo{author}{Roney, C.H.}, \bibinfo{author}{{Dos Santos}, R.W.},
  \bibinfo{author}{Wang, L.}, \bibinfo{year}{2020}.
\newblock \bibinfo{title}{{Creation and application of virtual patient cohorts
  of heart models.}}
\newblock \bibinfo{journal}{Philos. Trans. A. Math. Phys. Eng. Sci.}
  \bibinfo{volume}{378}, \bibinfo{pages}{20190558}.
\newblock \DOIprefix\doi{10.1098/rsta.2019.0558}.
\bibitem[{Padala et~al.(2021)Padala, Cabrera and
  Ellenbogen}]{padala2021anatomy}
\bibinfo{author}{Padala, S.K.}, \bibinfo{author}{Cabrera, J.A.},
  \bibinfo{author}{Ellenbogen, K.A.}, \bibinfo{year}{2021}.
\newblock \bibinfo{title}{Anatomy of the cardiac conduction system}.
\newblock \bibinfo{journal}{Pacing and Clinical Electrophysiology}
  \bibinfo{volume}{44}, \bibinfo{pages}{15--25}.
\bibitem[{Payer et~al.(2017)Payer, {\v{S}}tern, Bischof and
  Urschler}]{payer2017multi}
\bibinfo{author}{Payer, C.}, \bibinfo{author}{{\v{S}}tern, D.},
  \bibinfo{author}{Bischof, H.}, \bibinfo{author}{Urschler, M.},
  \bibinfo{year}{2017}.
\newblock \bibinfo{title}{Multi-label whole heart segmentation using cnns and
  anatomical label configurations}, in: \bibinfo{booktitle}{International
  Workshop on Statistical Atlases and Computational Models of the Heart},
  \bibinfo{organization}{Springer}. pp. \bibinfo{pages}{190--198}.
\bibitem[{Pezzuto et~al.(2017)Pezzuto, Kal'avsk{\'{y}}, Potse, Prinzen,
  Auricchio and Krause}]{pezzuto2017:_ecg}
\bibinfo{author}{Pezzuto, S.}, \bibinfo{author}{Kal'avsk{\'{y}}, P.},
  \bibinfo{author}{Potse, M.}, \bibinfo{author}{Prinzen, F.W.},
  \bibinfo{author}{Auricchio, A.}, \bibinfo{author}{Krause, R.},
  \bibinfo{year}{2017}.
\newblock \bibinfo{title}{{Evaluation of a rapid anisotropic model for {ECG}
  simulation}}.
\newblock \bibinfo{journal}{Front. Physiol.} \bibinfo{volume}{8}.
\newblock \DOIprefix\doi{10.3389/fphys.2017.00265}.
\bibitem[{Piersanti et~al.(2021)Piersanti, Africa, Fedele, Vergara, Dedè,
  Corno and Quarteroni}]{Piersanti2021}
\bibinfo{author}{Piersanti, R.}, \bibinfo{author}{Africa, P.C.},
  \bibinfo{author}{Fedele, M.}, \bibinfo{author}{Vergara, C.},
  \bibinfo{author}{Dedè, L.}, \bibinfo{author}{Corno, A.F.},
  \bibinfo{author}{Quarteroni, A.}, \bibinfo{year}{2021}.
\newblock \bibinfo{title}{Modeling cardiac muscle fibers in ventricular and
  atrial electrophysiology simulations}.
\newblock \bibinfo{journal}{Computer Methods in Applied Mechanics and
  Engineering} \bibinfo{volume}{373}, \bibinfo{pages}{113468}.
\newblock \URLprefix \url{https://doi.org/10.1016/j.cma.2020.113468},
  \DOIprefix\doi{10.1016/j.cma.2020.113468}. \bibinfo{note}{publisher: Elsevier
  B.V.}
\bibitem[{Plank et~al.(2021)Plank, Loewe, Neic, Augustin, Huang, Gsell,
  Karabelas, Nothstein, Prassl, S{\'a}nchez et~al.}]{plank2021opencarp}
\bibinfo{author}{Plank, G.}, \bibinfo{author}{Loewe, A.},
  \bibinfo{author}{Neic, A.}, \bibinfo{author}{Augustin, C.},
  \bibinfo{author}{Huang, Y.L.}, \bibinfo{author}{Gsell, M.A.},
  \bibinfo{author}{Karabelas, E.}, \bibinfo{author}{Nothstein, M.},
  \bibinfo{author}{Prassl, A.J.}, \bibinfo{author}{S{\'a}nchez, J.}, et~al.,
  \bibinfo{year}{2021}.
\newblock \bibinfo{title}{The opencarp simulation environment for cardiac
  electrophysiology}.
\newblock \bibinfo{journal}{Computer methods and Programs in Biomedicine}
  \bibinfo{volume}{208}, \bibinfo{pages}{106223}.
\bibitem[{Platonov et~al.(2002)Platonov, Mitrofanova, Chireikin and
  Olsson}]{platonov2002morphology}
\bibinfo{author}{Platonov, P.}, \bibinfo{author}{Mitrofanova, L.},
  \bibinfo{author}{Chireikin, L.}, \bibinfo{author}{Olsson, S.},
  \bibinfo{year}{2002}.
\newblock \bibinfo{title}{Morphology of inter-atrial conduction routes in
  patients with atrial fibrillation}.
\newblock \bibinfo{journal}{Europace} \bibinfo{volume}{4},
  \bibinfo{pages}{183--192}.
\bibitem[{Platonov(2007)}]{platonov2007interatrial}
\bibinfo{author}{Platonov, P.G.}, \bibinfo{year}{2007}.
\newblock \bibinfo{title}{Interatrial conduction in the mechanisms of atrial
  fibrillation: from anatomy to cardiac signals and new treatment modalities}.
\newblock \bibinfo{journal}{Europace} \bibinfo{volume}{9},
  \bibinfo{pages}{vi10--vi16}.
\bibitem[{Potse(2018)}]{potse2018scalable}
\bibinfo{author}{Potse, M.}, \bibinfo{year}{2018}.
\newblock \bibinfo{title}{Scalable and accurate {ECG} simulation for
  reaction-diffusion models of the human heart}.
\newblock \bibinfo{journal}{Frontiers in physiology} \bibinfo{volume}{9},
  \bibinfo{pages}{370}.
\bibitem[{Potse et~al.(2006)Potse, Dub{\'e}, Richer, Vinet and
  Gulrajani}]{potse2006comparison}
\bibinfo{author}{Potse, M.}, \bibinfo{author}{Dub{\'e}, B.},
  \bibinfo{author}{Richer, J.}, \bibinfo{author}{Vinet, A.},
  \bibinfo{author}{Gulrajani, R.M.}, \bibinfo{year}{2006}.
\newblock \bibinfo{title}{A comparison of monodomain and bidomain
  reaction-diffusion models for action potential propagation in the human
  heart}.
\newblock \bibinfo{journal}{IEEE Transactions on Biomedical Engineering}
  \bibinfo{volume}{53}, \bibinfo{pages}{2425--2435}.
\bibitem[{Prassl et~al.(2009)Prassl, Kickinger, Ahammer, Grau, Schneider,
  Hofer, Vigmond, Trayanova and Plank}]{prassl09:_tarantula}
\bibinfo{author}{Prassl, A.J.}, \bibinfo{author}{Kickinger, F.},
  \bibinfo{author}{Ahammer, H.}, \bibinfo{author}{Grau, V.},
  \bibinfo{author}{Schneider, J.E.}, \bibinfo{author}{Hofer, E.},
  \bibinfo{author}{Vigmond, E.J.}, \bibinfo{author}{Trayanova, N.A.},
  \bibinfo{author}{Plank, G.}, \bibinfo{year}{2009}.
\newblock \bibinfo{title}{Automatically generated, anatomically accurate meshes
  for cardiac electrophysiology problems.}
\newblock \bibinfo{journal}{IEEE transactions on bio-medical engineering}
  \bibinfo{volume}{56}, \bibinfo{pages}{1318--1330}.
\newblock \DOIprefix\doi{10.1109/TBME.2009.2014243}.
\bibitem[{Reddy et~al.(2023)Reddy, Kong, Petru, Maan, Funasako, Minami,
  Ruppersberg, Dukkipati and Neuzil}]{reddy2023electrographic}
\bibinfo{author}{Reddy, V.Y.}, \bibinfo{author}{Kong, M.H.},
  \bibinfo{author}{Petru, J.}, \bibinfo{author}{Maan, A.},
  \bibinfo{author}{Funasako, M.}, \bibinfo{author}{Minami, K.},
  \bibinfo{author}{Ruppersberg, P.}, \bibinfo{author}{Dukkipati, S.},
  \bibinfo{author}{Neuzil, P.}, \bibinfo{year}{2023}.
\newblock \bibinfo{title}{Electrographic flow mapping of persistent atrial
  fibrillation: intra-and inter-procedure reproducibility in the absence of
  ‘ground truth’}.
\newblock \bibinfo{journal}{Europace} \bibinfo{volume}{25},
  \bibinfo{pages}{euad308}.
\bibitem[{Roberts et~al.(1979)Roberts, Hersh and Scher}]{roberts1979influence}
\bibinfo{author}{Roberts, D.E.}, \bibinfo{author}{Hersh, L.T.},
  \bibinfo{author}{Scher, A.M.}, \bibinfo{year}{1979}.
\newblock \bibinfo{title}{Influence of cardiac fiber orientation on wavefront
  voltage, conduction velocity, and tissue resistivity in the dog.}
\newblock \bibinfo{journal}{Circulation research} \bibinfo{volume}{44},
  \bibinfo{pages}{701--712}.
\bibitem[{Rodero et~al.(2021)Rodero, Strocchi, Marciniak, Longobardi, Whitaker,
  O’Neill, Gillette, Augustin, Plank, Vigmond, Lamata and
  Niederer}]{rodero2021linking}
\bibinfo{author}{Rodero, C.}, \bibinfo{author}{Strocchi, M.},
  \bibinfo{author}{Marciniak, M.}, \bibinfo{author}{Longobardi, S.},
  \bibinfo{author}{Whitaker, J.}, \bibinfo{author}{O’Neill, M.D.},
  \bibinfo{author}{Gillette, K.}, \bibinfo{author}{Augustin, C.},
  \bibinfo{author}{Plank, G.}, \bibinfo{author}{Vigmond, E.J.},
  \bibinfo{author}{Lamata, P.}, \bibinfo{author}{Niederer, S.A.},
  \bibinfo{year}{2021}.
\newblock \bibinfo{title}{Linking statistical shape models and simulated
  function in the healthy adult human heart}.
\newblock \bibinfo{journal}{PLOS Computational Biology} \bibinfo{volume}{17},
  \bibinfo{pages}{e1008851}.
\newblock \DOIprefix\doi{10.1371/journal.pcbi.1008851}.
\bibitem[{Roney et~al.(2023)Roney, Solis~Lemus and et~al}]{roney2023:_bia_vol}
\bibinfo{author}{Roney, C.}, \bibinfo{author}{Solis~Lemus, J.A.},
  \bibinfo{author}{et~al}, \bibinfo{year}{2023}.
\newblock \bibinfo{title}{Constructing bilayer and volumetric atrial models at
  scale}.
\newblock \bibinfo{journal}{Interface Focus} , \bibinfo{pages}{20230038}.
\bibitem[{Roney et~al.(2021)Roney, Bendikas, Pashakhanloo, Corrado, Vigmond,
  McVeigh, Trayanova and Niederer}]{roney2021constructing}
\bibinfo{author}{Roney, C.H.}, \bibinfo{author}{Bendikas, R.},
  \bibinfo{author}{Pashakhanloo, F.}, \bibinfo{author}{Corrado, C.},
  \bibinfo{author}{Vigmond, E.J.}, \bibinfo{author}{McVeigh, E.R.},
  \bibinfo{author}{Trayanova, N.A.}, \bibinfo{author}{Niederer, S.A.},
  \bibinfo{year}{2021}.
\newblock \bibinfo{title}{Constructing a human atrial fibre atlas}.
\newblock \bibinfo{journal}{Annals of biomedical engineering}
  \bibinfo{volume}{49}, \bibinfo{pages}{233--250}.
\bibitem[{Roney et~al.(2019a)Roney, Pashaei, Meo, Dubois, Boyle, Trayanova,
  Cochet, Niederer and Vigmond}]{roney2019uac}
\bibinfo{author}{Roney, C.H.}, \bibinfo{author}{Pashaei, A.},
  \bibinfo{author}{Meo, M.}, \bibinfo{author}{Dubois, R.},
  \bibinfo{author}{Boyle, P.M.}, \bibinfo{author}{Trayanova, N.A.},
  \bibinfo{author}{Cochet, H.}, \bibinfo{author}{Niederer, S.A.},
  \bibinfo{author}{Vigmond, E.J.}, \bibinfo{year}{2019}a.
\newblock \bibinfo{title}{{Universal atrial coordinates applied to
  visualisation, registration and construction of patient specific meshes}}.
\newblock \bibinfo{journal}{Med. Image Anal.} \bibinfo{volume}{55},
  \bibinfo{pages}{65--75}.
\newblock \DOIprefix\doi{10.1016/j.media.2019.04.004},
  \href{http://arxiv.org/abs/1810.06630}{\tt arXiv:1810.06630}.
\bibitem[{Roney et~al.(2019b)Roney, Pashaei, Meo, Dubois, Boyle, Trayanova,
  Cochet, Niederer and Vigmond}]{RONEY201965}
\bibinfo{author}{Roney, C.H.}, \bibinfo{author}{Pashaei, A.},
  \bibinfo{author}{Meo, M.}, \bibinfo{author}{Dubois, R.},
  \bibinfo{author}{Boyle, P.M.}, \bibinfo{author}{Trayanova, N.A.},
  \bibinfo{author}{Cochet, H.}, \bibinfo{author}{Niederer, S.A.},
  \bibinfo{author}{Vigmond, E.J.}, \bibinfo{year}{2019}b.
\newblock \bibinfo{title}{Universal atrial coordinates applied to
  visualisation, registration and construction of patient specific meshes}.
\newblock \bibinfo{journal}{Medical Image Analysis} \bibinfo{volume}{55},
  \bibinfo{pages}{65--75}.
\newblock \URLprefix
  \url{https://www.sciencedirect.com/science/article/pii/S1361841518308089},
  \DOIprefix\doi{https://doi.org/10.1016/j.media.2019.04.004}.
\bibitem[{Sakamoto et~al.(2005)Sakamoto, Nitta, Ishii, Miyagi, Ohmori and
  Shimizu}]{sakamoto2005interatrial}
\bibinfo{author}{Sakamoto, S.I.}, \bibinfo{author}{Nitta, T.},
  \bibinfo{author}{Ishii, Y.}, \bibinfo{author}{Miyagi, Y.},
  \bibinfo{author}{Ohmori, H.}, \bibinfo{author}{Shimizu, K.},
  \bibinfo{year}{2005}.
\newblock \bibinfo{title}{Interatrial electrical connections: the precise
  location and preferential conduction}.
\newblock \bibinfo{journal}{Journal of cardiovascular electrophysiology}
  \bibinfo{volume}{16}, \bibinfo{pages}{1077--1086}.
\bibitem[{Sakata et~al.(2024)Sakata, Bradley, Prakosa, Yamamoto, Ali, Loeffler,
  Tice, Boyle, Kholmovski, Yadav, Sinha, Marine, Calkins, Spragg and
  Trayanova}]{sakata2024:_assessing}
\bibinfo{author}{Sakata, K.}, \bibinfo{author}{Bradley, R.P.},
  \bibinfo{author}{Prakosa, A.}, \bibinfo{author}{Yamamoto, C.A.P.},
  \bibinfo{author}{Ali, S.Y.}, \bibinfo{author}{Loeffler, S.},
  \bibinfo{author}{Tice, B.M.}, \bibinfo{author}{Boyle, P.M.},
  \bibinfo{author}{Kholmovski, E.G.}, \bibinfo{author}{Yadav, R.},
  \bibinfo{author}{Sinha, S.K.}, \bibinfo{author}{Marine, J.E.},
  \bibinfo{author}{Calkins, H.}, \bibinfo{author}{Spragg, D.D.},
  \bibinfo{author}{Trayanova, N.A.}, \bibinfo{year}{2024}.
\newblock \bibinfo{title}{Assessing the arrhythmogenic propensity of fibrotic
  substrate using digital twins to inform a mechanisms-based atrial
  fibrillation ablation strategy}.
\newblock \bibinfo{journal}{Nature Cardiovascular Research}
  \bibinfo{volume}{3}, \bibinfo{pages}{857--868}.
\newblock \URLprefix \url{https://www.nature.com/articles/s44161-024-00489-x},
  \DOIprefix\doi{10.1038/s44161-024-00489-x}. \bibinfo{note}{publisher: Nature
  Publishing Group}.
\bibitem[{Schneider et~al.(2018)Schneider, Hu, Dumas, Gao, Panozzo and
  Zorin}]{schneider2018decoupling}
\bibinfo{author}{Schneider, T.}, \bibinfo{author}{Hu, Y.},
  \bibinfo{author}{Dumas, J.}, \bibinfo{author}{Gao, X.},
  \bibinfo{author}{Panozzo, D.}, \bibinfo{author}{Zorin, D.},
  \bibinfo{year}{2018}.
\newblock \bibinfo{title}{Decoupling simulation accuracy from mesh quality}.
\newblock \bibinfo{journal}{ACM transactions on graphics} .
\bibitem[{Schotten et~al.(2011)Schotten, Verheule, Kirchhof and
  Goette}]{schotten2011pathophysiological}
\bibinfo{author}{Schotten, U.}, \bibinfo{author}{Verheule, S.},
  \bibinfo{author}{Kirchhof, P.}, \bibinfo{author}{Goette, A.},
  \bibinfo{year}{2011}.
\newblock \bibinfo{title}{Pathophysiological mechanisms of atrial fibrillation:
  a translational appraisal}.
\newblock \bibinfo{journal}{Physiological reviews} \bibinfo{volume}{91},
  \bibinfo{pages}{265--325}.
\bibitem[{Schuler et~al.(2021)Schuler, Pilia, Potyagaylo and
  Loewe}]{schuler2021cobiveco}
\bibinfo{author}{Schuler, S.}, \bibinfo{author}{Pilia, N.},
  \bibinfo{author}{Potyagaylo, D.}, \bibinfo{author}{Loewe, A.},
  \bibinfo{year}{2021}.
\newblock \bibinfo{title}{Cobiveco: Consistent biventricular coordinates for
  precise and intuitive description of position in the heart--with matlab
  implementation}.
\newblock \bibinfo{journal}{Medical Image Analysis} \bibinfo{volume}{74},
  \bibinfo{pages}{102247}.
\bibitem[{Schuler et~al.(2019)Schuler, Tate, Oostendorp, MacLeod and
  Dössel}]{schuler2019:_spatial}
\bibinfo{author}{Schuler, S.}, \bibinfo{author}{Tate, J.D.},
  \bibinfo{author}{Oostendorp, T.F.}, \bibinfo{author}{MacLeod, R.S.},
  \bibinfo{author}{Dössel, O.}, \bibinfo{year}{2019}.
\newblock \bibinfo{title}{Spatial {Downsampling} of {Surface} {Sources} in the
  {Forward} {Problem} of {Electrocardiography}}, in:
  \bibinfo{editor}{Coudière, Y.}, \bibinfo{editor}{Ozenne, V.},
  \bibinfo{editor}{Vigmond, E.}, \bibinfo{editor}{Zemzemi, N.} (Eds.),
  \bibinfo{booktitle}{Functional {Imaging} and {Modeling} of the {Heart}},
  \bibinfo{publisher}{Springer International Publishing},
  \bibinfo{address}{Cham}. pp. \bibinfo{pages}{29--36}.
\newblock \DOIprefix\doi{10.1007/978-3-030-21949-9_4}.
\bibitem[{Seemann et~al.(2006)Seemann, Höper, Sachse, Dössel, Holden and
  Zhang}]{seemann2006:_atria}
\bibinfo{author}{Seemann, G.}, \bibinfo{author}{Höper, C.},
  \bibinfo{author}{Sachse, F.B.}, \bibinfo{author}{Dössel, O.},
  \bibinfo{author}{Holden, A.V.}, \bibinfo{author}{Zhang, H.},
  \bibinfo{year}{2006}.
\newblock \bibinfo{title}{Heterogeneous three-dimensional anatomical and
  electrophysiological model of human atria.}
\newblock \bibinfo{journal}{Philosophical transactions. Series A, Mathematical,
  physical, and engineering sciences} \bibinfo{volume}{364},
  \bibinfo{pages}{1465--81}.
\newblock \URLprefix \url{http://www.ncbi.nlm.nih.gov/pubmed/16766355},
  \DOIprefix\doi{10.1098/rsta.2006.1781}.
\bibitem[{Szili-Torok et~al.(2023)Szili-Torok, Neuzil, Langbein, Petru,
  Funasako, Dinshaw, Wijchers, Bhagwandien, Rillig, Spitzer
  et~al.}]{szili2023electrographic}
\bibinfo{author}{Szili-Torok, T.}, \bibinfo{author}{Neuzil, P.},
  \bibinfo{author}{Langbein, A.}, \bibinfo{author}{Petru, J.},
  \bibinfo{author}{Funasako, M.}, \bibinfo{author}{Dinshaw, L.},
  \bibinfo{author}{Wijchers, S.}, \bibinfo{author}{Bhagwandien, R.},
  \bibinfo{author}{Rillig, A.}, \bibinfo{author}{Spitzer, S.G.}, et~al.,
  \bibinfo{year}{2023}.
\newblock \bibinfo{title}{Electrographic flow--guided ablation in redo patients
  with persistent atrial fibrillation (flow-af): design and rationale}.
\newblock \bibinfo{journal}{Heart Rhythm O2} \bibinfo{volume}{4},
  \bibinfo{pages}{391--400}.
\bibitem[{Tan et~al.(2000)Tan, Chan and Choi}]{tan2000detection}
\bibinfo{author}{Tan, K.}, \bibinfo{author}{Chan, K.}, \bibinfo{author}{Choi,
  K.}, \bibinfo{year}{2000}.
\newblock \bibinfo{title}{Detection of the qrs complex, p wave and t wave in
  electrocardiogram}, in: \bibinfo{booktitle}{2000 First International
  Conference Advances in Medical Signal and Information Processing (IEE Conf.
  Publ. No. 476)}, \bibinfo{organization}{IET}. pp. \bibinfo{pages}{41--47}.
\bibitem[{Thaler et~al.(2021)Thaler, Payer, Bischof and
  Stern}]{thaler2021efficient}
\bibinfo{author}{Thaler, F.}, \bibinfo{author}{Payer, C.},
  \bibinfo{author}{Bischof, H.}, \bibinfo{author}{Stern, D.},
  \bibinfo{year}{2021}.
\newblock \bibinfo{title}{Efficient multi-organ segmentation using spatial
  configuration-net with low gpu memory requirements}.
\newblock \bibinfo{journal}{arXiv preprint arXiv:2111.13630} .
\bibitem[{Varela et~al.(2017)Varela, Morgan, Theron, Dillon-Murphy, Chubb,
  Whitaker, Henningsson, Aljabar, Schaeffter, Kolbitsch
  et~al.}]{varela2017novel}
\bibinfo{author}{Varela, M.}, \bibinfo{author}{Morgan, R.},
  \bibinfo{author}{Theron, A.}, \bibinfo{author}{Dillon-Murphy, D.},
  \bibinfo{author}{Chubb, H.}, \bibinfo{author}{Whitaker, J.},
  \bibinfo{author}{Henningsson, M.}, \bibinfo{author}{Aljabar, P.},
  \bibinfo{author}{Schaeffter, T.}, \bibinfo{author}{Kolbitsch, C.}, et~al.,
  \bibinfo{year}{2017}.
\newblock \bibinfo{title}{Novel mri technique enables non-invasive measurement
  of atrial wall thickness}.
\newblock \bibinfo{journal}{IEEE transactions on medical imaging}
  \bibinfo{volume}{36}, \bibinfo{pages}{1607--1614}.
\bibitem[{Viceconti et~al.(2020)Viceconti, Pappalardo, Rodriguez, Horner,
  Bischoff and Musuamba~Tshinanu}]{viceconti2020:_in_silico}
\bibinfo{author}{Viceconti, M.}, \bibinfo{author}{Pappalardo, F.},
  \bibinfo{author}{Rodriguez, B.}, \bibinfo{author}{Horner, M.},
  \bibinfo{author}{Bischoff, J.}, \bibinfo{author}{Musuamba~Tshinanu, F.},
  \bibinfo{year}{2020}.
\newblock \bibinfo{title}{In silico trials: {Verification}, validation and
  uncertainty quantification of predictive models used in the regulatory
  evaluation of biomedical products}.
\newblock \bibinfo{journal}{Methods}
  \DOIprefix\doi{10.1016/j.ymeth.2020.01.011}.
\bibitem[{Vigmond and Clements(2007)}]{vigmond2006:_sawtooth}
\bibinfo{author}{Vigmond, E.J.}, \bibinfo{author}{Clements, C.},
  \bibinfo{year}{2007}.
\newblock \bibinfo{title}{Construction of a computer model to investigate
  sawtooth effects in the {Purkinje} system}.
\newblock \bibinfo{journal}{IEEE Transactions on Biomedical Engineering}
  \bibinfo{volume}{54}, \bibinfo{pages}{389--399}.
\newblock \DOIprefix\doi{10.1109/TBME.2006.888817}. \bibinfo{note}{iSBN:
  0018-9294 (Print){\textbackslash}r0018-9294 (Linking)}.
\bibitem[{Wachter et~al.(2015)Wachter, Loewe, Krueger, D{\"o}ssel and
  Seemann}]{wachter2015mesh}
\bibinfo{author}{Wachter, A.}, \bibinfo{author}{Loewe, A.},
  \bibinfo{author}{Krueger, M.W.}, \bibinfo{author}{D{\"o}ssel, O.},
  \bibinfo{author}{Seemann, G.}, \bibinfo{year}{2015}.
\newblock \bibinfo{title}{Mesh structure-independent modeling of
  patient-specific atrial fiber orientation}.
\newblock \bibinfo{journal}{Current Directions in Biomedical Engineering}
  \bibinfo{volume}{1}, \bibinfo{pages}{409--412}.
\bibitem[{Whitaker et~al.(2016)Whitaker, Rajani, Chubb, Gabrawi, Varela,
  Wright, Niederer and O'Neill}]{whitaker2016role}
\bibinfo{author}{Whitaker, J.}, \bibinfo{author}{Rajani, R.},
  \bibinfo{author}{Chubb, H.}, \bibinfo{author}{Gabrawi, M.},
  \bibinfo{author}{Varela, M.}, \bibinfo{author}{Wright, M.},
  \bibinfo{author}{Niederer, S.}, \bibinfo{author}{O'Neill, M.D.},
  \bibinfo{year}{2016}.
\newblock \bibinfo{title}{The role of myocardial wall thickness in atrial
  arrhythmogenesis}.
\newblock \bibinfo{journal}{Ep Europace} \bibinfo{volume}{18},
  \bibinfo{pages}{1758--1772}.
\bibitem[{Whittaker et~al.(2020)Whittaker, Clerx, Lei, Christini and
  Mirams}]{whittaker2020calibration}
\bibinfo{author}{Whittaker, D.G.}, \bibinfo{author}{Clerx, M.},
  \bibinfo{author}{Lei, C.L.}, \bibinfo{author}{Christini, D.J.},
  \bibinfo{author}{Mirams, G.R.}, \bibinfo{year}{2020}.
\newblock \bibinfo{title}{Calibration of ionic and cellular cardiac
  electrophysiology models}.
\newblock \bibinfo{journal}{Wiley Interdisciplinary Reviews: Systems Biology
  and Medicine} \bibinfo{volume}{12}, \bibinfo{pages}{e1482}.
\bibitem[{Yang et~al.(2015)Yang, Li, Campbell and Jia}]{yang2015go}
\bibinfo{author}{Yang, J.}, \bibinfo{author}{Li, H.},
  \bibinfo{author}{Campbell, D.}, \bibinfo{author}{Jia, Y.},
  \bibinfo{year}{2015}.
\newblock \bibinfo{title}{Go-icp: A globally optimal solution to 3d icp
  point-set registration}.
\newblock \bibinfo{journal}{IEEE transactions on pattern analysis and machine
  intelligence} \bibinfo{volume}{38}, \bibinfo{pages}{2241--2254}.
\bibitem[{Yang et~al.(2013)Yang, Li and Jia}]{yang2013go}
\bibinfo{author}{Yang, J.}, \bibinfo{author}{Li, H.}, \bibinfo{author}{Jia,
  Y.}, \bibinfo{year}{2013}.
\newblock \bibinfo{title}{Go-icp: Solving 3d registration efficiently and
  globally optimally}, in: \bibinfo{booktitle}{Proceedings of the IEEE
  International Conference on Computer Vision}, pp.
  \bibinfo{pages}{1457--1464}.
\bibitem[{Yushkevich et~al.(2006)Yushkevich, Piven, Cody~Hazlett, Gimpel~Smith,
  Ho, Gee and Gerig}]{py06nimg}
\bibinfo{author}{Yushkevich, P.A.}, \bibinfo{author}{Piven, J.},
  \bibinfo{author}{Cody~Hazlett, H.}, \bibinfo{author}{Gimpel~Smith, R.},
  \bibinfo{author}{Ho, S.}, \bibinfo{author}{Gee, J.C.},
  \bibinfo{author}{Gerig, G.}, \bibinfo{year}{2006}.
\newblock \bibinfo{title}{User-guided {3D} active contour segmentation of
  anatomical structures: Significantly improved efficiency and reliability}.
\newblock \bibinfo{journal}{Neuroimage} \bibinfo{volume}{31},
  \bibinfo{pages}{1116--1128}.
\bibitem[{Zappon et~al.(2024)Zappon, Manzoni and Quarteroni}]{ZAPPON2024112815}
\bibinfo{author}{Zappon, E.}, \bibinfo{author}{Manzoni, A.},
  \bibinfo{author}{Quarteroni, A.}, \bibinfo{year}{2024}.
\newblock \bibinfo{title}{A non-conforming-in-space numerical framework for
  realistic cardiac electrophysiological outputs}.
\newblock \bibinfo{journal}{Journal of Computational Physics}
  \bibinfo{volume}{502}, \bibinfo{pages}{112815}.
\newblock \URLprefix
  \url{https://www.sciencedirect.com/science/article/pii/S0021999124000640},
  \DOIprefix\doi{https://doi.org/10.1016/j.jcp.2024.112815}.
\bibitem[{Zheng et~al.(2021)Zheng, Azzolin, S\'{a}nchez, D\"{o}ssel and
  Loewe}]{ZhengAzzolinSanchezDosselLoewe2021}
\bibinfo{author}{Zheng, T.}, \bibinfo{author}{Azzolin, L.},
  \bibinfo{author}{S\'{a}nchez, J.}, \bibinfo{author}{D\"{o}ssel, O.},
  \bibinfo{author}{Loewe, A.}, \bibinfo{year}{2021}.
\newblock \bibinfo{title}{An automate pipeline for generating fiber orientation
  and region annotation in patient specific atrial models}.
\newblock \bibinfo{journal}{Current Directions in Biomedical Engineering}
  \bibinfo{volume}{7}, \bibinfo{pages}{136--139}.
\newblock \URLprefix \url{https://doi.org/10.1515/cdbme-2021-2035},
  \DOIprefix\doi{doi:10.1515/cdbme-2021-2035}.

\end{thebibliography}
	\FloatBarrier

\end{document}